\documentclass[a4paper,reqno,10pt,oneside]{amsart}
\usepackage[left=2.61cm,right=2.61cm,top=2.72cm,bottom=2.72cm]{geometry}
\usepackage{hyperref}
\usepackage{color}
\usepackage[english]{babel}
\usepackage[toc,page,header]{appendix}
\usepackage{mathrsfs,amssymb,amsmath, amsfonts, amsthm,enumitem,latexsym}

\usepackage{ esint }

\allowdisplaybreaks[1]

\numberwithin{equation}{section}

\DeclareFontFamily{U}{mathb}{\hyphenchar\font45}
\DeclareFontShape{U}{mathb}{m}{n}{ <-6> matha5 <6-7> matha6 <7-8>
mathb7 <8-9> mathb8 <9-10> mathb9 <10-12> mathb10 <12-> mathb12 }{}
\DeclareSymbolFont{mathb}{U}{mathb}{m}{n}

\DeclareMathAccent{\abxring}{0}{mathb}{"38}

\DeclareFontFamily{U}{mathb}{\hyphenchar\font45}
\DeclareFontShape{U}{mathb}{m}{n}{ <-6> matha5 <6-7> matha6 <7-8>
mathb7 <8-9> mathb8 <9-10> mathb9 <10-12> mathb10 <12-> mathb12 }{}
\DeclareSymbolFont{mathb}{U}{mathb}{m}{n}

\renewcommand{\O}{\operatorname{O}}

\renewcommand{\(}{\left(}
\renewcommand{\)}{\right)}
\renewcommand{\[}{\left[}
\renewcommand{\]}{\right]}

\newtheorem{theorem}{Theorem}[section]
\newtheorem{proposition}[theorem]{Proposition}
\newtheorem{corollary}[theorem]{Corollary}
\newtheorem{lemma}[theorem]{Lemma}
\theoremstyle{definition}
\newtheorem{remark}[theorem]{Remark}
\theoremstyle{definition}
\newtheorem{definition}[theorem]{Definition}
\theoremstyle{definition}

\renewcommand{\le}{\leqslant}
\renewcommand{\ge}{\geqslant}

\renewcommand{\O}{\mathcal{O}}

\newcommand{\A}{{\mathcal A}}

\newcommand{\N}{\mathbb{N}}

\newcommand{\E}{\mathcal{E}}

\newcommand{\V}{\mathcal{V}}

\newcommand{\eps}{\varepsilon}

\newcommand{\beq}{\begin{equation}}
\newcommand{\eeq}{\end{equation}}
\newcommand{\beqs}{\begin{equation*}}
\newcommand{\eeqs}{\end{equation*}}
\newcommand{\beqn}{\begin{eqnarray}}
\newcommand{\eeqn}{\end{eqnarray}}
\newcommand{\beqns}{\begin{eqnarray*}}
\newcommand{\eeqns}{\end{eqnarray*}}
\newcommand{\bdoc}{\begin{document}}
\newcommand{\edoc}{\end{document}}
\newcommand{\be}{\begin{enumerate}}
\newcommand{\ee}{\end{enumerate}}
\newcommand{\bdescr}{\begin{description}}
\newcommand{\edescr}{\end{description}}
\newcommand{\ba}{\begin{array}}
\newcommand{\ea}{\end{array}}
\newcommand{\intR}{\int_{\mathbb R^N}}
\newcommand{\R}{\mathbb R}
\newcommand{\RN}{\mathbb{R}^N}
\newcommand{\B}{\mathbb B}
\newcommand{\C}{\mathcal C}
\renewcommand{\H}{\mathcal H}
\renewcommand{\L}{\mathbb L}
\newcommand{\parallelsum}{\mathbin{\!/\mkern-5mu/\!}}
\newcommand{\e}{\varepsilon}
\newcommand{\Co}{\C_\omega}
 \renewcommand{\(}{\left(}
\renewcommand{\)}{\right)}
\renewcommand{\[}{\left[}
\renewcommand{\]}{\right]}
\renewcommand{\appendixpagename}{\centering Appendix}

\begin{document}
\title[A shape optimization problem in cylinders and related overdetermined problems]{A shape optimization problem in cylinders\\and related overdetermined problems}

\author{P\MakeLowercase{aolo} Caldiroli, A\MakeLowercase{lessandro} Iacopetti, F\MakeLowercase{ilomena} Pacella}

\thanks{\emph{Acknowledgements.} Research partially supported by Gruppo Nazionale per l'Analisi Matematica, la Pro\-ba\-bi\-li\-t\`a e le loro Applicazioni (GNAMPA) of the Istituto Nazionale di Alta Matematica (INdAM). Paolo Caldiroli is partially supported by the PRIN 2022 PNRR project 20227HX33Z \emph{Pattern formation in nonlinear phenomena}, founded by the European Union - Next Generation EU. Alessandro Iacopetti is partially supported by the PRIN 2022 PNRR project 2022R537CS \emph{$NO^3$ - Nodal Optimization, NOnlinear elliptic equations, NOnlocal geometric problems, with a focus on regularity},  founded by the European Union - Next Generation EU and by the GNAMPA 2024 project E53C23001670001: ``Aspetti geometrici e analitici di alcuni problemi locali e non-locali in mancanza di compattezza". Filomena Pacella is partially supported by the PRIN 2022 PNRR project 2022AKNSE4
\emph{Variational and Analytical aspects of Geometric PDEs
}, founded by the European Union - Next Generation EU}

\address[Paolo Caldiroli]{Dipartimento di Matematica ``Giuseppe Peano", Universit\`a degli Studi di Torino, Via Carlo Alberto 10, 10123 Torino, Italy}
\email{paolo.caldiroli@unito.it}

\address[Alessandro Iacopetti]{Dipartimento di Matematica ``Giuseppe Peano", Universit\`a degli Studi di Torino, Via Carlo Alberto 10, 10123 Torino, Italy}
\email{alessandro.iacopetti@unito.it}

\address[Filomena Pacella]{Dipartimento di Matematica, Sapienza Universit\`a di Roma,  P.le Aldo Moro 2, 00185 Roma, Italy}
\email{filomena.pacella@mat.uniroma1.it}

\begin{abstract}
In this paper, we study a shape optimization problem for the torsional energy associated with a domain contained in an infinite cylinder, under a volume constraint. We prove that a minimizer exists for all fixed volumes and show some of its geometric and topological properties. As this issue is closely related to the question of characterizing domains in cylinders that admit solutions to an overdetermined problem, our minimization result allows us to deduce interesting consequences in that direction. In particular, we find that, for some cylinders and some volumes, the ``trivial" domain given by a bounded cylinder is not the only domain where the overdetermined problem has a solution. Moreover, it is not even a minimizer, which indicates that solutions with flat level sets are not always the best candidates for optimizing the torsional energy.
\end{abstract}

\maketitle

\section{Introduction}
Let $N\in\N$, $N\geq 2$, let $\omega$ be a bounded domain in $\R^{N-1}$ with Lipschitz boundary and let $\C_\omega$ be the cylinder spanned by $\omega$, namely
\begin{equation*}%\label{def:cylspanned}
\C_\omega:=\{x=(x^\prime,x_N) \in \R^{N}; \ x^\prime\in\omega,\, x_N\in\R\}=\omega\times\R.
\end{equation*}
For a bounded domain $\Omega\subset \C_\omega$ we set:
$$
\Gamma_\Omega:=\partial \Omega \cap \C_\omega, \quad \Gamma_{1,\Omega}:=\partial\Omega\cap\partial\C_\omega.
$$
Usually $\Gamma_\Omega$ is called the relative (to $\C_\omega$) or free boundary of $\Omega$.
We consider the torsion problem with mixed boundary conditions:
\beq
\label{eq:torsionmixedbc}
\begin{cases}
-\Delta u = 1 & \text{in $\Omega$},\\
 u = 0 & \text{on $\Gamma_\Omega$},\\
\frac{\partial u}{\partial \nu} = 0 & \text{on $\Gamma_{1,\Omega}$}.\\
\end{cases}
\eeq
It is easy to see that \eqref{eq:torsionmixedbc} has a unique weak solution $u_\Omega$ in the space $H_0^1(\Omega\cup\Gamma_{1,\Omega})$ which is the subspace of functions in $H^1(\Omega)$ whose trace vanishes on $\Gamma_\Omega$ (see \cite[Chap. 1]{DP}). Indeed, $u_\Omega$ is the unique minimizer of the energy functional
$$%\label{eq:defenergyJ}
J_\Omega(v):= \frac{1}{2}\int_\Omega |\nabla v|^2 \ dx - \int_\Omega v \ dx, \ \ v\in H_0^1(\Omega\cup\Gamma_{1,\Omega}),
$$
and $u_\Omega>0$ in $\Omega$, by the maximum principle. 

Usually, the function $u_\Omega$ is called {\em torsion function} or {\em energy function}  of $\Omega$ and its energy $J_\Omega(u_\Omega)$ represents the {\em torsional energy} of the domain $\Omega$. This allows us to consider the functional
$$%\label{eq:deftorsionalenergy}
\E(\Omega; \C_\omega):=J_\Omega(u_\Omega),
$$
for every bounded domain $\Omega\subset \C_\omega$, as described above.

From the weak formulation of \eqref{eq:torsionmixedbc} we have
$$\int_\Omega |\nabla u_\Omega|^2\, dx=\int_\Omega u_\Omega\, dx,$$
so that
$$
%\label{eq:rewritetorsionalenergy}
\E(\Omega; \C_\omega)=-\frac{1}{2} \int_\Omega |\nabla u_\Omega|^2\, dx= -\frac{1}{2} \int_\Omega u_\Omega\, dx.
$$

In the present paper we address the question of minimizing the torsional energy functional $\E$ among all domains inside $\C_\omega$ with a fixed measure $c>0$.

This is a classical shape-optimization problem reframed in the case of domains $\Omega$ which are constrained to be inside the cylinder $\C_\omega$. 

Our main result shows that, whatever cylinder we consider, a minimizer for $\E(\Omega ; \C_\omega)$ always exists, for any prescribed volume.

To prove it we need to formulate our minimization problem in the larger class of quasi-open sets contained in $\C_\omega$ (see Sect. \ref{S:preliminaries} and \ref{S:sequences}). This can be done since the existence of a unique weak solution to \eqref{eq:torsionmixedbc} also holds for quasi-open sets in a suitable Sobolev space (see Lemma \ref{energy-function}).
Then, for any fixed $c>0$ we define
\begin{gather*}
%\begin{split}
\A_{\omega,c}:= \{\Omega\subset \C_\omega;\ \Omega \ \hbox{quasi-open} \ \hbox{and} \ |\Omega|\leq c\}, 
%\label{eq:defclassadmissibleintro}
\\[3pt]
\mathcal{O}_c(\C_\omega):=\inf\{\E(\Omega;\C_\omega); \  \Omega \in \A_{\omega; c}\},
%\label{def:inftorsionvolumeconstrintro}
%\end{split}
\end{gather*}
and we have:

\begin{theorem}\label{mainteo1}
Let $\omega$ be any bounded domain in $\R^{N-1}$ with Lipschitz boundary and let $c>0$. Then there exists a quasi-open set $\Omega^*\in \A_{\omega,c}$ which attains $\O_c(\C_\omega)$.
Moreover $\Omega^*$ is a bounded, connected open set and $|\Omega^*|=c$. 
\end{theorem}
We can also show that if $c>0$ is small enough or if $c$ is large enough and $\Omega^{*}$ satisfies some additional assumptions, then $\Omega^{*}$ touches the boundary of the container on a set of positive $(N-1)$-dimensional Hausdorff measure, that is $\mathcal{H}^{N-1}(\Gamma_{1,\Omega^{*}})>0$ (see Corollary \ref{cor1:sect6} (ii) and (iii)). 
 
Let us observe that, since $\C_\omega$ is an unbounded set, the existence of a minimizer is not obvious. To prove it we use the concentration-compactness principle (see \cite{LIO}) as done in a similar optimization problem in cones \cite{IPW} (see also \cite{BU}).

However in \cite{IPW} the proof strongly exploited the scaling invariance of the problem in cones, in particular to exclude the loss of compactness by dichotomy and also to show the connectedness of minimizers. Here we do not have this property, but we exploit other arguments, especially we make use of Steiner symmetrization which is appropriate for domains in cylinders. This allows to show that one minimizer can always be found in the class of domains which are convex in the $x_N$-direction and symmetric with respect to the hyperplane $$T_0:=\{x=(x^\prime,x_N)\in \R^N; \ x_N=0\}.$$ We also point out that, in fact, every minimizer of $\E(\cdot;\C_{\omega})$ whose relative boundary does not contain vertical parts turns out to be convex in the $x_N$-direction and symmetric with respect to some horizontal hyperplane, in view of a result by Cianchi and Fusco \cite{CF} (see Remark \ref{R:CF}). 
\begin{remark}
If instead of taking the cylinder $\C_\omega$ we consider the half-cylinder spanned by $\omega$: 
$$\C^+_\omega:=\{x=(x^\prime,x_N) \in \R^{N}; \ x^\prime\in\omega,\, x_N>0\}=\omega\times\left]0,+\infty\right[,$$
we would get an equivalent minimization problem. This is a consequence of the symmetry of a minimizer in $\C_\omega$ with respect to the hyperplane $T_0$. We refer to Sect. \ref{S:half} for details.
\end{remark}

Concerning the regularity of the relative boundary of minimizers of $\Omega\mapsto \E(\Omega; \C_\omega)$ with $\Omega\in \A_{\omega,c}$ we have the following:

\begin{theorem}\label{teoreg}
There exists a critical dimension $d^*$ which can be either $5,6$ or $7$, such that, for any minimizer $\Omega^*$ of $\E(\cdot; \C_\omega)$ in $\A_{\omega,c}$ , it holds that:
 \begin{itemize}
\item[$(i)$] if $N<d^*$ then $\Gamma_{\Omega^*}$ is locally a smooth manifold;
\item[$(ii)$] if $N=d^*$ then $\Gamma_{\Omega^*}$ is locally a smooth manifold except possibly for a countable set of isolated singularities;
\item[$(iii)$] if $N>d^*$ then $\Gamma_{\Omega^*}$ is locally a smooth manifold except possibly for a set of Hausdorff dimension $N-d^*$.
\end{itemize}
\end{theorem}

This result is known for the volume-constrained minimizers of the first Dirichlet Laplacian eigenvalue and it is widely believed---even if not explicitly proved in the literature--- to also hold for the torsional energy problem, in a relative setting, such as the one considered here. For the sake of completeness we provide a proof in Sect. 5. This is accomplished with the arguments and a careful adaptation of the quite sophisticated techniques and methods developed in  \cite{B, BL, DS, KN, LS2, MTV, V, W1}.

The shape optimization problem discussed in Theorem \ref{mainteo1} is closely related to overdetermined problems of the form
\begin{equation}\label{over}
\begin{cases}
-\Delta u = 1 & \text{in $\Omega$},\\
 u = 0 & \text{on $\Gamma_\Omega$},\\
\frac{\partial u}{\partial \nu} = 0 & \text{on $\Gamma_{1,\Omega}$},\\
\frac{\partial u}{\partial \nu} = \alpha<0  & \text{on $\Gamma_{\Omega}$},\\
\end{cases}
\end{equation}
for some constant $\alpha$. Indeed, using domain derivative technique, as for similar problems in shape-optimization theory (see \cite{HP}), we can prove that if $\Omega^*\subset\C_{\omega}$ is a minimizer of $\E(\cdot; \C_\omega)$ in $\A_{\omega,c}$, for some fixed $c>0$, then its corresponding energy function $u_{\Omega^*}$ solves \eqref{over} (see Proposition \ref{lem2:sect5}). In fact, any critical point of $\E(\cdot; \C_\omega)$, with respect to volume-preserving deformations which leave the cylinder invariant is a domain whose corresponding torsion function satisfies \eqref{over}. 

Let us observe that, among the bounded domains contained in $\C_\omega$, we could consider a bounded cylinder 
\begin{equation}
\label{bdd-cyl}
\Omega_{\omega,h}:=\omega\times\left]-\frac h2,\frac h2\right[,
\end{equation}
with height $h>0$. The relative boundary of $\Omega_{\omega,h}$ is given by the union of two flat surfaces lying on the horizontal hyperplanes $x_N=\pm \frac h2$. It is easy to see that the torsion function corresponding to $\Omega_{\omega,h}$ is 
$$
u_{\Omega_{\omega,h}}(x)=u_{\Omega_{\omega,h}}(x_N)=\frac{h^2-4x_N^2}{8}
$$
which depends only on the $x_N$-variable and so has flat level sets. In fact, $u$ satisfies the overdetermined problem \eqref{over}. 

Now, the question is whether $\Omega_{\omega, h}$ is the only domain, for a fixed volume $c_h=h |\omega|$, for which the overdetermined torsion problem \eqref{over} has a solution, and whether $\Omega_{\omega, h}$ is a global minimizer of $\E(\cdot;\C_{\omega})$ in $\A_{\omega,c_{h}}$.  In this direction we recall the following theorem proved in \cite{AIP}.

\begin{theorem}\label{thm:sharp_stability_for_torsion}
Let $\lambda_1(\omega)$ be the first nontrivial Neumann eigenvalue of the Laplace operator $- \Delta_{\mathbb R^{N - 1}}$ in the domain $\omega \subset \mathbb R^{N - 1}$. Then there exists a number $\beta \approx 1,439$ (defined as the unique solution of $\sqrt s\,\tanh\sqrt s=1$) such that:
    \begin{itemize}
        \item[(i)] if $\lambda_1(\omega) < \frac{4\beta}{h^2}$ then $\Omega_{\omega,h}$ is not a local minimizer for the energy $\E(\cdot; \C_\omega)$;\vspace{6pt}
        \item[(ii)] if $\lambda_1(\omega) > \frac{4\beta}{h^2}$ then $\Omega_{\omega,h}$ is a local minimizer for the energy $\E(\cdot; \C_\omega)$.
    \end{itemize}
\end{theorem}
Note that in \cite{AIP} the above result is stated for $h=1$ and the semicylinder $\C_\omega^+$ (see \cite[Theorem 1.4]{AIP}) but is easy to reformulate it in $\C_\omega$ for any $h>0$. We also remark that Theorem \ref{thm:sharp_stability_for_torsion} does not state the existence of a minimizer but, combining it with Theorem \ref{mainteo1}, we obtain: 

\begin{theorem}\label{mainteo2}
Let $\omega \subset\R^{N-1}$ be a smooth bounded domain, and let $h>0$. 
If $\lambda_1(\omega) < \frac{4\beta}{h^2}$ then there exists a domain $\Omega_h^*$ in $\C_\omega$, with $|\Omega_h^*|=c_h=h |\omega|$, and different from $\Omega_{\omega,h}$, which is a minimizer of $\E(\cdot; \C_\omega)$ for the volume $c_h$. Moreover, in $\Omega_h^*$, the overdetermined problem \eqref{over} admits a solution.
\end{theorem}

In the recent paper \cite{PRS}, by a bifurcation argument, the existence of domains different from $\Omega_{\omega,h}$ for which \eqref{over} admits a solution is proved. However this result does not give any information on the existence of minimizers.
\medskip

When $\omega$ is a (bounded) domain in $\R^{N-1}$ of class $C^{2,\alpha}$ and the energy function $u_{\Omega^{*}}$ corresponding to a minimizer $\Omega^{*}$ for $\E(\cdot;\C_{\omega})$ given by Theorem \ref{mainteo1} has some more regularity, more precisely $u_{\Omega^{*}}\in W^{1,\infty}(\Omega^{*})$, we can provide some estimates on the $(N-1)$-dimensional Hausdorff measure $\H^{N-1}$ of the relative boundary of $\Omega^{*}$. 
 
In particular, by making a comparison with a bounded cylinder like \eqref{bdd-cyl} and using fine estimates (see Proposition \ref{P:bdr-estimate-1}), we can show that
\beq\label{eq:estimateHmintro}
\H^{N-1}(\Gamma_{\Omega^{*}})\leq 2\sqrt{3}\,\H^{N-1}(\omega).
\eeq
We stress that in \eqref{eq:estimateHmintro}, the $(N-1)$-dimensional Hausdorff measure of $\Gamma_{\Omega^{*}}$ is bounded from above by a constant independent of $c$. 
 %This means that for large $c$, the boundary of $\Omega^*$ cannot become too jagged. 
 We also note that, apart from the factor $\sqrt{3}$ in \eqref{eq:estimateHmintro}, the constant $2\H^{N-1}(\omega)=2|\omega|$ represents the measure of the relative boundary of the bounded cylinder $\Omega_{\omega,h}$ defined in \eqref{bdd-cyl}, for any $h>0$. In fact, motivated by this and the estimate obtained in \cite[Theorem 3.7]{AIP1}, we conjecture that for large $c$, the optimal domain $\Omega^{*}$ of Theorem \ref{mainteo1} is closer and closer to a bounded cylinder like $\Omega_{\omega, h}$ with height $h=\frac{c}{|\omega|}$.

Instead, a comparison with an (almost) half ball in $\C_{\omega}$ centered at a point of $\partial\C_{\omega}$ yields an estimate on the measure of the relative boundary of $\Omega^{*}$ with an upper bound vanishing as $c^{1-\frac1N}$, as $c\to 0^{+}$ (see Proposition \ref{P:bdr-estimate-2}). Hence, in this limit, the bounded cylinder \eqref{bdd-cyl} with $h=\frac{c}{|\omega|}$ cannot be the optimal domain $\Omega^{*}$ of Theorem \ref{mainteo1}, and this is consistent with part (i) of Theorem \ref{thm:sharp_stability_for_torsion}. In fact, we conjecture that for small $c$, the optimal domain $\Omega^{*}$ of Theorem \ref{mainteo1} becomes closer and closer to a half ball in $\C_{\omega}$, with volume $c$, centered at some point of $\partial\C_{\omega}$ which is a maximum point of the mean curvature of $\partial\C_{\omega}$. 

In the 2-dimensional case, the container is $\C_{\omega}=\left]0,a\right[\times\R$ for some $a>0$, and we can better describe the minimizers of $\E(\cdot; \C_\omega)$ in $\A_{\omega, c}$. More precisely, in Proposition \ref{prop:charelboundaryN2} we show that their closure always intersects $\partial\C_{\omega}$ along a set with positive 1-dimensional measure. Moreover, if $\Omega^{*}$ is a minimizer which is symmetric with respect to the horizontal axis $x_2=0$ and convex in the vertical direction and its torsion function $u_{\Omega^{*}}$ belongs to $W^{1,\infty}(\Omega^{*})$, then $\Omega^{*}$ can be either a half-disk attached to one of the boundary vertical lines $x_{1}=0$ or $x_{1}=a$, or a connected set whose relative boundary is the union of two symmetric curves, with respect to $x_2=0$, joining the two lines $x_{1}=0$ and $x_{1}=a$. Furthermore, when $c>\frac{3a^{2}}\pi$, only the second situation can occur. 

A natural question is whether also in higher dimension, at least in the case of convex cylinders, for large values of $c$, the projection of Steiner-symmetric optimal domains on the hyperplane $x_{N}=0$ completely covers the horizontal section $\omega$ of the container $\C_{\omega}$. 
\vspace{6pt}
 
 The paper is organized as follows. In Sect. \ref{S:preliminaries} we present some basic results on quasi-open sets and Steiner symmetrization. In Sect. \ref{S:sequences} we introduce the minimization problem and prove some preliminary results. Sect. \ref{S:minimizer} is devoted to the proof of the existence of a minimizer for the functional $\E(\cdot; \C_\omega)$ with a fixed volume. In Sect. \ref{S:properties} we prove the qualitative properties of a minimizer, as stated in Theorem \ref{mainteo1}. The connections with the overdetermined problem are emphasized in Sect. \ref{S:odtpb}, where we show also the estimates on the Hausdorff measure of the free boundary of minimizers, as well the study of the 2-dimensional case. Finally in Sect. \ref{S:half} we discuss the equivalence between the shape optimization problem in $\C_\omega$ and the analogous one in the half-cylinder $\C_\omega^+$.
\vspace{6pt}

For reader's convenience, let us recall some notation used in this work.
\begin{itemize}[leftmargin=15pt]
\item[$\bullet$] $\omega$ is a bounded domain in $\R^{N-1}$ ($N\ge 2$) with Lipschitz boundary, $\C_{\omega}=\omega\times\R$;\vspace{1.8pt}
\item[$\bullet$] a point $x\in\C_{\omega}$ is often written as $x=(x',x_{N})$ with $x'\in\omega$ and $x_{N}\in\R$; \vspace{1.8pt}
\item[$\bullet$] for $x\in\R^{N}$ and $r>0$ the open Euclidean ball centered at $x$ and with radius $r$ is denoted $B_{r}(x)$; \vspace{1.8pt}
\item[$\bullet$] $\mathrm{cap}(E)=$ capacity of a set $E\subset\R^{N}$ with respect to the $H^1$-norm;\vspace{1.8pt}
\item[$\bullet$] ``a.e.'' means almost everywhere with respect to the Lebesgue measure;\vspace{1.8pt}
\item[$\bullet$] ``q.e.'' means quasi everywhere, i.e., up to sets of zero capacity;\vspace{1.8pt}
\item[$\bullet$] $\mathcal{L}^{k}$ denotes the $k$-dimensional Lebesgue measure; we will often write $|E|$ instead of $\mathcal{L}^{k}(E)$ for $E$ measurable subset of $\R^{k}$;\vspace{1.8pt}
\item[$\bullet$] $\H^{k}$ denotes the $k$-dimensional Hausdorff measure and $\mathrm{dim}_\H(A)$ the Hausdorff dimension of a set $A\subset \RN$;\vspace{1.8pt}
\item[$\bullet$]  $d\sigma$ denotes the $(N-1)$-dimensional surface area element;\vspace{1.8pt}
\item[$\bullet$] $H^{1}(\C_{\omega})$ denotes the standard Sobolev space endowed with the norm $\|u\|_{H^{1}(\C_\omega)}^{2}=\|u\|_{L^{2}(\C_\omega)}^{2}+\|\nabla u\|_{L^{2}(\C_\omega)}^{2}$;\vspace{1.8pt}
\item[$\bullet$] for $\Omega\subset\C_{\omega}$ quasi-open, $H^{1}_{0}(\Omega;\C_{\omega})=\{u\in H^{1}(\C_{\omega})\,;~u=0~\text{q.e. in }\C_{\omega}\setminus\Omega\}$;\vspace{1.8pt}
\item[$\bullet$] for $E\subset\R^{N}$, $E^{\sharp}$ denotes the Steiner symmetrization of $E$ with respect to the hyperplane $\{x_N=0\}$;\vspace{1.8pt}
\item[$\bullet$] for $u\in H^{1}_{0}(\Omega;\C_{\omega})$ the energy of $u$ is $J_{\Omega}(u)=\int_{\C_{\omega}}\left(\frac12|\nabla u|^{2}-u\right)\,dx$;\vspace{1.8pt}
\item[$\bullet$] for $\Omega$ quasi-open subset of $\C_{\omega}$ with $|\Omega|<+\infty$ the torsional energy of $\Omega$ is $\E(\Omega;\C_{\omega})=\displaystyle\inf\{J_{\Omega}(u)\,;~u\in H^{1}_{0}(\Omega;\C_{\omega})\}$ and a function $u_{\Omega}\in H^{1}_{0}(\Omega;\C_{\omega})$ such that $J_{\Omega}(u_{\Omega})=\E(\Omega;\C_{\omega})$ is the corresponding energy function or torsion function;\vspace{1.8pt}
\item[$\bullet$] $\mathcal{A}_{\omega,c}=\{\Omega\subset\C_{\omega}\,;~\Omega\text{ quasi-open, }|\Omega|\le c\}$ and $\mathcal{O}_{c}(\C_{\omega})=\inf\{\E(\Omega;\C_{\omega})\,;~\Omega\in \mathcal{A}_{\omega,c}\}$, for $c>0$;\vspace{1.8pt}
\item[$\bullet$] for $\Omega$ quasi-open subset of $\C_{\omega}$ the relative boundary of $\Omega$ is $\Gamma_{\Omega}=\partial\Omega\cap\C_{\omega}$; moreover $\Gamma_{1,\Omega}=\partial\Omega\cap\partial\C_{\omega}$;\vspace{1.8pt}
%\item[$\bullet$] $\mathrm{Reg}(\Gamma_{\Omega^*})$,  $\mathrm{Sing}(\Gamma_{\Omega^*})$ denote, respectively, the regular and the singular part of  $\Gamma_{\Omega^*}$ (see the proof of Theorem \ref{teoreg}).
%\item[$\bullet$] for $\Omega$ bounded open subset of $\C_{\omega}$, $\partial\Omega^{\mathrm{reg}}$ is the regular part of $\partial{\Omega}$. 
\item In Sect. \ref{S:half} we will study the problem in the half-cylinder $\C^{+}_{\omega}=\omega\times\left]0,+\infty\right[$. The superscript $+$ will also be used with obvious meaning in other notation, like $J^{+}_{\Omega}$, $\A^{+}_{\omega,c}$.
\end{itemize}

\section{Preliminaries on the functional framework and Steiner symmetrization}\label{S:preliminaries}

In the first part of this Section we recall the notions of capacity and of quasi-open set. Let $\omega$ be a bounded domain of $\R^{N-1}$ with Lipschitz boundary and let $\mathcal{C}_\omega:=\omega\times\R$ be the cylinder generated by $\omega$.
% (see \eqref{def:cylspanned}).
\begin{definition}\label{def:capacity}
For a generic set $E \subset \RN$, we define the capacity of $E$ (with respect to the $H^1$-norm), and we denote it by $\mathrm{cap}(E)$, the number
$$\mathrm{cap}(E)=\mathrm{inf}\left\{\|u\|_{H^1(\RN)}^2; \ u\in H^1(\RN), \ u\geq 1 \ \hbox{in a neighborhood of $E$}\right\},$$
where $\|\cdot\|_{H^1(\RN)}$ is the standard norm in $H^1(\RN)$.
\end{definition}
We recall that for every $\mathcal{L}^{N}$-measurable set $E\subset\R^{N}$ one has $\mathrm{cap}(E)\ge|E|:=\mathcal{L}^{N}(E)$, and thus, in particular, the sets of zero capacity are negligible under the Lebesgue measure. 
\begin{definition}\label{def:quasiopen}
We say that $\Omega\subset \C_\omega$ is quasi-open, if for any $\varepsilon>0$, there exists an open set $\Lambda_\e \subset \C_\omega$ such that $\mathrm{cap}(\Lambda_\e)\leq \e$ and $\Omega\cup\Lambda_\e$ is open.
\end{definition}

\begin{definition}
We say that $u:\C_\omega\to \R$ is quasi-continuous, if for any $\varepsilon>0$, there exists an open set $\Lambda_\e \subset \C_\omega$ such that $\mathrm{cap}(\Lambda_\e)\leq \e$ and the restriction of $u$ on the set $\C_\omega\setminus\Lambda_\e$ is continuous.
\end{definition}

We recall that any $u\in H^1(\C_\omega)$ has a quasi-continuous representative (see e.g. \cite[Sect. 4.8]{EG}). In the present work we will always identify any $u\in H^1(\C_\omega)$ with its quasi-continuous representative, so that, the superlevel set $\{u>t\}$ is a quasi-open for any $t\in \R$ (see \cite[Sect. 2]{BV}).
For any quasi-open set $\Omega\subset \C_\omega$ we consider the Sobolev space:
$$
H^1_0(\Omega; \C_\omega):=\left\{u \in H^1(\C_\omega); \ \ u=0 \ \ \hbox{q.e. on} \ \C_\omega\setminus\Omega  \right\},
$$
where q.e. means quasi-everywhere, i.e., up to sets of zero capacity.

The space $H^1_0(\Omega; \C_\omega)$ is a closed subspace of $H^1(\C_\omega)$ and it is the appropriate functional space to study the torsion problem
\begin{equation}\label{eq:mixbvprobquasiopen}
\begin{cases}
-\Delta u = 1 & \text{in $\Omega$,}\\
 u = 0 & \text{on $\partial\Omega \cap \C_\omega$,}\\
 \frac{\partial u}{\partial \nu} = 0 & \text{on $\partial \C_\omega$.}
\end{cases}
\end{equation}

\begin{definition}
Given a quasi-open set $\Omega\subset \C_\omega$, a (weak) solution of \eqref{eq:mixbvprobquasiopen} is a function $u\in H^1_0(\Omega; \C_\omega)$ such that
$$
\int_{\C_\omega} \nabla u \boldsymbol{\cdot} \nabla v \ dx = \int_{\C_\omega} v \ dx \ \ \ \ \forall v \in H^1_0(\Omega; \C_\omega).
$$
\end{definition}
Hence, a (weak) solution of \eqref{eq:mixbvprobquasiopen} is a critical point of the functional $J_\Omega\colon H^1_0(\Omega; \C_\omega)\to \R$ defined by
\beq\label{eq:deftorsion}
J_\Omega(u)=\frac{1}{2}\int_{\C_\omega} |\nabla u|^2 \ dx - \int_{\C_\omega} u \ dx.
\eeq
Notice that the dependence of $J_{\Omega}$ on $\Omega$ appears in the domain $H^1_0(\Omega; \C_\omega)$. Using the definitions of capacity and of the space $H^1_0(\Omega; \C_\omega)$ and recalling that if $|\Omega|<+\infty$ the inclusion $H^1_0(\Omega; \C_\omega)\hookrightarrow L^2(\C_\omega)$ is compact (see \cite[Proposition 2.3-(i)]{BV}), with standard arguments one plainly obtains the following result:
 
\begin{lemma}\label{energy-function}
For every quasi-open set $\Omega\subset \C_{\omega}$ with $|\Omega|<+\infty$, there exists a unique (weak) solution of \eqref{eq:mixbvprobquasiopen}, denoted by $u_{\Omega}$ and characterized by
$$
J_{\Omega}(u_{\Omega})=\min_{u\in H^1_0(\Omega; \C_\omega)}J_{\Omega}(u).
$$
Such a function $u_{\Omega}$ is called \emph{energy function} or \emph{torsion function of} $\Omega$ and satisfies
\begin{gather}
\nonumber
\Omega=\{u_\Omega>0\} \text{ up to a set of zero capacity,}\\
\label{def:torsionquasiopen}
J_\Omega(u_\Omega)=-\frac{1}{2} \int_\Omega |\nabla u_\Omega|^2 \ dx=-\frac{1}{2} \int_\Omega u_\Omega \ dx.
\end{gather}
\end{lemma}

%% \begin{proof}
%% If $|\Omega|<+\infty$ the inclusion $H^1_0(\Omega; \C_\omega)\hookrightarrow L^2(\C_\omega)$ is compact (see \cite[Proposition 2.3-(i)]{BV}). This implies that $J_{\Omega}$ has a minimizer $u_\Omega \in H^1_0(\Omega; \C_\omega)$. Then $u_{\Omega}$ is a (weak) solution of \eqref{eq:mixbvprobquasiopen}. Hence 
%% $$J_{\Omega}(u_{\Omega}+v)-J_{\Omega}(u_{\Omega})=\frac12\int_\Omega |\nabla v|^2\,  dx~~\forall v\in H^1_0(\Omega; \C_\omega)$$
%% and this implies uniqueness. Property \eqref{u>0} is proved in \cite[Proposition 2.8-(e)]{BV}. Finally \eqref{def:torsionquasiopen} is obtained by taking $v=u_{\Omega}$ in \eqref{weaksol}.
%% \end{proof}

\begin{remark}\label{uniqueness}
We point out that in general, $u\in H^1_0(\Omega; \C_\omega)$ and $v\in H^{1}(\C_{\omega})$ with $v=u$ a.e. in $\C_{\omega}$ does not imply $v\in H^1_0(\Omega; \C_\omega)$. This is true when $v=u$ q.e. in $\C_{\omega}$. Hence the uniqueness of the weak solution to \eqref{eq:mixbvprobquasiopen} is up to a set of null capacity. In particular, we can assume that the energy function $u_{\Omega}$ satisfies $u_{\Omega}(x)=0$ for every $x\in\C_{\omega}\setminus\Omega$. Indeed, if $u_{\Omega}\in H^1_0(\Omega; \C_\omega)$ is an energy function corresponding to $\Omega$, there exists a set $\Lambda\subset\C_{\omega}$ with $\mathrm{cap}(\Lambda)=0$ such that $u_{\Omega}(x)=0$ for every $x\in\C_{\omega}\setminus(\Omega\cup\Lambda)$. Now, taking
$$
\tilde{u}_{\Omega}(x)=\begin{cases}u_{\Omega}(x)&\text{if $x\in\C_{\omega}\setminus\Lambda$}\\0&\text{if $x\in\Lambda$,}\end{cases}
$$
also $\tilde{u}_{\Omega}\in H^1_0(\Omega; \C_\omega)$ and $\tilde{u}_{\Omega}=u_{\Omega}$ q.e. in $\C_{\omega}$. In particular, $\tilde{u}_{\Omega}$ is an energy function satisfying $\tilde{u}_{\Omega}(x)=0$ for every $x\in\C_{\omega}\setminus\Omega$. From now on, for every quasi-open set $\Omega\subset\C_{\omega}$ with $|\Omega|<+\infty$, we denote $u_{\Omega}$ an energy function such that $u_{\Omega}(x)=0$ for every $x\in\C_{\omega}\setminus\Omega$.
\end{remark}

In the sequel we are going to study the problem of minimizing the functional 
$$
\E(\Omega;\C_\omega):=J_{\Omega}(u_{\Omega})
$$ 
among quasi-open sets of uniformly bounded measure. To this aim, the following estimates are useful. 

\begin{proposition}\label{prop1:sect2}
Let $c>0$ and let $\omega\subset\R^{N-1}$ be a bounded domain with Lipschitz boundary. Then, there exists a positive constant $C$ depending only on $N$, $\omega$ and $c$ such that for any quasi-open subset $\Omega$ of $\C_\omega$ with $|\Omega| \leq c$, it holds:\vspace{3pt}
\begin{itemize}
\item[(i)] $u_\Omega$ is bounded and $\|u_\Omega\|_{L^\infty(\C_\omega)}\leq C |\Omega|^{2/N}$;\vspace{6pt}
\item[(ii)] $\displaystyle\int_{\C_\omega} |\nabla u_\Omega|^2\  dx \leq C |\Omega|^{\frac{N+2}{N}}$;\vspace{3pt}
\item[(iii)] $\displaystyle\int_{\C_\omega} u_\Omega^2\  dx \leq C |\Omega|^{\frac{N+4}{N}}$.
\end{itemize}
 \end{proposition}  
\begin{proof}
The proof relies on \cite[Lemma 2.5]{BV} and it is the same as in \cite[Proposition 6.4]{IPW} with slight adjustments. In order to apply  \cite[Lemma 2.5]{BV}  we point out that $\C_\omega$ is a domain with uniformly Lipschitz boundary (because $\omega$ is bounded and has Lipschitz boundary).
\end{proof}

The following additivity property will be useful in the sequel. 

\begin{lemma}\label{additivity}
Let $\omega\subset\R^{N-1}$ be a bounded domain with Lipschitz boundary and let $\Omega_{1},\Omega_{2}\subset\C_{\omega}$ be two disjoint quasi-open sets with $|\Omega_{1}|<+\infty$ and $|\Omega_{2}|<+\infty$. Then
\begin{equation}\label{subadditivity-energy}
\E(\Omega_{1}\cup\Omega_{2};\C_{\omega})\leq\E(\Omega_{1};\C_{\omega})+\E(\Omega_{2};\C_{\omega}).
\end{equation}
If in addition $\mathrm{cap}(\partial\Omega_{1}\cap\partial\Omega_{2}\cap\C_{\omega})=0$,  %and $|\Omega_{1}\cup\Omega_{2}|<+\infty$, 
then 
\begin{equation}\label{additivity-energy}
\E(\Omega_{1}\cup\Omega_{2};\C_{\omega})=\E(\Omega_{1};\C_{\omega})+\E(\Omega_{2};\C_{\omega}).
\end{equation}
\end{lemma}
\begin{proof} For $i=1,2$, let $u_{\Omega_{i}}\in H^{1}_{0}(\Omega_{i};\C_{\omega})$ be the energy function corresponding to $\Omega_{i}$ ($i=1,2$). It holds that $u_{\Omega_1}+u_{\Omega_2}\in H^{1}_{0}(\Omega_{1}\cup\Omega_{2};\C_{\omega})$, and since $\Omega_{1}\cap\Omega_{2}=\varnothing$, we have
$$
\E(\Omega_{1}\cup\Omega_{2};\C_{\omega})\leq J_{\Omega_{1}\cup\Omega_{2}}(u_{\Omega_1}+u_{\Omega_2})=J_{\Omega_{1}}(u_{\Omega_1})+J_{\Omega_{2}}(u_{\Omega_2})=\E(\Omega_{1};\C_{\omega})+\E(\Omega_{2};\C_{\omega}).
$$
To prove \eqref{additivity-energy}, consider the function $\hat{u}=u_{\Omega_{1}\cup\Omega_{2}}$ which minimizes $J_{\Omega_{1}\cup\Omega_{2}}$ in $H^{1}_{0}(\Omega_{1}\cup\Omega_{2};\C_{\omega})$ and for $i=1,2$ define
$$
\hat{u}_{i}=\begin{cases}\hat{u}&\text{ in }\Omega_{i}\\ 0&\text{ in }\C_{\omega}\setminus\Omega_{i}.\end{cases}
$$
Then, since $\mathrm{cap}(\partial\Omega_{1}\cap\partial\Omega_{2}\cap\C_{\omega})=0$, by Remark \ref{uniqueness}, we have that $\hat{u}_{i}\in H^{1}_{0}(\Omega_{i};\C_{\omega})$ and 
$$
\E(\Omega_{1};\C_{\omega})+\E(\Omega_{2};\C_{\omega})\le J_{\Omega_{1}}(\hat{u}_{1})+ J_{\Omega_{2}}(\hat{u}_{2})=J_{\Omega_{1}\cup\Omega_{2}}(\hat{u})=\E(\Omega_{1}\cup\Omega_{2};\C_{\omega}).
$$
This, together with \eqref{subadditivity-energy}, proves \eqref{additivity-energy}.\end{proof}

We now recall the definition of Steiner symmetrization with respect to the hyperplane $\{x_N = 0\}$ and review some properties that will be useful throughout the paper.
\medskip

Let $E$ be a subset of $\R^{N}$. For every $x'\in\R^{N-1}$ let 
$$
E_{x'}:=\{x_{N}\in\R\,;~x=(x',x_{N})\in E\}
$$
be the one-dimensional cross-section of $E$ parallel to the $x_{N}$-axis and set
$$
\ell_{E}(x'):=\begin{cases}|E_{x'}|=\mathcal{L}^{1}(E_{x'})&\text{ if $E_{x'}$ is $\mathcal{L}^{1}$-measurable}\\ 0&\text{ otherwise.}\end{cases}
$$
\begin{definition}\label{def:SteinerSymmetrization}
We define the Steiner symmetrization of $E$ with respect to the hyperplane $\{x_N=0\}$ the set
$$
E^{\sharp}:=%\bigcup_{x^\prime\in \R^{N-1}} (\{x^\prime\}\times E_{x'})
\left\{(x',x_{N})\in\R^{N}\,;~|x_{N}|<\frac12\ell_{E}(x')\right\}.
$$
\end{definition}

By Fubini's Theorem, if $E$ is a $\mathcal{L}^N$-measurable set, then $E_{x'}$ is a $\mathcal{L}^1$-measurable subset of $\R$ for $\mathcal{L}^{N-1}$-almost every $x'\in\R^{N-1}$ and the mapping $\ell_{E}\colon\R^{N-1}\to\left[0,\infty\right]$ is measurable. Hence also $E^{\sharp}$ is a $\mathcal{L}^N$-measurable set because $E^{\sharp}=\{h_{E}<0\}$ where 
$$
h_{E}(x)=|x_{N}|-\frac12\ell_{E}(x')~~\forall x=(x',x_{N})\in\R^{N}.
$$
In addition, 
\begin{equation}
\label{eq:steinerpreservesvolume}
|E|=|E^{\sharp}|
%\int_{E}dx=\int_{\R^{N-1}}\left[\int_{E_{x'}}dx_{N}\right]dx'=\int_{\R^{N-1}}\ell_{E}(x')\,dx'=\int_{\R^{N-1}}\left[\int_{-\frac12\ell_{E}(x')}^{\frac12\ell_{E}(x')}dx_{N}\right]dx'=\int_{E^{\sharp}}dx=
\end{equation}
(see \cite[Sect. 6]{BA} or \cite{BUR}).

Other possible equivalent choices for the definition of the Steiner symmetrization of a set are available in the literature (see e.g. \cite[Sect. 6]{BA} or \cite[(3.8)]{MA}). For our purposes we follow \cite[Sect. 6.1.2]{HP}, as we take advantage from the following further properties:

\begin{lemma}\label{opensharp}
If $\Omega$ is an open subset of $\R^{N}$, then:\vspace{3pt}
\begin{itemize}
\item[(i)]
$\Omega_{x'}:=\{x_{N}\in\R\,;~(x',x_{N})\in\Omega\}$ is an open subset of $\R$ for every $x'\in\R^{N-1}$;\vspace{6pt}
\item[(ii)]
if $\Omega$ is bounded, the mapping $x'\mapsto\ell_{\Omega}(x'):=|\Omega_{x'}|$ is lower semicontinuous in $\R^{N-1}$;\vspace{6pt}
\item[(iii)]
$\Omega^{\sharp}$ is an open subset of $\R^{N}$. If in addition $\Omega\subset\C_{\omega}$, then also $\Omega^{\sharp}\subset\C_{\omega}$. 
\end{itemize}
\end{lemma}

\begin{proof} (i) Fix $x'\in\R^{N-1}$. If $x_{N}\in \Omega_{x'}$ then $(x',x_{N})\in\Omega$ and since $\Omega$ is open, there exists $r>0$ such that $B_{r}(x',x_{N})\subset\Omega$, where $B_{r}(x',x_{N})=\{(y',y_{N})\in\R^{N}\,;~|y'-x'|^{2}+|y_{N}-x_{N}|^{2}<r^{2}\}$. In particular $(x',y_{N})\in\Omega~\forall y_{N}\in\left]x_{N}-r,x_{N}+r\right[$, namely $\left]x_{N}-r,x_{N}+r\right[\subset \Omega_{x'}$. Thus (i) is proved. 
\medskip

\noindent
(ii) Since $\Omega$ is bounded, for every $x'\in\R^{N}$ the set $\Omega_{x'}$ is a bounded subset of $\R$. Hence, by (i), $|\Omega_{x'}|\in\left[0,\infty\right[$. Fixing $x'\in\R^{N-1}$ such that $\Omega_{x'}\ne\varnothing$, for every $\varepsilon>0$ there exists a compact subset $K$ of $\Omega_{x'}$ such that $|K|>|\Omega_{x'}|-\varepsilon$. Then $\{x'\}\times K$ is a compact subset of $\Omega$. Hence there exists $\delta>0$ such that $B'_{\delta}(x')\times K\subset \Omega$ where $B'_{\delta}(x')=\{y'\in\R^{N-1}\,;~|y'-x'|<\delta\}$. Then
$$
|\Omega_{y'}|\ge |K|> |\Omega_{x'}|-\varepsilon~~\forall y'\in B'_{\delta}(x'),
$$
that is, the mapping $\ell_{\Omega}$ is lower semicontinuous. 
\medskip

\noindent
(iii) Firstly, suppose $\Omega$ bounded. By (ii), the function $h_{\Omega}\colon\R^{N}\to\R$ given by
$
h_{\Omega}(x):=|x_{N}|-\frac12\ell_{\Omega}(x')~~\forall x=(x',x_{N})\in\R^{N}
$
is upper semicontinuous on $\R^{N}$. Then $\Omega^{\sharp}=\{h_{\Omega}<0\}$ is open. Let us discuss the general case. In particular, let us prove that 
\begin{equation}\label{omegaunbounded}
\Omega^{\sharp}=\bigcup_{r>0}\left(\Omega\cap B_{r}\right)^{\sharp}\quad\text{where~~}B_{r}=\{x\in\R^{N}~:~|x|<r\}.
\end{equation}
Indeed $\Omega\cap B_{r}\subseteq\Omega$ implies $\ell_{\Omega\cap B_{r}}\le\ell_{\Omega}$, then $\left(\Omega\cap B_{r}\right)^{\sharp}\subseteq\Omega^{\sharp}$ for every $r>0$ and thus $\bigcup_{r>0}\left(\Omega\cap B_{r}\right)^{\sharp}\subseteq\Omega^{\sharp}$. If $x=(x',x_{N})\in\Omega^{\sharp}$ then $|x_{N}|<\frac12\ell_{\Omega}(x')$. In particular $\ell_{\Omega}(x')>0$. Since 
$$
\ell_{\Omega\cap B_{r}}(x')=\left|\left(\Omega\cap B_{r}\right)_{x'}\right|\to\left|\Omega_{x'}\right|=\ell_{\Omega}(x')\quad\text{as~}r\to\infty,
$$ 
there exists $r>0$ such that $|x_{N}|<\frac12\ell_{\Omega\cap B_{r}}(x')\le\frac12\ell_{\Omega}(x')$. That is, $x\in \left(\Omega\cap B_{r}\right)^{\sharp}$. Therefore $\bigcup_{r>0}\left(\Omega\cap B_{r}\right)^{\sharp}\supseteq\Omega^{\sharp}$ and \eqref{omegaunbounded} is proved. By the first part of the proof, $\left(\Omega\cap B_{r}\right)^{\sharp}$ is open for every $r>0$. Hence, by \eqref{omegaunbounded}, also $\Omega^{\sharp}$ is open. In addition, if $\Omega\subset\C_{\omega}$ then $\Omega_{x'}=\varnothing$ for every $x'\in\R^{N-1}\setminus\omega$ and consequently $\Omega^{\sharp}\subset\C_{\omega}$. 
\end{proof}

We recall now the definition of Steiner symmetrization of a non-negative measurable function.
\medskip

Let $u\colon\C_\omega\to \R^+\cup\{0\}$ be a non-negative $\mathcal{L}^N$-measurable function, for $x^\prime\in\omega$ we define the slice function $u^{x^\prime}\colon\R\to \R^+\cup\{0\}$ by
$$u^{x^\prime}(x_N):=u(x^\prime,x_N), \ \ x_N\in\R,$$
and assume that for $\mathcal{L}^{N-1}$-almost every $x^\prime\in\omega$ the slice function $u^{x^\prime}$ satisfies the finiteness condition
\beq\label{eq:deffinitenesscond}
 \mathcal{L}^1(\{u^{x^\prime}>t\})<+\infty \ \ \forall t>0.
\eeq
Notice that if \eqref{eq:deffinitenesscond} is satisfied, the one-dimensional symmetric decreasing rearrangement of $u^{x^\prime}$, denoted by $(u^{x^\prime})^\sharp$, is well-defined. 

We are now in position to give the following.
\begin{definition}
Let $u\colon\C_\omega\to \R^+\cup\{0\}$ be a non-negative $\mathcal{L}^N$-measurable function, satisfying \eqref{eq:deffinitenesscond}. We define the Steiner symmetrization of $u$ on $\C_\omega$ as the function $u^\sharp\colon\C_\omega\to \R^+\cup\{0\}$ by setting
$$
u^\sharp(x^\prime, x_N):= \begin{cases} (u^{x^\prime})^\sharp(x_N) & \hbox{if $u^{x^\prime}$ satisfies \eqref{eq:deffinitenesscond}},\\
0 & \hbox{otherwise},
\end{cases}
$$
 where $(u^{x^\prime})^\sharp$ is the one-dimensional symmetric decreasing rearrangement of $u^{x^\prime}$.
\end{definition}
If $u \in L^p(\C_\omega; \R^+\cup\{0\})$, with $p\in [1, +\infty[$, then by Fubini's theorem and the layer-cake formula it follows that \eqref{eq:deffinitenesscond} is satisfied for $\mathcal{L}^{N-1}$-almost every $x^\prime\in\omega$, and it is well known that
\begin{equation}\label{eq:Lp-normpreservedSteiner}
\|u\|_{L^p(\C_\omega)}=\|u^\sharp\|_{L^p(\C_\omega)}.
\end{equation}
If in addition $u$ belongs to the Sobolev space $W^{1,p}$, then the Steiner symmetrization decreases the $p$-Dirichlet integrals, namely we have the following.
\begin{proposition}\label{prop2:sect2}
Let $p\in [1, +\infty[$ and let $u\in W^{1,p}(\C_\omega; \R^+\cup\{0\})$ be a non-negative function. Then $u^\sharp \in W^{1,p}(\C_\omega; \R^+\cup\{0\})$ and
$$\int_{\C_\omega} |\nabla u^\sharp|^p\, dx \leq \int_{\C_\omega} |\nabla u|^p\, dx.$$
\end{proposition}
\begin{proof}
The proof is the same of \cite[Theorem 6.19]{BA} (see also \cite{BUR} or \cite[Sect. 6.1.2]{HP} and the references therein).
\end{proof}
As a consequence of Proposition \ref{prop2:sect2} we obtain the following result that will be used throughout the paper.
\begin{lemma}\label{lem1:sect2}
Let $\Omega$ be an open subset of $\C_\omega$ and let $u_\Omega\in H_0^1(\Omega;\C_\omega)$ be its energy function. Then:\vspace{3pt}
\begin{itemize}
\item[(i)] $u_\Omega^\sharp \in H_0^1(\Omega^\sharp;\C_\omega)$;\vspace{6pt}
\item[(ii)] $\E(\Omega^\sharp; \C_\omega) \leq \E(\Omega; \C_\omega)$.
\end{itemize}
\end{lemma}
\begin{proof}
(i) Let $u_\Omega\in H_0^1(\Omega;\C_\omega)$ be the energy function corresponding to an open set $\Omega\subset\C_{\omega}$. Then $u_\Omega\in H^1(\C_\omega)$ and $u_\Omega(x)=0$ for every $x\in\C_\omega\setminus\Omega$ (see Remark \ref{uniqueness}). Thanks to Proposition \ref{prop2:sect2}, $u_\Omega^\sharp\in H^1(\C_\omega)$, and by definition of Steiner symmetrization, $u_\Omega^\sharp=0$ in $\C_\omega\setminus \Omega^{\sharp}$. Thus, $u_\Omega^\sharp\in H_0^1(\Omega^\sharp;\C_\omega)$. 
\medskip

\noindent
(ii)
Using \eqref{eq:Lp-normpreservedSteiner}, Proposition \ref{prop2:sect2} and part (i), we infer that $u_{\Omega}^{\sharp}$ is a minimizer of $J_{\Omega^\sharp}$ in $H^{1}_{0}(\Omega^\sharp; \C_\omega)$. Hence, by uniqueness, $u_{\Omega}^{\sharp}=u_{\Omega^{\sharp}}$ and then
$
\E(\Omega^\sharp; \C_\omega) =J_{\Omega^{\sharp}}(u_{\Omega^{\sharp}})= J_{\Omega^{\sharp}}(u_{\Omega}^{\sharp})\le J_{\Omega}(u_{\Omega})=\E(\Omega; \C_\omega).
$
\end{proof}

\section{The minimization problem}
\label{S:sequences}

In this section we consider the minimization problem for the torsional energy, under a volume constraint, when the container is a cylinder or a half-cylinder spanned by a Lipschitz bounded domain $\omega$ in $\R^{N-1}$. We study the properties of the minimizing sequences and determine a relation between the energy levels of the two problems.
\medskip 

Fixing $c>0$ we denote by $\A_{\omega,c}$ the class of quasi-open sets in $\C_\omega$ of measure less or equal than $c$, namely
\begin{equation}\label{eq:defclassadmissiblecylinder}
\A_{\omega,c}:= \{\Omega\subset \C_\omega;\ \Omega \ \hbox{quasi-open} \ \hbox{and} \ |\Omega|\leq c\},
\end{equation}
and define
$$
\mathcal{O}_c(\C_\omega):=\inf\{\E(\Omega;\C_\omega); \  \Omega \in \A_{\omega; c}\}.
$$
\begin{remark}\label{remark:negativityenergy}
For any $c>0$ it holds $\O_c(\C_\omega)<0$, as the torsional energy of the flat bounded cylinder of volume $c$ is negative. Moreover, by Proposition \ref{prop1:sect2} it immediately follows that $\O_c(\C_\omega)>-\infty$. %The same conclusions hold for $\O_c(\C_\omega^+)$.
\end{remark}

Following \cite{BV}, let us recall the following

\begin{definition}[$\gamma$-convergence]
Let $\Omega_{n}\subset\C_{\omega}$ be a sequence of quasi-open sets of finite measure. We say that $\Omega_{n}$ $\gamma$-converges to the quasi-open set $\Omega\subset\C_{\omega}$ if the sequence of energy functions $u_{\Omega_{n}}\in H^{1}_{0}(\Omega_{n};\C_{\omega})$ converges strongly in $L^{2}(\C_{\omega})$ to the energy function $u_{\Omega}\in H^{1}_{0}(\Omega;\C_{\omega})$.
\end{definition}

\begin{lemma}\label{lem:tech2}
Let $c>0$, let $\omega\subset\R^{N-1}$ be a Lipschitz bounded domain and let $\Omega\in \A_{\omega, c}$. Then, for any sequence of positive real numbers $(\eta_n)_n$ such that $\eta_n\to 0^+$, as $n\to +\infty$,  there exists a sequence $(\Omega_n)_n$ of open subsets of $\C_\omega$ such that $|\Omega_n|\leq c+\eta_n$ for all $n\in\N$ and $\Omega_n$ $\gamma$-converges to $\Omega$, as $n\to +\infty$. 
\end{lemma}

\begin{proof}
Let $c > 0$ and $\omega$ as in the statement. Let $\Omega \subset \C_\omega$ be a quasi-open set with $|\Omega|\leq c$ and let $u_\Omega \in H^1_0(\Omega;\C_\omega)$ be its energy function. By definition of quasi-open set, given $(\eta_n)_n$ such that $\eta_n\to 0^+$, as $n\to +\infty$, we find a sequence of open sets $(\Lambda_n)_n $ of  $\C_\omega$ such that $\Omega_n:=\Omega\cup \Lambda_n$ is an open set and $\mathrm{cap}(\Lambda_n)\leq \eta_n$ for all $n\in \N^+$. Accordingly, by definition of capacity, we find a sequence $(v_n)_n \subset H^1(\RN)$ such that $v_n \geq 1$ in a neighborhood of $\Lambda_n$ and $\|v_n\|_{H^1(\RN)}\to 0$, as $n\to +\infty$. Passing to a subsequence, still denoted $(v_{n})_{n}$, we can also assume that $v_{n}\to 0$ a.e. in $\RN$ as $n\to+\infty$. 
Then, setting $$\varphi_n:=\min\{1,|v_n|\}$$ we readily check that $0\leq \varphi_n\leq 1$ a.e. in $\RN$, $\varphi_n \equiv 1$ in a neighborhood of $\Lambda_n$, for any $n\in \N$. Moreover, by construction, $(\varphi_n)_n$ is a bounded sequence in $H^1(\RN)$ and, as $n\to +\infty$, it holds
\begin{equation}\label{eq-1:lemtech2}
 \|\varphi_n\|_{H^1(\RN)}\to 0\quad\text{and}\quad \varphi_{n}\to 0\quad\text{a.e. in }\RN.
 \end{equation}
Now, let $u_{\Omega_n} \in H_0^1(\Omega_n; \C_\omega)$ be the energy function of $\Omega_n$. Recalling that $|\Lambda_n|\leq\mathrm{cap}(\Lambda_n)$ (see \cite[Sect. 2.1]{BV}) we have that $|\Omega_n|\leq c+\mathrm{cap}(\Lambda_n)\leq c+ \eta_n$ and thus by Proposition \ref{prop1:sect2} we infer that there exists a positive constant $C_1$ such that
\begin{equation}\label{eq0:lemtech2}
\|u_{\Omega_n} \|_{L^\infty(\C_\omega)} +\|u_{\Omega_n} \|_{H^1(\C_\omega)} \leq C_1 \ \ \forall n\in \N^+.
\end{equation}

We aim to show that $u_{\Omega_n}\to u_\Omega$ strongly in $L^2(\C_\omega)$, as $n\to +\infty$. To this end, let us first write
\begin{equation}\label{eq-3:lemtech2}
u_{\Omega_n}=u_{\Omega_n} \varphi_n + u_{\Omega_n} (1-\varphi_n)
\end{equation}
and observe that by H\"older's inequality, \eqref{eq-1:lemtech2}, \eqref{eq0:lemtech2} it follows that
\begin{equation}\label{eq-2:lemtech2}
\|u_{\Omega_n} \varphi_n\|_{L^2(\C_\omega)}\to 0,
\end{equation}
as $n\to +\infty$. Hence, to conclude, it remains to show that $u_{\Omega_n} (1-\varphi_n) \to u_\Omega$ strongly in $L^2(\C_\omega)$.
\medskip

To prove this, observe that by construction and by \eqref{eq0:lemtech2} we have that $u_{\Omega_n} (1-\varphi_n) \in H_0^1(\Omega; \C_\omega)$, for all $n\in\N$, and $(u_{\Omega_n} (1-\varphi_n))_n$ is a bounded sequence in $H_0^1(\Omega; \C_\omega)$. Then, up to a subsequence, still indexed by $n$, we deduce that 
\begin{equation}\label{eq1:lemtech2}
u_{\Omega_n} (1-\varphi_n) \rightharpoonup u,
\end{equation}
for some $u\in H_0^1(\Omega; \C_\omega)$. Moreover, as $\C_\omega$ is a uniformly Lipschitz domain and $\Omega$ is of finite measure, then the imbedding $H_0^1(\Omega; \C_\omega)\hookrightarrow L^2(\C_\omega)$ is compact (see \cite[(2.2)]{BV}) and thus, up to a further subsequence, we get that
\begin{equation}\label{eq2:lemtech2}
 u_{\Omega_n} (1-\varphi_n) \to u
\end{equation}
strongly in $L^2(\C_\omega)$, as $n\to +\infty$. We claim that $u=u_\Omega$.
\medskip

 Indeed, for any test function $\varphi \in H_0^1(\Omega; \C_\omega)$ it holds

\begin{equation}\label{eq3:lemtech2}
\begin{array}{lll}
&&\displaystyle\int_{\Omega} \nabla\left[u_{\Omega_n} (1-\varphi_n)\right]\cdot \nabla \varphi\, dx\\[12pt]
&=&\displaystyle\int_{\Omega} \nabla u_{\Omega_n}\cdot \nabla \varphi \, dx-\underbrace{\int_{\Omega} \left(\nabla u_{\Omega_n}\cdot \nabla \varphi\right) \varphi_n \, dx}_{\mathbf{I}_{1,n}} -\underbrace{\int_{\Omega} \left(\nabla \varphi_n \cdot \nabla \varphi\right) u_{\Omega_n}\  \, dx}_{\mathbf{I}_{2,n}}\\[-9pt]
&=&\displaystyle\int_{\Omega} \varphi\, dx +o(1)
\end{array}
\end{equation}
where in the last equality we used that $u_{\Omega_n}$ is the energy function of $\Omega_n$ and $\Omega\subset \Omega_n$, while to estimate the terms $\mathbf{I}_{1,n}$, $\mathbf{I}_{2,n}$ we used H\"older's inequality, \eqref{eq-1:lemtech2} and \eqref{eq0:lemtech2}. On the other hand, from \eqref{eq1:lemtech2}, as $n\to +\infty$, we have
\begin{equation}\label{eq4:lemtech2}
\int_{\Omega} \nabla\left[u_{\Omega_n} (1-\varphi_n)\right]\cdot \nabla \varphi\, dx = \int_{\Omega} \nabla u \cdot \nabla \varphi\, dx+o(1).
\end{equation}
Hence, comparing \eqref{eq3:lemtech2} and \eqref{eq4:lemtech2} and taking the limit as $n\to +\infty$ we obtain that
$$ \int_{\Omega} \nabla u \cdot \nabla \varphi\, dx =  \int_{\Omega} \varphi\, dx,$$
and by uniqueness of the weak solution to \eqref{eq:mixbvprobquasiopen} we conclude that $u=u_\Omega$. Finally, from \eqref{eq2:lemtech2} we readily deduce that
\begin{equation}\label{eq5:lemtech2}
 u_{\Omega_n} (1-\varphi_n) \to u_\Omega,
\end{equation} 
 strongly in $L^2(\C_\omega)$, as $n\to +\infty$.
Summing up, from \eqref{eq-3:lemtech2}, \eqref{eq-2:lemtech2} and \eqref{eq5:lemtech2} we infer that $u_{\Omega_n}\to u_\Omega$ strongly in $L^2(\C_\omega)$, as $n\to +\infty$, namely $\Omega_n$ $\gamma$-converges to $\Omega$, and by construction $(\Omega_n)_n$ has the desired properties. The proof is complete.
\end{proof}

%As a first result we show that we can always find a minimizing sequence for $\O_c(\C_\omega)$ (or $\O_c(\C_\omega^+)$ made by convex domains with respect to the variable $x_N$.

\begin{lemma}\label{lem1:sect3}
Let $c>0$ and let $\omega\subset\R^{N-1}$ be a Lipschitz bounded domain. For any sequence $(\e_n)_n$ of positive real numbers such that $\e_n\to 0^+$, as $n\to +\infty$, there exist a sequence $(\Omega_n)_n$ of open subsets of $\C_\omega$ and a sequence $(\delta_n)_n$ of positive real numbers such that $\delta_n\to 0^+$, as $n\to +\infty$, and for any $n\in \N$ it holds:
\begin{itemize}
\item[(i)] $|\Omega_n|\leq c+\delta_n$;
\item[(ii)]  $\Omega_n$ is convex with respect to the variable $x_N$;
\item[(iii)]  $\Omega_n$ is symmetric with respect to the hyperplane $\{x_N=0\}$;
\item[(iv)] $\O_{c+\delta_n}(\C_\omega)\leq \E(\Omega_n; \C_\omega) \leq \O_c(\C_\omega)+2\e_n$.
\end{itemize}
Moreover, for all $n\in \N$, denoting by $u_{\Omega_n}$ the energy function of $\Omega_n$, it holds
\begin{itemize}
\item[(v)] $u_{\Omega_n}$ is even with respect to $x_N$, a.e. in $\C_\omega$;
\item[(vi)] $\frac{\partial u_{\Omega_n}}{\partial x_N} (x^\prime, x_N) \leq 0 $ for a.e. $(x^\prime, x_N)\in \C_\omega\cap\{x_N> 0\}$.
\end{itemize}
%The same result holds for $\O_{c}(\C^+_{\omega})$, apart from properties (iii) and (iv).
\end{lemma}
\begin{proof}
Let us fix $c > 0$ and $\omega$ as in the statement.  By definition of $\O_{c}(\C_{\omega})$, for any $\e>0$ there exists a quasi-open set $\widetilde\Omega_\e\in\A_{\omega, c}$ such that
$$\O_{c}(\C_{\omega})\leq \E(\widetilde\Omega_\e; \C_{\omega}) \leq \O_{c}(\C_{\omega}) + \e.$$
By Lemma \ref{lem:tech2}, fixing a sequence $(\eta_m)_m$ of positive real numbers such that $\eta_m\to 0^+$, as $m\to +\infty$, we find a sequence $(\widetilde\Omega_{\e,m})_m$ of open subsets of $\C_\omega$ such that 
\begin{equation}\label{eq1:techlem3}
|\widetilde\Omega_{\e,m}| \leq c+\eta_m, \ \ \forall m\in \N
\end{equation}
and
\begin{equation}\label{eq2:techlem3}
u_{\widetilde\Omega_{\e,m}}\to u_{\widetilde\Omega_\e},\text{strongly in $L^2(\C_\omega)$, as $m\to +\infty$.}
\end{equation}
Then, as $u_{\widetilde\Omega_{\e,m}}\in H^1_0(\widetilde\Omega_{\e,m}; \C_\omega)$, $u_{\widetilde\Omega_{\e}}\in H^1_0(\widetilde\Omega_{\e}; \C_\omega)$, taking into account \eqref{eq1:techlem3}, that $|\widetilde\Omega_\e|\leq c$, and Proposition \ref{prop1:sect2}, from \eqref{eq2:techlem3} and H\"older's inequality we deduce that $u_{\widetilde\Omega_{\e,m}}\to u_{\widetilde\Omega_\e}$ strongly in $L^1(\C_\omega)$,  as $m\to +\infty$. Hence, since $\E(\widetilde\Omega_{\e,m}; \C_{\omega})=-\frac{1}{2} \int_{\C_\omega} u_{\widetilde\Omega_{\e,m}} \, dx$, we get that
$$\lim_{m\to +\infty} \E(\widetilde\Omega_{\e,m}; \C_{\omega})=\E(\widetilde\Omega_\e; \C_{\omega}).$$
In particular, for any $\e>0$ we find $\bar m=\bar m(\e)\in \N^+$ such that
\begin{equation}\label{eq3:techlem3}
 \O_{c}(\C_{\omega}) -\e  \leq \E(\widetilde\Omega_{\e,m}; \C_{\omega})  \leq \O_{c}(\C_{\omega}) + 2\e, \ \ \forall m\geq \bar m(\e).
\end{equation}
Therefore, given a sequence of positive real numbers $(\e_n)_n$ such that $\e_n\to 0^+$, as $n\to +\infty$, and setting  for every $n\in\N$
$$
m_{n}=\max\left\{\bar m(\e_{n}),\left[\e_n^{-1}\right]+1\right\}\,,~~
\widetilde\Omega_n:=\widetilde\Omega_{\e_n,m_{n}}, \ \ \ \delta_n:=\eta_{m_{n}},
$$
then $m_{n}\to+\infty$ and $\delta_{n}\to 0$ as $n\to+\infty$ and, by construction, \eqref{eq1:techlem3} and \eqref{eq3:techlem3}, $(\widetilde\Omega_n)_n$ is a  sequence of open sets in $\C_\omega$ such that $|\widetilde\Omega_n|\leq c+\delta_n$ for all $n\in\N$ and 
$$
\lim_{n\to+\infty} \E(\widetilde\Omega_n; \C_\omega)= \O_{c}(\C_{\omega}).
$$
Now, setting 
$$
\Omega_n:=\widetilde\Omega_n^{\sharp},\ n\in\N, 
$$
where $\widetilde\Omega_n^{\sharp}$ is the Steiner symmetrization of $\widetilde\Omega_n$ with respect to the hyperplane $\{x_N=0\}$, then, by construction, \eqref{eq:steinerpreservesvolume} and thanks to Lemma \ref{opensharp}, $\Omega_n$ is an open subset of $\C_\omega$ and satisfies (i)--(iii). Moreover, by definition, Lemma \ref{lem1:sect2}-(ii) and \eqref{eq3:techlem3} we infer that
$$
\E(\Omega_n; \C_\omega)\leq \E(\widetilde\Omega_n;\C_\omega)\leq \O_c(\C_\omega) + 2\e_n.
$$
On the other hand, as $\Omega_n\in \A_{\omega, c+\delta_n}$, we have
$$
\O_{c+\delta_n}(\C_\omega) \leq \E(\Omega_n; \C_\omega),
$$
and thus also (iv) is fulfilled.
\medskip

%Let us consider the Steiner symmetrization of $u_{\widetilde\Omega_n}^\sharp$. By Proposition......for all $n\in\N$ it holds $u_{\widetilde\Omega_n}^\sharp \in H_0^1(\S(\widetilde\Omega_n); \C_\omega)$ and by Proposition .... we infer that
%\beq\label{eq1:lem1sect3}
%J_{\S(\widetilde\Omega_n)}(u_{\widetilde\Omega_n}^\sharp)\leq J_{\widetilde\Omega_n}%(u_{\widetilde\Omega_n}).
%\eeq
%Moreover, as $\E(\S(\widetilde\Omega_n); \C_{\omega})\leq J_{\S(\widetilde\Omega_n)}(u_{\widetilde\Omega_n}^\sharp)$, from \eqref{eq1:lem1sect3}, and since $%\S(\widetilde\Omega_n)\in \A_{\omega,c}$, we obtain that
%\beq\label{eq1tris:lem1sect3}
% \O_{c+\delta_n}(\C_{\omega})\leq\E(\S(\widetilde\Omega_n); \C_{\omega})\leq %\E(\widetilde\Omega_n; \C_{\omega}).
 %\eeq
 %From this we readily infer that $(\S(\widetilde\Omega_n))_n$ is a minimizing sequence for $\O_c(\C_\omega)$ in $\A_{\omega,c}$, and by definition of Steiner symmetrization we check that $\S(\widetilde\Omega_n)$ satisfies (ii) and (iii) for all $n\in\N$.\\ 
 
Finally, for (v) and (vi), let us consider the Steiner symmetrization $u_{\Omega_n}^\sharp$ of $u_{\Omega_n}$. Since $\widetilde\Omega_n^{\sharp}=\Omega_n$, from Lemma \ref{lem1:sect2}-(ii) we infer that $u_{\Omega_n}^\sharp\in H_0^1(\Omega_n; \C_\omega)$, and by \eqref{eq:Lp-normpreservedSteiner} and Proposition \ref{prop2:sect2} we obtain
$$
J_{\Omega_n}(u_{\Omega_n}^\sharp)\leq J_{\Omega_n}(u_{\Omega_n})=\E(\Omega_n; \C_{\omega}).
$$
Hence, by uniqueness of the minimizer for $J_{\Omega_n}$ in $H^{1}_{0}(\Omega_{n};\C_{\omega})$, we conclude that 
$$
u_{\Omega_n}^\sharp=u_{\Omega_n},
$$
and thus, by the definition of $u_{\Omega_n}^\sharp$, we infer that $u_{\Omega_n}$ is even with respect to the variable $x_N$, a.e. in $\C_\omega$ and
$$\frac{\partial u_{\Omega_n}}{\partial x_N} (x^\prime, x_N) \leq 0  \ \ \hbox{for a.e. $(x^\prime, x_N)\in \C_\omega\cap\{x_N>0\}$}.$$
The proof is complete.
\end{proof}

\section{Existence of minimizers for the torsional energy in the cylinder}\label{S:minimizer}

In this section we prove the existence of a quasi-open set in the cylinder which minimizes the torsional energy under a volume constraint. Namely, we prove the existence of $\Omega^*$ in the class $\A_{\omega,c}$ (see \eqref{eq:defclassadmissiblecylinder}) such that $\E(\Omega^*; \C_\omega)=\O_c(\C_\omega)$. This is the content of the following.

\begin{theorem}\label{teo:mainteoexistmin}
Let $c>0$ and let $\omega\subset\R^{N-1}$ be a bounded domain with Lipschitz boundary. Then there exists a quasi-open set $\Omega^*\in \A_{\omega,c}$ which is a minimizer of $\E(\cdot; \C_\omega)$ in $\A_{\omega,c}$.
\end{theorem}

\begin{proof}
Let us fix $c > 0$ and $\omega$ as in the statement. Let $(\e_n)_n$ be a sequence of positive real numbers such that $\e_n\to 0^+$, as $n\to +\infty$, and let $(\Omega_n)_n$ be the sequence of open subsets of $\C_\omega$ provided by Lemma \ref{lem1:sect3} and let $(u_{\Omega_n})_n $ be the corresponding sequence of energy functions. 
For simplicity we set $u_n:=u_{\Omega_n}$, for all $n\in\N$, and observe that by Lemma \ref{lem1:sect3}-(i) and Proposition \ref{prop1:sect2} (ii)-(iii) it follows that $(u_n)_n$ is a bounded sequence in $H^1(\C_\omega)$. In particular, since $(\|u_n\|^2_{L^2(\C_\omega)})_n$ is a bounded sequence in $\R^+\cup\{0\}$, then, up to a subsequence (still denoted for convenience by $(u_n)_n$), it holds that

$$
\|u_n\|^2_{L^2(\C_\omega)} \to \lambda,\ \ \hbox{for $n\to +\infty$},
$$
for some $\lambda\geq 0$. We claim that $\lambda>0$.

Indeed, if $\lambda=0$, then,  by H\"older's inequality, exploiting that $u_n \in H_0^1(\Omega_n; \C_\omega)$, taking into account Lemma \ref{lem1:sect3}-(i), we get that
$$
\|u_n\|_{L^1(\C_\omega)}\to 0, \ \hbox{as $n\to +\infty$}.
$$
and thus, by \eqref{def:torsionquasiopen}, we have $\E(\Omega_n; \C_\omega)\to 0$, as $n\to +\infty$.
On the other hand, by Lemma \ref{lem1:sect3}-(iv) we have
\beq\label{eq-1:teo1sect4}
\E(\Omega_n; \C_\omega) \leq \O_c(\C_\omega)+2\e_n,
\eeq
and taking the limit as $n\to+\infty$ we conclude that $\O_c(\C_\omega)\geq 0$, contradicting $\O_c(\C_\omega)<0$ (see Remark \ref{remark:negativityenergy}).

Therefore, up to taking $\tilde u_n:=\lambda \frac{u_n}{\|u_n\|_{L^2(\C_\omega)}}$, we can assume without loss of generality that 
\beq\label{eq0:teo1sect4}
\|u_n\|_{L^2(\C_\omega)}=\lambda>0 \ \ \ \forall n\in \N.
\eeq
Now, since $(u_n)_n$ is a bounded sequence in $H^1(\C_\omega)$ such that \eqref{eq0:teo1sect4} holds then it satisfies the hypotheses of the concentration-compactness principle, as formulated in \cite[Lemma III.1]{LIO}, with minor adjustments due to the fact that our working functional space is $H^1(\C_\omega)$ instead of $H^1(\R^N)$. However, for our purposes, it is not sufficient to apply directly \cite[Lemma III.1]{LIO} but we need to combine the construction performed in the proof of P. L. Lions together with the properties of the sequence $(u_n)$ stated in Lemma \ref{lem1:sect3}. 

To this end, following the proof of \cite[Lemma I.1, Lemma III.1]{LIO}  we consider the sequence of concentration functions $(Q_n)_n$ associated to $(u_n)_n$, where, for all $n\in\N$, $Q_n\colon\R^+\cup\{0\}\to \R^+\cup\{0\}$ is the function defined by
$$Q_n(R):=\sup_{y\in\C_\omega} \int_{B_R(y) \cap \C_\omega} u_n^2\, dx.$$
We notice that, since $u_n\in L^2(\C_\omega)$, the supremum in the definition of $Q_n(R)$ is actually a maximum over $\overline{\C_\omega}$, namely, for all $n\in \N$, $R\geq 0$ there exists $y_{n,R}
\in\overline{\C_\omega}$ such that 
\beq\label{eq:remarkQn}
Q_n(R)=\int_{B_R(y_{n,R}) \cap \C_\omega} u_n^2\, dx.
\eeq
By construction we have that $(Q_n)_n$ is a sequence of non-decreasing, non-negative and uniformly bounded functions defined in $\R^+\cup\{0\}$. Then, by Helly's selection theorem there exists a subsequence $(n_k)_k$ in $\N$ such that $(Q_{n_k})_k$ converges pointwise in $\R^+\cup\{0\}$. Defining the function $Q:\R^+\cup\{0\}\to \R^+\cup\{0\}$ by
$$
Q(R):=\lim_{k\to +\infty} Q_{n_k}(R),
$$
we readily check that $Q$ is non-decreasing, non-negative and bounded. Hence there exists 
\beq\label{eq:defalpha}
\alpha:=\lim_{R\to+\infty} Q(R), \ \hbox{and $\alpha\in[0,\lambda]$.}
\eeq
Then, only one of the following three alternatives can occur: 
\begin{itemize}
\item[(i)]  $\alpha=\lambda$;
\item[(ii)]  $\alpha=0$;
\item[(iii)]  $0<\alpha<\lambda$.
\end{itemize}

%As seen in the proof of  \cite[Lemma I.1]{LIO}, the cases $\alpha=\lambda$ and $\alpha=0$, imply, respectively, that
%\begin{itemize}
%\item[(i)] there exists a sequence $(y_{n_k})_k$ of points in $\overline{\C_\omega}$  such that $$ \forall \e>0 \ \exists R>0 \ \ \hbox{such that} \  \ \int_{B_R(y_{n_k})\cap \C_\omega} u_{n_k}^2 \ dx \geq \lambda - \eps \ \ \forall k\in \N;\vspace{-4pt}$$
%\item[(ii)] for all $R>0$ it holds $\displaystyle \lim_{k\to +\infty} \sup_{y \in \C_\omega} \int_{B_R(y)\cap \C_\omega} u_{n_k}^2 \ dx=0$;\\[2pt]
%\item[(iii)] there exists $\alpha\in]0,\lambda[$ such that for all $\e>0$, there exist $R_1=R_1(\e)>0$, $k_0=k_0(\e)\in\N^+$,
%a sequence of positive real numbers $(R_k)_{k}$ such that $R_k\to +\infty$, as $k\to +\infty$, and a sequence $(y_k)_k$ of points in $\overline{C_\omega}$ such that, fixing $\xi, \varphi \in C_c^\infty(\R^N,[0,1])$ functions such that $\xi\equiv 1$, $\varphi\equiv 0$ on $\overline{B_1(0)}$ and $\xi\equiv 0$, $\varphi\equiv 1$  in $(\overline{B_2(0)})^\complement$, and setting
%\end{itemize}

In the next two steps we prove that the only admissible case is (i).
\medskip
 
\noindent
\textbf{Step 1:} the ``vanishing'' case (ii) cannot happen.
\medskip
 
The proof is essentially the same of  \cite[Theorem 6.8-Step 1]{IPW}. For the sake of completeness we provide it.
 Assume by contradiction that $\alpha=0$. By definition, and since $Q$ is non-decreasing, it follows that that for all $R>0$ it holds
$$
\lim_{k\to +\infty} \sup_{y \in \C_\omega} \int_{B_R(y)\cap \C_\omega} u_{n_k}^2 \ dx=0.
$$
Then, since $(u_{n_k})_k$ is a bounded sequence in $H^1(\C_\omega)$, by \cite[Lemma 1.21]{WI} (which holds also in $H^1(\C_\omega)$ with the same proof), we infer that $u_{n_k}\to 0$ in  $L^p(\C_\omega)$, for any $2<p<\frac{2N}{N-2}$, if $N\geq 3$, or for any $p>2$ if $N=2$. Hence, since $u_{n_k}\in H_0^1(\Omega_{n_k}; \C_\omega)$ and $|\Omega_{n_k}|\leq c+\delta_n$, for all $k\in\N$, then by H\"older's inequality we get that $u_{n_k}\to 0$ in $L^2(\C_\omega)$, as $k\to +\infty$, contradicting \eqref{eq0:teo1sect4}.
\medskip

\noindent\textbf{Step 2:} the ``dichotomy" case (iii) cannot happen.
\medskip

Assume by contradiction that (iii) holds. Let  $\xi, \theta \in C_c^\infty(\R^N,[0,1])$ be two cut-off functions such that $\xi\equiv 1$, $\theta\equiv 0$ on $\overline{B_1(0)}$ and $\xi\equiv 0$, $\theta\equiv 1$  in $(\overline{B_2(0)})^\complement$, and fix a positive real number $L$ such that 
\begin{equation}
\label{un-bounded}
\sup_{n\in \N}\|u_n\|_{H^1(\C_\omega)}\le L.
\end{equation}
For $\beta>0$ we define the rescaled functions $\xi_\beta, \theta_\beta \in C_c^\infty(\R^N,[0,1])$ by 
$$\xi_\beta(x):=\xi \left(\frac{x}{\beta}\right), \ \ \theta_\beta(x):=\theta \left(\frac{x}{\beta}\right),\ \  x\in \R^N.$$
Following verbatim the proof of \cite[Lemma III.1]{LIO} we have that, for all sufficiently small $\e>0$ there exist $R_1=R_1(\e,L,\theta,\xi)>0$, $C=C(L)>0$ and $k_0=k_0(\e)\in \N^+$ such that for all $v\in H^1(\C_\omega)$, with $\|v\|_{H^1(\C_\omega)}\leq L$, and for all $R\geq R_1$ it holds
$$\left|\int_{\C_\omega} |\nabla (\xi_R v) |^2\, dx - \int_{\C_\omega} \xi_R^2 |\nabla v|^2\, dx \right|\leq \e, \ \ \left|\int_{\C_\omega} |\nabla (\theta_R v) |^2\, dx - \int_{\C_\omega} \theta_R^2 |\nabla v|^2\, dx \right|\leq \e$$
and for any $k\geq k_0$
\begin{equation}\label{eq2:teo1sect4}
\alpha-2\e\leq Q_{n_k}(R_1)\leq \alpha-\e.
\eeq
As observed in \eqref{eq:remarkQn}, for all $k\geq k_0$, we find $y_k:=y_{n_k,R_1}\in\overline{\C_\omega}$ such that
 \begin{equation}\label{eq:defyk:teo1sect4}
Q_{n_k}(R_1)= \int_{B_{R_1}(y_k)\cap \C_\omega} u_{n_k}^2\, dx.
\eeq
Since $u_{n_k}$ is even with respect to the variable $x_N$ (see Lemma \ref{lem1:sect3}-(v)), writing $y_k=(y_k^\prime, (y_k)_N)$, where $y_k^\prime\in \overline\omega$, and $(y_k)_N \in \R$ denotes the $N$-th component of $y_k$, we can assume without loss of generality that for all $k\geq k_0$ it holds
\beq\label{eq:wlogsignykN}
 (y_k)_N\geq 0.
 \eeq
Now, setting for $k\geq k_0$
\begin{equation}\label{eq3:teo1sect4}
u_{1,k}(x):=\xi_{R_1}\left(x-y_k\right)u_{n_k}(x), \ \ x\in \C_\omega,
\eeq
by \eqref{un-bounded}, the sequence $(u_{1,k})_{k\geq k_0}$ is bounded in $H^1(\C_\omega)$, and, thanks to \eqref{eq2:teo1sect4}, for all $k\geq k_0$ it holds
\begin{equation}\label{eq4:teo1sect4}
\left| \int_{\C_\omega} u_{1,k}^2 \ dx -\alpha \right| \leq 2\e.
 \eeq
 Moreover, from the definition of $\alpha$ (see \eqref{eq:defalpha}) we find $k_1=k_1(\e)\in \N^+$ (without loss of generality we can assume $k_1\geq k_0$) and a sequence of positive real numbers $(R_k)_{k\geq k_1}$ such that $R_k\to +\infty$, as $k\to +\infty$, and 
 \begin{equation}\label{eq5:teo1sect4}
  Q_{n_k}(2R_k)\leq \alpha+2\e \ \ \forall k\geq k_1.
  \eeq
For $k\geq k_1$ we set
\begin{equation}\label{eq6:teo1sect4}
u_{2,k}(x):=\theta_{R_k}\left(x-y_k\right) u_{n_k}(x),  \ x\in \C_\omega.
\eeq
By construction $(u_{2,k})_{k\geq k_1}$ is a bounded sequence in $H^1(\C_\omega)$ and it holds
$$\textrm{dist}(\textrm{supp}(u_{1,k}), \textrm{supp}(u_{2,k})) \to +\infty,\ \hbox{as $k\to+\infty$.}$$
In addition, denoting by $A_{R_1,R_k}(y_k):=\{x\in \R^N; \ R_1\leq |x-y_k|\leq 2R_k\}$ the annulus of radii $R_1$, $R_k$, centred at $y_k$, then from \eqref{eq2:teo1sect4}--\eqref{eq3:teo1sect4} and \eqref{eq5:teo1sect4}--\eqref{eq6:teo1sect4}, for all $k\geq k_1$, we infer that
\begin{equation}\label{eq7:teo1sect4}
\int_{\C_\omega} |u_{n_k} - u_{1,k}-u_{2,k}|^2\, dx \leq \int_{\C_\omega\cap A_{R_1,R_k}(y_k)} u_{n_k}^2\, dx \leq Q_{n_k}(2R_k)- Q_{n_k}(R_1) \leq 4\e.
\eeq
From \eqref{eq4:teo1sect4} and \eqref{eq7:teo1sect4} we have
\begin{equation}\label{eq8:teo1sect4}
\lim_{k\to +\infty}\int_{\C_\omega} |u_{n_k} - u_{1,k}-u_{2,k}|^2\, dx=0,\ \  \lim_{k\to +\infty} \int_{\C_\omega} u_{1,k}^2 \ dx =\alpha.
\eeq
Then, since the supports of $u_{1,k}$, $u_{2,k}$ are disjoint for large $k$, from \eqref{eq8:teo1sect4} and (iii), we check that
\beq\label{eq9:teo1sect4}
\lim_{k\to +\infty} \int_{\C_\omega} u_{2,k}^2 \ dx =\lambda-\alpha>0.
\eeq

In order to conclude the proof of \textit{Step 2}, we need a couple of intermediate steps. We begin by analyzing the asymptotic behaviour of the sequence of points $(y_k)_{k\geq k_0}$ in $\overline{\C_\omega}$ appearing in the definition of $u_{1,k}$ (and also in $u_{2,k}$). More precisely, the following property holds true.
\medskip

\textbf{Claim I}: let $\e>0$ such that $\alpha -2\e>0$ and let $(y_k)_{k\geq k_0}$ be a sequence in $\overline{\C_\omega}$ satisfying \eqref{eq:defyk:teo1sect4}. Then  $$\sup_{k\geq k_0} |y_k| < +\infty.$$

Assume by contradiction that the thesis is false. Then, there exist $\e>0$ such that $\alpha -2\e>0$, and a subsequence, still indexed by $k$ for convenience, let us say defined for $k\geq k_2$, with $k_2\geq k_0$, such that $$\lim_{k\to +\infty} |y_k|=+\infty.$$
Since $\C_\omega=\omega\times\R$ and $\omega$ is a bounded domain in $\R^{N-1}$, then, from \eqref{eq:wlogsignykN} we infer that
\beq\label{eq10:teo1sect4}
\lim_{k\to +\infty} (y_k)_N=+\infty.
\eeq
Moreover, let $m_k \in \N$, $r_k\in\R^+\cup\{0\}$ be such that $0\leq r_k<2R_1$ and
$$(y_k)_N=2m_k R_1+r_k,$$
where $R_1=R_1(\e, L, \xi, \theta)$ is the positive real number appearing in \eqref{eq:defyk:teo1sect4}. 
By definition $r_k$ is uniformly bounded, hence from \eqref{eq10:teo1sect4} it follows that, as $k\to +\infty$, it holds
\beq\label{eq10m:teo1sect4}
 m_k\to +\infty.
\eeq
Now, observing that $B_{2R_1}(y_k) \cap \C_\omega \subset \omega \times [(y_k)_N-2R_1, (y_k)_N+2R_1]$, and setting $$\tau_1:=\alpha-2\e>0,$$ 
then, by definition of $u_{1,k}$ (see \eqref{eq3:teo1sect4}) and thanks to \eqref{eq4:teo1sect4}, we infer that for all $k\geq k_2$
\beq\label{eq11:teo1sect4}
\int_{\omega\times [(y_k)_N-2R_1, (y_k)_N+2R_1]} u_{1,k}^2\,dx \geq \tau_1.
\eeq
On the other hand, recalling that by Lemma \ref{lem1:sect3} (v)-(vi) the function $u_{n_k}$ is monotone decreasing with respect to $x_N$, a.e. in $\C_\omega\cap\{x_N>0\}$, even with respect to $x_N$ a.e. in $\C_\omega$, we deduce that for all sufficiently large $k$
%Then, thanks to the monotonicity property of $u_{n_k}$ stated in Lemma \ref{lem1:sect3} (vi), we deduce that
\beq\label{eq12:teo1sect4}
\begin{split}
\displaystyle\frac{1}{2}\int_{\C_\omega} u_{n_k}^2\, dx&=\displaystyle\int_{\C_\omega\cap\{x_N\geq 0\}} u_{n_k}^2\, dx\geq\displaystyle \int_{\omega\times [0, (y_k)_N+R_1]} u_{n_k}^2\, dx\\
&\geq\displaystyle \sum_{j=0}^{m_k-1}\int_{\omega\times [ (y_k)_N-2jR_1-R_{1}, (y_k)_N-2jR_1+R_{1}]} u_{n_k}^2\, dx\\
&\geq\displaystyle m_k\int_{\omega\times [ (y_k)_N-R_1, (y_k)_N+R_1]} u_{n_k}^2\, dx\geq\displaystyle m_k\int_{\omega\times [ (y_k)_N-R_1, (y_k)_N+R_1]} u_{k,1}^2\, dx\,,
\end{split}
\eeq
where in the last inequality we used that $u_{n_k}\geq u_{1,k}$ in $\C_\omega$, as it readily follows by definition.
Finally, combining \eqref{eq11:teo1sect4} and \eqref{eq12:teo1sect4} we conclude that for all sufficiently large $k$ it holds that
\beq\label{eq13:teo1sect4}
\frac{1}{2}\int_{\C_\omega} u_{n_k}^2\, dx \geq m_k \tau_1,
\eeq
and then, thanks to \eqref{eq10m:teo1sect4}, taking the limit, as $k\to +\infty$, we get that $\|u_{n_k}\|_{L^2(\C_\omega)}\to +\infty$, contradicting \eqref{eq0:teo1sect4}. The proof of \textit{Claim I} is complete.
\medskip

As a second intermediate step we analyze $(u_{2,k})_k$ and show that it cannot exhibit vanishing. More precisely, the following fact holds:
\medskip

\textbf{Claim II}: there exist $\overline R>0$, $\tau_2>0$, $\bar k\in \N^+$ such that for all $k\geq \bar k$.
\beq\label{eq14:teo1sect4}
\sup_{y\in \C_\omega} \int_{B_{\overline R}(y)\cap \C_\omega} u_{2,k}^2\, dx> \tau_2.\\
\eeq
If \eqref{eq14:teo1sect4} does not hold then, for all $\overline R>0$ we find a subsequence $(k_h)_h$ in $\N^+$, such that $k_h\to +\infty$, as $h\to +\infty$, and
$$
\lim_{h\to +\infty}\sup_{y\in \C_\omega} \int_{B_{\overline R}(y)\cap \C_\omega} u_{2,k_h}^2\, dx= 0.
$$
Since $(u_{2,k})_k$ is a bounded sequence in $H^1(\C_\omega)$, then, from \eqref{eq11:teo1sect4} and \cite[Lemma 1.21]{WI}, we get that $u_{2,k_h}\to 0$ strongly in  $L^p(\C_\omega)$, as $h\to +\infty$, for any $2<p<\frac{2N}{N-2}$, if $N\geq 3$, or for any $p>2$ if $N=2$. Observing that by construction it holds $u_{2,k_h}\in H_0^1(\Omega_{{n_k}_h}; \C_\omega)$, and recalling that $|\Omega_{{n_k}_h}|\leq c+\delta_{{n_k}_h}$, then, arguing as in the proof of \textit{Step 1} we conclude that $u_{2,k_h}\to 0$ strongly in $L^2(\C_\omega)$, as $h\to +\infty$, contradicting \eqref{eq9:teo1sect4}. The proof of \textit{Claim II} is complete.
\medskip

We now conclude the proof of \textit{Step 2}. In view of \textit{Claim II}, for all $k\geq \bar k$, we find $z_k \in \C_\omega$ such that
$$
\int_{B_{\overline R}(z_k)\cap \C_\omega} u_{2,k}^2\, dx\geq \frac{\tau_2}{2}.
$$
In particular, writing $z_k =(z_k^\prime,(z_k)_N)$, with $z_k^\prime\in \omega$, since the set $B_{\overline R}(z_k)\cap \C_\omega$ is contained in the bounded cylinder $\omega\times[(z_k)_N-\overline R, (z_k)_N+\overline R]$, then, for all $k\geq \bar k$ it holds
\beq\label{eq16:teo1sect4}
\int_{\omega\times[(z_k)_N-\overline R, (z_k)_N+\overline R]} u_{2,k}^2\, dx\geq \frac{\tau}{2}.
\eeq
Moreover, by definition of $u_{2,k}$ we have $u_{2,k}\equiv 0$ in $B_{R_k}(y_k)\cap \C_\omega$  (see \eqref{eq6:teo1sect4}), and thus from \eqref{eq11:teo1sect4}, we readily deduce that for all $k\geq \bar k$ 
 $$(\omega\times[(z_k)_N-\overline R, (z_k)_N+\overline R]) \cap (B_{R_k}(y_k))^\complement\neq \varnothing.$$ 
 In particular, thanks to \textit{Claim I}, taking into account \eqref{eq:wlogsignykN}, since $\overline R$ is independent of $k$ and $R_k\to +\infty$, as $k\to +\infty$, we deduce that
$$\lim_{k\to +\infty}(z_k)_N=+\infty.$$
Then, as in the proof of \textit{Claim I}, writing
$$(z_k)_N=2m_k^\prime \overline R+r_k^\prime,$$
for some $m_k^\prime \in \N$, $r_k^\prime\in\R^+\cup\{0\}$ such that $0\leq r_k^\prime<2\overline R$, we infer that $m_k^\prime \to +\infty$, as $k\to +\infty$, and arguing as in \eqref{eq12:teo1sect4}--\eqref{eq13:teo1sect4}, exploiting \eqref{eq16:teo1sect4}, we conclude that
\beq\label{eq:17:teo1sect4}
\frac{1}{2}\int_{\C_\omega} u_{n_k}^2\, dx \geq m_k^\prime\int_{\omega\times [ (z_k)_N-\overline R, (z_k)_N+\overline R]} u_{k,2}^2\, dx \geq m_k^\prime\frac{\tau_2}{2}.
\eeq
Finally, taking the limit as $k\to +\infty$ in \eqref{eq:17:teo1sect4} we get that $\|u_{n_k}\|_{L^2(\C_\omega)}\to +\infty$, but this contradicts \eqref{eq0:teo1sect4}. The proof of \textit{Step 2} is complete.
\medskip

Thanks to \textit{Step 1} and \textit{Step 2} we know that the only possibility is the ``compactness" case (i), i.e. $\alpha=\lambda$. In the next step we prove that this implies the pre-compactness in $L^2(\C_\omega)$ of the sequence $(u_{n_k})_k$. The proof is similar to that of \cite[Theorem 6.8-Step 5]{IPW}, but some adjustments are needed.\medskip

\noindent\textbf{Step 3:} The sequence $(u_{n_k})_k$ admits a subsequence which strongly converges in $L^2(\C_\omega)$.
\medskip

From (i) we know that there exists a sequence $(y_{n_k})_k$ of points in $\overline{\C_\omega}$ satisfying the following property:
\beq\label{eq:18:teo1sect4}
 \forall \e>0 \ \exists R>0 \ \hbox{s.t.} \ \forall k\in \N\ \ \int_{B_R(y_{n_k})\cap \C_\omega} u_{n_k}^2 \ dx \geq \lambda - \eps.
\eeq
As a first intermediate step we show the following property.
\medskip

\noindent\textbf{Claim III:} The sequence $(y_{n_k})_k$ is bounded.
\medskip

Indeed, writing $y_{n_k}=(y_{n_k}^\prime, (y_{n_k})_N)$, with $y_{n_k}^\prime\in\overline\omega$, $(y_{n_k})_N\in\R$, as $\omega$ is a bounded domain, it is sufficient to prove that $((y_{n_k})_N)_k$ is a bounded sequence. For this, thanks to \eqref{eq:18:teo1sect4}, fixing $\e>0$ such that $\lambda - \eps>\frac{\lambda}{2}$, we find $R=R(\e)>0$ independent of $k$, such that for all $k$ we have
\beq\label{eq19:teo1sect4}
\int_{\omega\times[R-(y_{n_k})_N, R+(y_{n_k})_N]} u_{n_k}^2 \ dx \geq \frac{\lambda}{2}.
\eeq

Then, assuming by contradiction that, up to a subsequence, $(y_{n_k})_N\to +\infty$, as $k\to +\infty$, and arguing as in the proof of \textit{Claim I}, taking into account \eqref{eq19:teo1sect4}, we infer that for all sufficiently large $k$ it holds 
\beq\label{eq:20:teo1sect4}
\frac{1}{2}\int_{\C_\omega} u_{n_k}^2\, dx \geq m_k^{\prime\prime}\int_{\omega\times[R-(y_{n_k})_N, R+(y_{n_k})_N]} u_{n_k}^2 \ dx \geq m_k^{\prime\prime}\frac{\lambda}{2},
\eeq
where  $(m_k^{\prime\prime})_k$ is a sequence of natural numbers such that $m_k^{\prime\prime}\to +\infty$, as $k\to +\infty$. Hence, taking the limit as $k\to +\infty$ in \eqref{eq:20:teo1sect4} we contradict \eqref{eq0:teo1sect4}. The proof of \textit{Claim III} is complete.
\medskip

Next, as a second intermediate step, we prove the following.
\medskip

\textbf{Claim IV}: For any $y_0 \in \C_\omega$ it holds 
\beq\label{eq21:teo1sect4}
\lim_{R\to+\infty} \sup_{k \in \N} \int_{B_R^\complement(y_0)\cap\C_\omega} u_{n_k}^2 \ dx=0.
\eeq

To prove the claim we argue by contradiction. Assume that \eqref{eq21:teo1sect4} is not true, then there exist $y_0\in \C_\omega $, $\eps^\prime>0$, a sequence $(R_h)_h$ of positive real numbers such that $R_h\to+\infty$, as $h\to +\infty$, and a subsequence $(n_{k_h})_h$ such that for all $h\in \N$ it holds
\beq\label{eq22:teo1sect4}
\int_{B_{R_h}^\complement(y_0)\cap\C_\omega} u_{n_{k_h}}^2 \ dx\geq \frac{\eps^\prime}{2}.
\eeq
On the other hand, taking $\e=\frac{\eps^\prime}{4}$ in \eqref{eq:18:teo1sect4} we find $R^\prime=R^\prime(\e)>0$ such that for all $k \in \N$
\beq\label{eq23:teo1sect4}
 \int_{B_{R^\prime}(y_k)\cap\C_\omega} u_{n_{k}}^2 \ dx\geq \lambda-\frac{\eps^\prime}{4}.
\eeq
Thanks to \textit{Claim III} we infer that there exists $R^{\prime\prime}>0$ independent of $k$, such that $B_{R^\prime}(y_{n_k})\subset B_{R^{\prime\prime}}(y_0)$ for all $k\in \N$. Hence, from \eqref{eq23:teo1sect4}, we obtain that for all $k\in \N$ 
\beq\label{eq24:teo1sect4}
\int_{B_{R^{\prime\prime}(y_0)}\cap\C_\omega} u_{n_{k}}^2 \ dx \geq \int_{B_{R^\prime}(y_{n_k})\cap\C_\omega} u_{n_{k}}^2 \ dx\geq \lambda-\frac{\eps^\prime}{4}.
\eeq
Now, since $R_h\to +\infty$, as $h\to +\infty$, then, for all sufficiently large $h$ we have $R_h>R^{\prime\prime}$. In addition, as
$$
\int_{\C_\omega} u_{n_{k_h}}^2 \ dx = \int_{B_{R_h}(y_0)\cap\C_\omega} u_{n_{k_h}}^2 \ dx+ \int_{B_{R_h}^\complement(y_0)\cap\C_\omega} u_{n_{k_h}}^2 \ dx,
$$
then, from \eqref{eq22:teo1sect4} and \eqref{eq24:teo1sect4}, we deduce that for all sufficiently large $h$
$$ \int_{\C_\omega} u_{n_{k_h}}^2 \ dx \geq \lambda-\frac{\eps^\prime}{4} + \frac{\eps^\prime}{2} = \lambda+\frac{\eps^\prime}{4},$$
which contradicts \eqref{eq0:teo1sect4}. The proof of \textit{Claim IV} is complete.
\medskip

We now conclude the proof of \textit{Step 3}. Let $\eps>0$ and let $y_0\in \C_\omega$. Then by \textit{Claim IV}, there exists $R>0$ such that for all $k\in\N$
\beq\label{eq25:teo1sect4}
\int_{B_R^\complement(y_0)\cap\C_\omega} u_{n_k}^2 \ dx< \eps.
\eeq
Let $k\in\N$ and let $v_k \in L^2(\C_\omega)$ be the function defined by $$v_k: = \chi_{B_R(y_0)}u_{n_k},$$ 
 where $\chi_{B_R(0)}$ denotes the characteristic function of $B_R(y_0)$. 
Then we readily check that $(v_k)_k$ is relatively compact in $L^2(\C_\omega)$. Indeed, since $(u_{n_k})_k$  is a bounded sequence in $H^1(\C_\omega)$ and thanks to the compactness of the embedding $H^1(B_R(y_0)\cap \C_\omega) \hookrightarrow L^2(B_R(y_0)\cap \C_\omega)$ we infer that the sequence $(u_{n_k}\big|_{B_R(0)})_k$ is relatively compact in $L^2(B_R(y_0)\cap \C_\omega)$, and thus by definition of $v_k$ we deduce that $(v_k)_k$ is relatively compact in $L^2(\C_\omega)$. 
Finally, from \eqref{eq25:teo1sect4} and by definition of $v_k$, for all $k\in\N$, we have 
\begin{equation}
  \label{eq:sufficient-crit-compactness}
\|u_{n_k}-v_k\|_{L^2(\C_\omega)}< \eps.
\end{equation}
Since $(v_k)_k$ is a relatively compact sequence in $L^2(\C_\omega)$, then \eqref{eq:sufficient-crit-compactness} implies that the set $\{u_{n_k};\ k \in \N \}$ is totally bounded in $L^2(\C_\omega)$, and thus, as $L^2(\C_\omega)$ is a Banach space it follows that $\{u_{n_k};\ k \in \N \}$ is relatively compact. The proof of \textit{Step 3} is complete.
\medskip

\noindent\textbf{Conclusion:} Existence of a minimizer of $\E(\cdot; \C_\omega)$ in $\A_{\omega,c}$.
\medskip

As remarked at the beginning of the proof, we know that $(u_n)_n$ is a bounded sequence in $H^1(\C_\omega)$ and, thanks to \textit{Step 3}, $(u_n)_n$ admits a subsequence which strongly converges in $L^2(\C_\omega)$. Hence, up to a subsequence, there exists $u \in H^1(\C_\omega)$ such that $u_n \rightharpoonup  u$ in $H^1(\C_\omega)$ and $u_n \to u$ in $L^2(\C_\omega)$, as $n\to +\infty$.

Let us set $$\Omega^*:=\{u>0\}.$$ 
We first observe that $\Omega^* \in \A_{\omega,c}$, i.e $\Omega^*$ belongs to the class of admissible sets (see \eqref{eq:defclassadmissiblecylinder}). Indeed, as $u\in H^1(\C_\omega)$ then $\Omega^*$ is a quasi-open subset of $\C_\omega$, and by definition we have that $u\in H_0^1(\Omega^*; \C_\omega)$.  It remains to show that $|\Omega^*|\leq c$. To this end, arguing as in \cite[Proof of Lemma 5.2]{BV1}, since $u_n \to u$ strongly in $L^2(\C_\omega)$, as $n\to +\infty$, and recalling that $|\Omega_n|\leq c+\delta_n$, with $\delta_n\to 0$, as $n\to +\infty$ (see  Lemma \ref{lem1:sect3}-(i)), then, applying Fatou's Lemma, we get that 
$$|\Omega^*|=\int_{\C_\omega} \chi_{\{ u>0\}} \ dx \leq \liminf_{n\to +\infty} \int_{\C_\omega} \chi_{\{ u_n>0\}} \ dx =  \liminf_{n\to +\infty} |\Omega_n| \leq \liminf_{n\to +\infty} (c+\delta_n) = c. $$

We now prove that $\Omega^*$ is a minimizer for $\E(\cdot; \C_\omega)$. 
To this end, we begin observing that, since $u_n \rightharpoonup u$ in $H^1(\C_\omega)$, 
\beq\label{eq26:teo1sect4}
\int_{\C_\omega} |\nabla u|^2 \ dx  \leq \liminf_{n\to+\infty} \int_{\C_\omega} |\nabla u_n|^2 \ dx.
\eeq
In addition, since $u_n \to u$ in $L^2(\C_\omega)$, taking into account that $u_n\in H_0^1(\Omega_n;\C_\omega)$, $u \in H_0^1(\Omega^*;\C_\omega)$, and that $|\Omega_n|\leq c+\delta_n$, $|\Omega^*|\leq c$, then, exploiting H\"older's inequality, we check that, as $n\to +\infty$, it holds
\beq\label{eq27:teo1sect4}
u_n \to u \ \hbox{strongly in $L^1(\C_\omega)$}.
\eeq

Hence, from \eqref{eq26:teo1sect4}, \eqref{eq27:teo1sect4}, and thanks to \eqref{eq-1:teo1sect4}, we obtain that
\beq\label{eqfinal:teo1sect4}
\E(\Omega^*; \C_\omega)\leq J_{\Omega^*}(u)\leq \liminf_{n\to+\infty} \left(\frac{1}{2} \int_{\C_\omega} |\nabla u_n|^2 \ dx - \int_{\C_\omega} u_n \ dx \right)=\liminf_{n\to+\infty} \E(\Omega_n; \C_\omega) \leq \mathcal{O}_c(\C_\omega),
\eeq
where $J_{\Omega^*}:H_0^1(\Omega^*; \C_\omega)\to \R$ is the functional defined in \eqref{eq:deftorsion}. Finally, as $\Omega^*\in \A_{\omega,c}$, from \eqref{eqfinal:teo1sect4} we conclude that $\E(\Omega^*; \C_\omega)=\mathcal{O}_c(\C_\omega)$. Hence $\Omega^*$ is a minimizer for $\E(\cdot; \C_\omega)$ and the proof is complete.
\end{proof}

\section{Qualitative and regularity properties of minimizers}\label{S:properties}

In this section we show some qualitative and topological properties of minimizers for $\E(\cdot; \C_\omega)$ in $\A_{\omega,c}$, with $c>0$ fixed. More precisely, we will prove that if $\Omega^*$ is a minimizer, then $\Omega^*$ is bounded (Proposition \ref{prop1:sect5}), it saturates the constraint, namely $|\Omega^*|=c$ (Corollary \ref{cor1:sect5}), and is open (Proposition \ref{prop2:sect5}). Moreover, we discuss the regularity of the relative boundary of $\Omega^*$, providing the proof of Theorem \ref{teoreg}. This is important in order to obtain that $\Omega^*$ is connected (Theorem \ref{connected}).

We begin with a preliminary technical lemma. 

\begin{lemma}\label{lem1:sect5}
Let $\Omega^*$ be a minimizer of $\E(\cdot; \C_\omega)$ in $\A_{\omega,c}$. Then for any $\widetilde \Omega \in \A_{\omega, c}$ such that $0<|\widetilde \Omega|\leq |\Omega^*|$ it holds
 \beq\label{eq:tesi:lem2}
\frac{\E(\widetilde\Omega; \C_\omega)}{|\widetilde\Omega|}\geq  \frac{\E(\Omega^*; \C_\omega)}{|\Omega^*|}.
\eeq
\end{lemma}
\begin{proof}
Let $\Omega^*$ be a minimizer for $\E(\cdot; \C_\omega)$ in $\A_{\omega,c}$. We observe that since $\O_c(\C_\omega)<0$ (see Remark \ref{remark:negativityenergy}) we have
\beq\label{eq0:prop1:sect5}
|\Omega^*|>0,  \ \ \E(\Omega^*; \C_\omega)<0.
\eeq
Indeed, if $|\Omega^*|=0$ then by Proposition \ref{prop1:sect2} we would get that $\O_c(\C_\omega)=\E(\Omega^*; \C_\omega)=0$ contradicting $\O_c(\C_\omega)<0$.

For any $t>0$ let us consider the diffeomorphism $F_t:\overline\C_\omega\to \overline\C_\omega$ defined by
 \beq\label{eq:defFt}
F_t(x^\prime, x_N):=(x^\prime, tx_N).
\eeq
It is immediate to check that for all $(x^\prime,x_N)\in \C_\omega$ it holds $$\mathrm{Jac(F_t)}(x^\prime,x_N)=\left[\begin{array}{cc}\mathbb{I}_{N-1}&(0^\prime)^T\\  0^\prime&t \end{array}\right],$$
where $\mathbb{I}_{N-1}$ is the identity matrix of order $N-1$, $0^\prime$ is the null vector in $\R^{N-1}$, $T$ is the transposition. In particular, $\textrm{Det}(\mathrm{Jac(F_t)}(x^\prime,x_N))=t$, for all $(x^\prime, x_N)\in\C_\omega$, and by the formula of change of variable we infer that for any $\mathcal{L}^N$-measurable set $E\subset\C_\omega$ we have
\beq\label{eq1:lem2sect5}
|F_t(E)| = t|E|.
\eeq
Now, let $\widetilde\Omega \in \A_{\omega, c}$ be such that $0<|\widetilde \Omega|\leq |\Omega^*|$ and let  $u_{\widetilde\Omega} \in H_0^1(\widetilde\Omega; \C_\omega)$ be the energy function of $\widetilde \Omega$, and choose
\beq\label{eq2:lem2sect5}
t:=\frac{|\Omega^* |}{|\widetilde\Omega|}.
\eeq
By construction we have $t\geq 1$. Moreover, by definition and by changing variable in the integral, we infer that
\beq\label{eq3:prop1:sect5}
\begin{array}{lll}
\displaystyle t \E(\widetilde\Omega; \C_\omega)&=&\displaystyle t \int_{\widetilde\Omega} \left(\frac{1}{2} |\nabla u_{\widetilde\Omega}|^2 - u_{\widetilde\Omega} \right)\, dx\\[9pt]
&=&\displaystyle  \int_{F_{t}(\widetilde\Omega)} \left(\frac{1}{2} \left|\nabla u_{\widetilde\Omega}\left(x^\prime, \frac{x_N}{t}\right)\right|^2 - u_{\widetilde\Omega}\left(x^\prime, \frac{x_N}{t}\right) \right)\, dx.
\end{array}
\eeq
Now, setting $\tilde u:= u_{\widetilde\Omega}\circ F_{t}^{-1}$, namely
\beq\label{eq4:prop1:sect5}
\tilde u(x^\prime, x_N)=u_{\widetilde\Omega}\left(x^\prime, \frac{x_N}{t}\right), \ \ \ \ (x^\prime, x_N)\in \C_\omega,
\eeq
we readily check that $\tilde u\in H_0^1(F_{t}(\widetilde\Omega); \C_\omega)$, and, denoting by $\nabla^\prime$ the gradient with respect to the variables $x_1,\ldots, x_{N-1}$, it holds
\beq\label{eq5:prop1:sect5}
\nabla^\prime \tilde u(x^\prime, x_N)=\nabla^\prime u_{\widetilde\Omega}\left(x^\prime, \frac{x_N}{t}\right), \ \ \ \frac{\partial  \tilde u }{\partial x_N} (x^\prime, x_N)=  \frac{1}{t} \frac{\partial  u_{\widetilde\Omega}}{\partial x_N} \left(x^\prime, \frac{x_N}{t}\right),
\eeq
for a.e. $(x^\prime, x_N)\in \C_\omega$. Then, from \eqref{eq3:prop1:sect5}--\eqref{eq5:prop1:sect5}, exploiting that $t\geq 1$, and since $\tilde u\in H_0^1(F_{t}(\widetilde\Omega); \C_\omega)$, we obtain
\begin{equation}\label{eq6:prop1:sect5}
\begin{split}
t \E(\widetilde\Omega; \C_\omega)
&=\int_{F_{t}(\widetilde\Omega)}\left(\frac{1}{2} \left|\nabla^\prime \tilde u \right|^2 + t^2 \left(\frac{\partial  \tilde u }{\partial x_N}\right)^2 -  \tilde u \right)\, dx\\
&\geq\int_{F_{t}(\widetilde\Omega)}\left(\frac{1}{2} \left|\nabla \tilde u \right|^2  -  \tilde u \right)\, dx
=J_{F_{t}(\widetilde\Omega)} (\tilde u)
\geq\E(F_{t}(\widetilde\Omega); \C_\omega).
\end{split}
\end{equation}
Finally, since by construction $|F_{t}(\widetilde\Omega)|=|\Omega^*|$ (see \eqref{eq1:lem2sect5}--\eqref{eq2:lem2sect5}), and as $\Omega^*$ is the minimizer of $\E(\cdot;\C_\omega)$, we have
\beq\label{eq7:prop1:sect5}
 \E(F_{t}(\widetilde\Omega); \C_\omega)\geq \E(\Omega^*; \C_\omega),
\eeq
and thus, from \eqref{eq6:prop1:sect5}--\eqref{eq7:prop1:sect5}, we obtain
%\beq\label{eq8:prop1:sect5}
$$
 t \E(\widetilde\Omega; \C_\omega) \geq \E(\Omega^*; \C_\omega),
$$
%\eeq
which readily gives \eqref{eq:tesi:lem2}. The proof is complete.
\end{proof}

We now show that any minimizer of $\E(\cdot;\C_\omega)$ is a bounded subset of $\C_\omega$.
We begin by recalling the definition of local shape subsolution for the energy introduced by D. Bucur in \cite[Sect. 2]{BU2}, adapted here for the relative torsion energy (see also \cite[Sect. 2.1.2]{BHS}).

\begin{definition}
We say that a quasi-open set $\Omega\subset\C_\omega$ of finite measure  is a \textit{local shape subsolution} for the relative torsion energy $\E(\cdot;\C_\omega)$ if there exist $\delta>0$ and $\Lambda>0$ such that for any quasi-open subset $\widetilde\Omega\subset \Omega$ with
$\|u_{\Omega}-u_{\widetilde\Omega}\|_{L^2(\C_\omega)}< \delta$ it holds
$$ \E(\Omega; \C_\omega)+\Lambda |\Omega| \leq  \E(\widetilde\Omega; \C_\omega)+\Lambda |\widetilde\Omega|.$$
\end{definition}

This definition enters to prove the following result.

\begin{proposition}\label{prop1:sect5}
Let $c>0$ and let $\omega\subset\R^{N-1}$ be a Lipschitz bounded domain.
If $\Omega^*$ is a minimizer for $\E(\cdot;\C_\omega)$ in $\A_{\omega,c}$ then $\Omega^*$ is bounded.
\end{proposition}
\begin{proof}
To prove the boundedness of $\Omega^*$, by \cite[Theorem 1]{BU2} it is sufficient to show that $\Omega^*$ is a local shape subsolution for the relative torsion energy $\E(\cdot;\C_\omega)$.  To this end assume by contradiction that there exist a sequence $(\Lambda_n)_n$ of strictly positive numbers, with $\Lambda_n\to 0^+$, as $n\to +\infty$, and an increasing sequence $(\widetilde \Omega_n)_n$ of quasi-open subsets of $\Omega^*$ such that $\|u_{\Omega^*}-u_{\widetilde \Omega_n}\|_{L^2(\C_\omega)}\to 0$, as $n\to +\infty$, and for all $n\in \N$ it holds
\beq\label{eq1:prop1:sect5}
\E(\Omega^*; \C_\omega)+\Lambda_n |\Omega^*| > \E(\widetilde\Omega_n; \C_\omega)+\Lambda_n |\widetilde\Omega_n|.
\eeq
Without loss of generality we can assume that 
\beq\label{eq1bis:prop1:sect5}
|\widetilde \Omega_n|>0
\eeq
for all $n\in \N$. Moreover we notice that since $\widetilde \Omega_n\subset \Omega^*$ for all $n\in\N$ we have $ \E(\Omega^*; \C_\omega)\leq \E(\widetilde\Omega_n; \C_\omega)$, and thus from \eqref{eq1:prop1:sect5} we readily obtain that
\beq\label{eq2:prop1:sect5}
|\widetilde\Omega_n|<|\Omega^*|.
\eeq
In view of \eqref{eq1bis:prop1:sect5}, \eqref{eq2:prop1:sect5} we can apply Lemma \ref{lem1:sect5} to $\widetilde\Omega_n$, for all $n\in \N$, obtaining that
\beq\label{eq9:prop1:sect5}
\E(\widetilde\Omega_n; \C_\omega) \geq \frac{|\widetilde\Omega_n|}{|\Omega^*|}\,\E(\Omega^*; \C_\omega).
\eeq
Then, from \eqref{eq1:prop1:sect5} %{eq8:prop1:sect5} 
and exploiting \eqref{eq9:prop1:sect5}, we infer that
$$
\Lambda_n \left( |\Omega^*| -  |\widetilde\Omega_n|\right)> \E(\widetilde\Omega_n; \C_\omega) - \E(\Omega^*; \C_\omega)\geq\left(\frac{|\widetilde\Omega_n|}{|\Omega^*|}-1\right)\E(\Omega^*; \C_\omega)= -\frac{|\Omega^*|-|\widetilde\Omega_n|}{|\Omega^*|}\,\E(\Omega^*; \C_\omega),
$$
from which we deduce that
\beq\label{eq10:prop1:sect5}
 \Lambda_n |\Omega^*| > -\E(\Omega^*; \C_\omega).
\eeq
Then, taking the limit as $n\to +\infty$ in \eqref{eq10:prop1:sect5}, recalling that $\Lambda_n\to 0^+$, and since $-\E(\Omega^*; \C_\omega)>0$ (see \eqref{eq0:prop1:sect5}) we obtain a contradiction. The proof is complete.
\end{proof}

\begin{corollary}\label{cor1:sect5}
Let $\Omega^*$ be a minimizer of $\E(\cdot;\C_\omega)$ in $\A_{\omega,c}$. Then $|\Omega^*|=c$. 
\end{corollary}
\begin{proof}
Let $\Omega^*\in\A_{\omega,c}$ be a minimizer of $\E(\cdot;\C_\omega)$. Thanks to Proposition \ref{prop1:sect5} we know that $\Omega^*$ is a bounded subset of $\C_\omega$. 
Assume by contradiction that $|\Omega^*|<c$. Then, as $\Omega^*$ is a bounded subset of $\C_\omega$, we find $r_0>0$ and $x_0\in \C_\omega\setminus \overline\Omega^*$ such that $\overline{B_r(x_0)} \subset \C_\omega\setminus \overline\Omega^*$ for all $0<r<r_0$. Therefore, fixing $0<r<r_0$ small enough so that $|B_r(x_0)|\leq c-|\Omega^*|$ and setting $$\Omega^*_r:=\Omega^*\cup B_r(x_0),$$
we have that $\Omega^*_r\in \A_{\omega,c}$. In addition, by construction and exploiting Lemma \ref{additivity} we have $$\E(\Omega^*_r; \C_\omega)=\E(\Omega^*; \C_\omega) + \E(B_r(x_0); \C_\omega) <\E(\Omega^*; \C_\omega),$$ because $\overline{\Omega^*}\cap \overline{B_r(x_0)}=\varnothing$ and $\E(B_r(x_0); \C_\omega) <0$. But this contradicts the minimality of $\Omega^*$ and the proof is complete.
\end{proof} 

\begin{remark}\label{rem1:sect5}
A straightforward consequence of Corollary \ref{cor1:sect5} is that the map $\Phi\colon \left]0,+\infty\right[\to \left]-\infty, 0\right[$ defined by $\Phi(c):=\O_c(\C_\omega)$ is strictly monotone decreasing. 
\end{remark}
In view of the previous remark the map $\Phi$ is continuous, up to a negligible set under $\mathcal{L}^1$. Actually we can prove that $\Phi$ is everywhere continuous. This is the content of the following.
\begin{proposition}\label{prop3:sect5}
Let $\omega\subset\R^{N-1}$ be a Lipschitz bounded domain.
The map $\Phi\colon\left]0,+\infty\right[\to \left]-\infty, 0\right[$ defined by $\Phi(c):=\O_c(\C_\omega)$ is continuous.
\end{proposition}
\begin{proof}
Let $c_0 \in ]0,+\infty[$. Since $c\mapsto \O_c(\C_\omega)$ is monotone decreasing and defined in the whole $]0,+\infty[$ then the limits $l_{c_0}^-:=\lim_{c\to c_0^-} \O_c(\C_\omega)$, $l_{c_0}^+:=\lim_{c\to c_0^+} \O_c(\C_\omega)$ exist and are finite. We first show that 
\beq\label{eq1:prop3:sect5}
l_{c_0}^-=\O_{c_0}(\C_\omega).
\eeq
Indeed, thanks to Theorem \ref{teo:mainteoexistmin} there exists $\Omega^* \in \A_{\omega, c_0}$ such that $\E(\Omega^*; \C_\omega)=\O_{c_0}(\C_\omega)$, and by Proposition \ref{prop1:sect5} and Proposition \ref{prop2:sect5} we know that $\Omega^*$ is a bounded open subset of $\C_\omega$. Let us fix $x_0\in \Omega^*$ and let $r_{\Omega^*}>0$ be such that $\overline{B_r(x_0)}\subset \Omega^*$ for all $0<r<r_{\Omega^*}$. For $0<r<r_{\Omega^*}$ we set
$$
 \Omega^*_{r}:=\Omega^* \cap \{x\in \RN;\, |x-x_0|\geq r\}.
$$
We observe that $\Omega^*_{r}$ is a quasi-open subset of $\C_\omega$ with $|\Omega^*_r|\leq c_0-\omega_Nr^N$. Indeed, for all $0<r<r_{\Omega^*}$, it holds $\Omega^*_r\cup B_r(x_0)=\Omega^*$, $\Omega^*$ and $B_r(x_0)$ are open sets, and $\mathrm{cap}(B_r(x_0))\to 0$ as $r\to 0^+$ (as it follows from \cite[(3.55)-(3.57)]{HP}). Then, arguing as in the proof of Lemma \ref{lem:tech2} we get that $\Omega^*_r$ $\gamma$-converges to $\Omega^*$ as $r\to 0^+$. This implies that $\E(\Omega^*_r; \C_\omega)\to \E(\Omega^*;\C_\omega)=\O_{c_0}(\C_\omega)$ as $r\to 0^+$. In addition, since $\Omega^*_r\in \A_{\omega, c_r}$, where $c_r:=c_0-\omega_Nr^N<c_0$, for all $0<r<r_{\Omega^*}$, then $\O_{c_r}(\C_\omega)\leq\E(\Omega^*_r; \C_\omega)$ and taking to the limit, as $r\to 0^+$, we infer that
$ l_{c_0}^- \leq \O_{c_0}(\C_\omega)$. On the other hand, by monotonicity of $c\mapsto \O_c(\C_\omega)$, it holds that $ l_{c_0}^- \geq \O_{c_0}(\C_\omega)$, and thus the only possibility is $l_{c_0}^-= \O_{c_0}(\C_\omega)$, as desired.

To conclude the proof we show that $l_{c_0}^-=l_{c_0}^+$. To this end we first observe that by monotonicity we have 
\beq\label{eq2:prop3:sect5}
l_{c_0}^-\geq l_{c_0}^+.
\eeq
To prove the reverse inequality, let us consider a sequence $(c_n)_n$ in $\R^+$ such that $c_n>c_0$ for all $n\in\N$ and $c_n\to c_0$, as $n\to +\infty$. In view of Theorem \ref{teo:mainteoexistmin} and Proposition \ref{prop1:sect5}, for all $n\in\N$, there exists an open set $\Omega_n$ such that $|\Omega_n|\leq c_n$ and $\E(\Omega_n; \C_\omega)=\O_{c_n}(\C_\omega)$. Moreover, up to considering the Steiner symmetrization of $\Omega_n$ and of its energy function $u_{\Omega_n}$, we can assume without loss of generality that $\Omega_n$ satisfies the properties (ii)-(iii) and $u_{\Omega_n}$ satisfies the properties (v)-(vi) of Lemma \ref{lem1:sect3}. Then, arguing as in the proof of Theorem \ref{teo:mainteoexistmin}, taking into account that $|\Omega_n|\leq c_n$, we infer that there exists $\Omega \in \A_{\omega,c}$ such that $u_{\Omega_n} \to u_{\Omega^*}$  strongly in $L^2(\C_\omega)$, as $n\to +\infty$, where $u_{\Omega^*}$ is the energy function associated $\Omega^*$, and $\E(\Omega_n;\C_\omega)\to \E(\Omega^*;\C_\omega)$. Hence, as $\Omega^*\in \A_{\omega, c_0}$, we have $\O_{c_0}(\C_\omega) \leq \E(\Omega^*;\C_\omega)$ and as $\E(\Omega_n; \C_\omega)=\O_{c_n}(\C_\omega)$, then, taking the limit, as $n\to +\infty$, we conclude that
\beq\label{eq3:prop3:sect5}
\O_{c_0}(\C_\omega) \leq l_{c_0}^+.
\eeq
Finally, combining \eqref{eq1:prop3:sect5}, \eqref{eq2:prop3:sect5} and \eqref{eq3:prop3:sect5} we conclude that $l_{c_0}^-= \O_{c_0}(\C_\omega)=l_{c_0}^+$, i.e. $\Phi$ is continuous at $c_0$, and, as $c_0$ is arbitrary, we are done. The proof is complete.
\end{proof}

\begin{proposition}\label{prop2:sect5}
If $\Omega^*$ is a minimizer for $\E(\cdot;\C_\omega)$ in $\A_{\omega,c}$, with $c>0$ fixed, then the energy function $u_{\Omega^*}\in H_0^1(\Omega^*; \C_\omega)$ is Lipschitz continuous in any Lipschitz domain $D \subset \C_\omega$ such that $\overline{D} \subset\C_\omega$ and $\Omega^*=\{u_{\Omega^*}>0\}$ is an open subset of $\C_\omega$.
\end{proposition}
\begin{proof}
The proof is similar to that of \cite[Theorem 2.14]{MPV} (see also \cite{BRP} for the Dirichlet case).
\end{proof}

We now give the proof of Theorem \ref{teoreg}.
\begin{proof}[Proof of Theorem \ref{teoreg}]
Let $\Omega^*$ be a minimizer for $\E(\cdot;\C_\omega)$ in $\A_{\omega,c}$ and let $u_{\Omega^*}\in H^1_0(\Omega; \C_\omega)$ be its energy function. To ease the notation we set $u_*:=u_{\Omega^*}$. Thanks Proposition \ref{prop1:sect2}, Proposition \ref{prop1:sect5}, Proposition \ref{prop2:sect5}, and Corollary \ref{cor1:sect5}, we have that $u_{*}$ is bounded, locally Lipschitz in $\C_\omega$, $\Omega^*=\{u_*>0\}$ is a bounded open set and $|\Omega^*|=c$. Since $u_*$ is a non-negative function, we have $\Omega_{u_{*}}=\Omega^*$, where $\Omega_u$ denotes, for a generic $u \in H^1(\C_\omega)$, the set $\Omega_u:=\{ u\neq 0\}$. In order to prove the result we divide the proof in some steps.\vspace{6pt}

\noindent\textbf{Step 1}:\ $u_{*}$ is a solution to the following  problem
\begin{equation}\label{eq:altProb}
\begin{cases}
u\in H^1(\C_\omega),\ |\Omega_{u}|= c,\\
J(u)\leq J(v)\ \forall v\in\mathbb{V}_c,
\end{cases}
\end{equation}
where $\mathbb{V}_c:=\{w\in H^1(\C_\omega);\, |\Omega_w|\leq c\}$ and $J:\mathbb{V}_c\to \R$ is the functional defined by
$$J(w):=\frac{1}{2}\int_{\C_\omega} |\nabla w|^2\, dx -\int_{\C_\omega} w\, dx.$$
\vspace{3pt}

Indeed, by definition $u_*\in H^1(\C_\omega)$ and saturates the constraint, namely $|\Omega_{u_{*}}|=c$. To conclude it suffices to show that $J(u_*)\leq J(v)$ for all $v\in\mathbb{V}_c$.%Clearly $u_{*}\in \mathbb{V}_c$ and thus $\inf\{J(u);\ u \in \mathbb{V}_c \}\leq \O_c(\C_\omega)$. We claim that $\inf\{J(u);\ u \in \mathbb{V}_c \}=\O_c(\C_\omega)$.

To this end, assume by contradiction that there exists $v\in \mathbb{V}_c$ such that 
\beq\label{eq1:sketchdimteoreg}
J(v)<J(u_{*})=\O_c(\C_\omega),
\eeq
and let us consider the function $v^+$, i.e. the positive part of $v$. Since $v^+\in H^1(\C_\omega)$ and $\Omega_{v^+}=\{v>0\}$ is a quasi-open subset of $\C_\omega$, with $|\Omega_{v^+}|\leq c$, then $\Omega_{v^+} \in \A_{\omega,c}$. Moreover, since $v\leq v^+$ and $|\nabla v^+|\leq |\nabla v|$ a.e. in $\C_\omega$, it holds that 
\beq\label{eq2:sketchdimteoreg}
J(v^+) \leq J(v). 
\eeq
In addition, by definition we check that $v^+ \in H_0^1(\Omega_{v^+}; \C_\omega)$ and thus we have
$$\E(\Omega_{v^+}; \C_\omega) \leq J(v^+),$$
which, together with \eqref{eq1:sketchdimteoreg} and \eqref{eq2:sketchdimteoreg}, implies that $\E(\Omega_{v^+}; \C_\omega)<\O_c(\C_\omega)$, a contradiction. The proof of \textit{Step 1} is complete.\vspace{6pt}

The next step states the existence of a non-negative Lagrange multiplier.\vspace{6pt}

\noindent\textbf{Step 2}: There exists a constant $\Lambda=\Lambda_{u_*}\geq 0$ such that for all $\Phi \in C^\infty_c(\C_\omega, \R^N)$ it holds
\beq\label{eq:EulerLagrangeeqproofteo13}
\int_{\C_\omega} \nabla u_*\boldsymbol{\cdot} D\Phi\nabla u_*\, dx - \frac{1}{2}\int_{\C_\omega} |\nabla u_*|^2 \mathrm{div} \Phi\, dx= \int_{\C_\omega} \nabla u_*\boldsymbol{\cdot} \Phi\, dx + \Lambda \int_{\Omega_{u_*}} \mathrm{div}\Phi\, dx. 
\eeq
Indeed, thanks to \textit{Step 1}, $u_*$ is a solution to Problem \eqref{eq:altProb}, and thus,  taking $f\equiv 1$ in \cite[Proposition 1.2]{B} we obtain the desired relation.\vspace{6pt}

In the next two steps we show that any solution to \eqref{eq:altProb} is also a solution of a ``pseudo-penalized" problem (see \cite[Sect. 4]{B}, \cite[Sect. 2]{BL} or \cite[Sect. 4.4.]{LS2}). We begin with a technical property:\vspace{6pt}

\noindent\textbf{Step 3}: Let $x_0\in \Gamma_{\Omega^*}= \partial\Omega_{u_*}\cap \C_\omega$ and let $r_0>0$ be such that $B_{r_0}(x_0)\subset\subset \C_\omega$. Then 
$$0<|\Omega_{u_*}\cap B_{r_0}(x_0)|<|B_{r_0}(x_0)|.$$

Clearly $|\Omega_{u_*}\cap B_{r_0}(x_0)|>0$ because $\Omega_{u_*}\cap B_{r_0}(x_0)$ is the intersection of two open sets with nonempty intersection (by construction). The second inequality $|\Omega_{u_*}\cap B_{r_0}(x_0)|<|B_{r_0}(x_0)|$ is a consequence of a standard argument (see the proof of \cite[Lemma 2.5]{BL}) that can be easily adapted to our setting. Namely, if by contradiction $|\Omega_{u_*}\cap B_{r_0}(x_0)|=|B_{r_0}(x_0)|$, then we can show that $-\Delta u_*=1$ in $B_{r_0}(x_0)$ and thus, in particular, it follows that $B_{r_0}(x_0)\subset \Omega_{u_*}$, which is not possible because $x_0$ is a boundary point for $\Omega_{u_*}$.

To prove this fact, assume that $|\Omega_{u_*}\cap B_{r_0}(x_0)|=|B_{r_0}(x_0)|$. This means that $u_*>0$ a.e. in $B_{r_0}(x_0)$. Let us consider $v\in H^1(\C_\omega)$ such that $v=u_*$ outside $B_{r_0}(x_0)$ (in particular $v-u_*\in H_0^1(B_{r_0}(x_0)$) and $-\Delta v= 1$ in $B_{r_0}(x_0)$. By the strong maximum principle we get that $v>0$ in $B_{r_0}(x_0)$ and thus we have $|\Omega_v|=|\Omega_{u_*}|$. Now, by minimality of $u_*$ and exploiting that $v=u_*$ outside $B_{r_0}(x_0)$, from $J(u_*)\leq J(v)$ we infer that
\beq\label{eq3:sketchdimteoreg}
\frac{1}{2}\int_{ B_{r_0}(x_0)} (|\nabla u_*|^2-|\nabla v|^2)\, dx-  \int_{ B_{r_0}(x_0)} (u_*-v) \, dx\leq 0.
\eeq
On the other hand, as $-\Delta v=1$ in $B_{r_0}(x_0)$, taking $v-u_*\in H_0^1(B_{r_0}(x_0))$ as test function, we get that
\beq\label{eq4:sketchdimteoreg}
\int_{B_{r_0}(x_0)} \nabla v \boldsymbol{\cdot}  \nabla u_* \, dx= \int_{B_{r_0}(x_0)} | \nabla v|^2\, dx + \int_{B_{r_0}(x_0)} (u_*-v)\, dx.
\eeq
From \eqref{eq3:sketchdimteoreg} and \eqref{eq4:sketchdimteoreg}, we obtain
$$\int_{ B_{r_0}(x_0)} |\nabla (u_*- v)|^2\, dx=\int_{ B_{r_0}(x_0)} (|\nabla u_*|^2-|\nabla v|^2)\, dx- 2 \int_{ B_{r_0}(x_0)} (u_*-v) \, dx\leq 0,$$
which implies that $u_*=v$ a.e. in $B_{r_0}(x_0)$, and by continuity, $u_*=v>0$ in $B_{r_0}(x_0)$. The proof of \textit{Step 3} is complete.\vspace{6pt}

Thanks to \textit{Step 3}, it follows that the conditions stated in \cite[(22)]{B} are satisfied with $D=\C_\omega$, $D_1=B_{r_0}(x_0)$. Indeed, as $0<|\Omega_{u_*}\cap B_{r_0}(x_0)|<|B_{r_0}(x_0)|$, then, there exists $h_0=h_0(x_0,r_0,u_*)>0$ such that $0<h_0<|\Omega_{u_*}\cap B_{r_0}(x_0)|<|B_{r_0}(x_0)|-h_0$ and $0<h_0< |(\C_\omega\setminus B_{r_0}(x_0))\cap \Omega_{u_*}|$. Moreover, up to taking a smaller radius we can assume without loss of generality that $u_*\not\equiv0$ on $\partial B_{r_0}(x_0)$ (otherwise if $u_*\equiv0$ on $\partial B_{r}(x_0)$, for all $0<r\leq r_0$ then $u_*\equiv0$ on $B_{r_0}(x_0)$, contradicting  $|\Omega_{u_*}\cap B_{r_0}(x_0)|>0$). Finally, the last condition in \cite[(22)]{B} is trivially satisfied since $|(\C_\omega\setminus B_{r_0}(x_0))\cap \Omega_{u_*}|\leq c<+\infty=|\C_\omega\setminus B_{r_0}(x_0)|$.

Then, similarly to \cite[Sect. 4]{B} (see also \cite[Sect. 4.4]{LS2}, \cite[Sect. 2]{BL}), we define the class
\beq\label{defsetF}
\mathcal{F}=\mathcal{F}(u_*,x_0,r_0):=\{v\in H^1(\C_\omega);\, v-u_*\in H_0^1(B_{r_0}(x_0))\},
\eeq
and for $h>0$ we set
\begin{eqnarray*}
\mu_{-}(h)&:=&\sup\{\mu\geq0;\, J(u_*)+\mu |\Omega_{u_*}| \leq J(v)+\mu |\Omega_{v}|\ \hbox{$\forall v\in \mathcal{F}$ \, such that }\, c-h \leq |\Omega_{v}| \leq c\},\\
\mu_{+}(h)&:=&\inf\{\mu\geq0;\, J(u_*)+\mu |\Omega_{u_*}| \leq J(v)+\mu |\Omega_{v}|\ \hbox{$\forall v\in \mathcal{F}$ \, such that }\, c\leq |\Omega_{v}| \leq c+h\}.
\end{eqnarray*}
For simplicity we omit the dependence on $x_0$, $r_0$, $u_*$ in the notation for $\mu_{-}$, $\mu_{+}$. The following asymptotic property for $\mu_{-}$, $\mu_{+}$ holds:\vspace{6pt}

\textbf{Step 4}:  Let $\Lambda=\Lambda_{u_*}\geq 0$ be the Lagrange multiplier given by \textit{Step 2}. Let $x_0\in \Gamma_{\Omega^*}= \partial\Omega_{u_*}\cap \C_\omega$ and let $r_0>0$ be such that $B_{r_0}(x_0)\subset\subset \C_\omega$. If $\Lambda>0$, then there exists $h_0>0$ such that
$$0<\mu_{-}(h) \leq \Lambda \leq \mu_+(h) < +\infty \ \ \forall h\in]0,h_0[. $$
Moreover, it holds that
$$\lim_{h\to 0^+} \mu_+(h)=\lim_{h\to 0^+} \mu_-(h)=\Lambda.$$

The proof is the same as that of \cite[Proposition 4.3]{B} (see also \cite[Theorem 1.5]{BL} or \cite[Proposition 4.12]{LS2}).\vspace{6pt}

\noindent\textbf{Step 5}: Let $\Lambda=\Lambda_{u_*}\geq 0$ be the number given by \textit{Step 2}. It holds that $\Lambda>0$ and $\Omega_{u_*}$ has locally finite perimeter.\vspace{6pt}

The proof of the positivity of the Lagrange multiplier $\Lambda$ follows from \cite[Proposition 6.1]{B} (see also the proof of \cite[Proposition 4.17]{LS2}), while by \cite[Theorem 2.4]{B} it follows that $\Omega_{u_*}$ has locally finite perimeter.\vspace{6pt}

The next two steps are about, respectively, the non-degeneracy property of the minimizer $u_*$ and a density estimate. Their proof is an adaptation of the arguments contained in \cite{AC}, \cite{GS} and can be found in \cite[Proposition 5.3, Proposition 5.5]{B}.\vspace{6pt}

\noindent\textbf{Step 6}: For all $\tau\in]0,1[$ there exist two positive numbers $C_1$, $\bar r_1$ such that for any ball $B_r(x_1)\subset B_{r_0}(x_0)$, with $0<r\leq \bar r_1$, if
$$\frac{1}{r}\fint_{\partial B_r(x_1)} u_*\, d\sigma \leq C_1,$$
then $u=0$ in $B_{\tau r}(x_1)$.\vspace{6pt}

\noindent\textbf{Step 7}: There exist positive numbers $\delta$, $\bar r_2$, with $\delta<1$, such that for any ball $B_r(x_1)\subset B_{r_0}(x_0)$, with $0<r\leq \bar r_2$, we have
$$ 
0<\delta\leq \frac{|B_r(x_1)\cap \Omega_{u_*}|}{|B_r(x_1)|}\leq 1-\delta<1.
$$

We now consider the Weiss boundary adjusted energy, introduced in \cite{W1}, and we prove a monotonicity formula which is the counterpart, for our problem, to that one obtained in \cite[Proposition 4.22]{LS2}. Following the notations of \cite[Chap. 9]{V}, for any $\Lambda\ge 0$ and $u\in H^1(B_1(0))$ we define 
$$W_{\Lambda}(u):=\int_{B_1(0)} |\nabla u|^2\, dx -\int_{\partial B_1(0)} u^2\, d\sigma +\Lambda |\{u>0\}\cap B_1(0)|.$$
If $u\in H^1(B_r(x_0))$ we define the rescaled function $u_{x_0,r}:B_1(0)\to \R$ as
\beq\label{eq:defblowup}
u_{x_0, r}(x):=\frac{1}{r}u(x_0+ rx),
\eeq
and the one-homogenous extension of $u_{x_0,r}$ in $B_1(0)$ as the function $z_{x_0,r}\colon B_1(0)\to \R$ defined by
$$
z_{x_0,r}(x):=|x|u_{x_0, r}\left(\frac{x}{|x|}\right).
$$
It is well known (see \cite[Chap 9, Sect. 1]{V}) that
$$
W_\Lambda(u_{x_0, r})=\frac{1}{r^N}\int_{B_r(x_0)} |\nabla u|^2\, dx -\frac{1}{r^{N+1}}\int_{\partial B_r(x_0)} u^2\, d\sigma +\frac{\Lambda}{r^N} |\{u>0\}\cap B_r(x_0)|.
$$

\noindent\textbf{Step 8}: Let $\Lambda=\Lambda_{u^*}$ be the number given by \textit{Step 2}. Let $x_0\in \Gamma_{\Omega^*}=\partial\Omega_{u_*}\cap \C_\omega$ and let $r_0>0$ be such that $B_{r_0}(x_0)\subset\subset \C_\omega$. There exists a positive number $C$ such that for all $0<r<r_0$ it holds
\beq\label{eq:Step5formula}
\frac{\partial W_\Lambda((u_*)_{x_0, r})}{\partial r}\geq \frac{2}{r^N} \int_{\partial B_r(x_0)}\left(\frac{u_*}{r} - (x\cdot\nabla {u_*})  \right)^2\, d\sigma - C.
\eeq
Moreover
$$\lim_{r\to 0^+} W_\Lambda((u_*)_{x_0, r})$$
exists and it is finite.\vspace{6pt}

Without loss of generality we can assume that $x_0=0$ and to ease the notation we write, respectively, $u$, $u_r$, $z_r$, instead of $u_*$, $(u_*)_{r,0}$, $z_{0,r}$. Then, by definition and by the computations contained in the proof of \cite[Lemma 9.2]{V} we get that for all $r\in]0,r_0[$ it holds
\begin{equation}\label{eq1:Step5Pteo13}
\begin{split}
\displaystyle W_\Lambda(z_r)-W_\Lambda(u_r)&=\displaystyle \frac{1}{N}\int_{\partial B_1(0)} \left(|\nabla u_r|^2+\Lambda \chi_{\{u_r>0\}} -(x\cdot \nabla u_r)^2 \right)\,d\sigma\\
&\quad - \frac{N-1}{N}\int_{\partial B_1(0)} u_r^2\,d\sigma\displaystyle -\int_{B_1(0)} (|\nabla u_r|^2+\Lambda \chi_{\{u_r>0\}})\, dx + \int_{\partial B_1(0)}  u_r^2\,d\sigma
 \end{split}
\end{equation}
because $|\{z_r>0\}\cap B_1|=\frac1N\left|\{u_r>0\}\cap \partial B_1\right|$. On the other hand, exploiting the Euler-Lagrange equation \eqref{eq:EulerLagrangeeqproofteo13} and arguing as in the proof of \cite[Lemma 9.8]{V}, we infer that
\begin{equation}\label{eq2:Step5Pteo13}
\begin{split} 
N\int_{B_r(0)} (|\nabla u|^2 &- 2 u + \Lambda \chi_{\{u>0\}})\, dx\\
&\quad=\displaystyle 2 \int_{B_r(0)} |\nabla u|^2 \, dx + r \int_{\partial B_r(0)} \left(|\nabla_\tau u|^2- \left|\frac{\partial u}{\partial \nu}\right|^2+\Lambda \chi_{\{u>0\}} - 2u\right)\,d\sigma\\
&\quad=\displaystyle 2 \int_{B_r(0)} u \, dx + r \int_{\partial B_r(0)} \left(|\nabla_\tau u|^2- \left|\frac{\partial u}{\partial \nu}\right|^2+\Lambda \chi_{\{u>0\}} - 2u + 2 \frac{u}{r} \frac{\partial u}{\partial \nu}\right)\,d\sigma,
 \end{split}
\end{equation}
where $\nabla_\tau u$ denotes the tangential component of $\nabla u$ on $\partial B_r(0)$ and in the last equality we applied the divergence theorem in $B_r(x_0)$ taking into account that $-\Delta u=1$ in $\{u>0\}$. Now, regrouping the terms in \eqref{eq2:Step5Pteo13} and rescaling we obtain
\begin{equation}\label{eq3:Step5Pteo13}
\begin{split}
\int_{B_1(0)} (|\nabla u_r|^2  &+ \Lambda \chi_{\{u_r>0\}})\, dx
=\displaystyle \frac{2(N+1)}{N} r \int_{B_1(0)} u_r \, dx - \frac{2}{N} r\int_{\partial B_1(0)} u_r \, d\sigma\\
&\qquad\displaystyle +\frac{1}{N} \int_{\partial B_1(0)} \left(|\nabla_\tau u_r|^2- (x\cdot \nabla u_r)^2+\Lambda \chi_{\{u_r>0\}} + 2u_r (x\cdot \nabla u_r)\right)\,d\sigma\,.
\end{split}
\end{equation}
Next, combining \eqref{eq1:Step5Pteo13} and \eqref{eq3:Step5Pteo13} we get that
\begin{equation}\label{eq4:Step5Pteo13}
\begin{split}
 W_\Lambda(z_r)-W_\Lambda(u_r)&=\displaystyle \frac{1}{N} \int_{\partial B_1(0)} \left(u_r- (x\cdot \nabla u_r)^2\right)\, d\sigma - \frac{2(N+1)r}{N}\int_{B_1(0)} u_r \, dx +\frac{2r}{N} \int_{\partial B_1(0)}  u_r\, d\sigma\,.
 \end{split}
\end{equation}
Then, from \eqref{eq4:Step5Pteo13} and \cite[Lemma 9.2]{V} we get that for all $0<r<r_0$ it holds
$$\frac{\partial W_\Lambda(u_r)}{\partial r} = \frac{2}{r} \int_{\partial B_1(0)} \left(u_r- (x\cdot \nabla u_r)\right)^2\, d\sigma - 2(N+1) \int_{B_1(0)} u_r \, dx + 2 \int_{\partial B_1(0)}  u_r\, d\sigma.$$
Finally, since $u$ is Lipschitz in ${B_{r_0}(0)}$, and vanishes at some point $x_1\in\partial B_{r_0}(0)$, we have the following estimate
$$2(N+1) \int_{B_1(0)} u_r \, dx\leq C_1\|\nabla u\|_{L^\infty(B_{r_0}(0))},$$
for some positive constant $C_1$ depending only on $N$, and thus setting $C=C_1\|\nabla u\|_{L^\infty(B_{r_0}(0))}$ we have
$$\frac{\partial W_\Lambda(u_r)}{\partial r} \geq \frac{2}{r} \int_{\partial B_1(0)} \left(u_r- (x\cdot \nabla u_r)\right)^2\, d\sigma - C,$$
and rescaling the integral in the right-hand side we obtain \eqref{eq:Step5formula}. In particular, the function $r\mapsto  W_\Lambda((u_*)_{x_0, r})+ Cr$ is monotone nondecreasing, and as $W_\Lambda((u_*)_{x_0, r})$ is bounded for $r\in]0,r_0[$ (because $u_*$ is Lipschitz in ${B_{r_0}(x_0)}$) we deduce that $W_\Lambda((u_*)_{x_0, r})$ has a finite limit as $r\to 0^+$. The proof of \textit{Step 8} is complete.\vspace{6pt}

In the next crucial step we show that the blow-up limit of $u_*$ at any $x_0\in \partial\Omega_{u_*}\cap \C_\omega$ converges to a 1-homogenous function $u_0$ which is a non-trivial global minimizer of the Alt-Caffarelli functional in $\RN$, as it will be discussed below. More precisely, we have:\vspace{6pt}

\noindent\textbf{Step 9}: Let $x_0\in \partial\Omega_{u_*}\cap \C_\omega$, let $r_0>0$ be such that $B_{r_0}(x_0)\subset\subset \C_\omega$ and consider the rescaled function $(u_*)_{x_0,r}:B_{r_0/r}(0)\to \R$ defined as in \eqref{eq:defblowup}, for $0<r\leq r_0$. There exists a sequence $(r_n)_{n\in\N}$ such that $r_n\to 0$, as $n\to +\infty$, and a non-negative Lipschitz continuous function $u_0:\RN\to \R$ such that $(u_*)_{x_0,r_n}\to u_0$ uniformly on compact subsets of $\RN$, as $n\to +\infty$. Moreover $u_0$ is a (positively) 1-homogenous function and it is a non-trivial global minimizer of the Alt-Caffarelli functional in $\RN$, that means that
\beq\label{eq:thesisStep9ProofTh13}
\int_{B_R(0)}|\nabla u_0|^2\, dx + 2\Lambda |\{u_0>0\} \cap B_R(0)| \leq \int_{B_R(0)}|\nabla w|^2\, dx + 2\Lambda |\{w>0\} \cap B_R(0)|,
\eeq
for all $R>0$ and $w\in H^1_{loc}(\RN)$ such that $w=u_0$ outside $B_R(0)$,  where $\Lambda$ is given by \textit{Step 2} (and it is positive in view of \textit{Step 5}).\vspace{6pt}

The existence of a sequence $(r_n)_{n\in\N}$ and $u_0$ satisfying the properties stated in first part of the statement is a standard fact and the proof can be found for instance in \cite[Proposition 4.25]{LS2} (see also \cite[Sect. 4]{MTV}). For the second part we adapt the proof of \cite[Proposition 4.26]{LS2}.

 In order to ease the notation we write $(u_*)_{x_0,r}$ instead of $(u_*)_{x_0,r_n}$, and take the limit as $r\to 0^+$. For any $R>0$, as in the first step of the proof of \cite[Proposition 4.26]{LS2} (see also the references therein), we have that $(u_*)_{x_0,r}\to u_0$ strongly in $H^1(B_R(0))$ and  $\chi_{\Omega_{(u_*)_{x_0,r}}}\to\chi_{\Omega_{u_0}}$ strongly in $L^1(B_R(0))$, as $r\to 0^+$.  Moreover, thanks to \textit{Step 6}, by a well-known argument, it follows that $\overline\Omega_{(u_*)_{x_0,r}}$, $\Omega_{(u_*)_{x_0,r}}^\complement$ converge, respectively, to $\overline\Omega_{u_0}$, $\Omega_{u_0}^\complement$ in the sense of the Hausdorff distance in $B_R(0)$, as $r\to 0^+$ (see the proofs of \cite[Proposition 4.26]{LS2}, \cite[Proposition 4.5]{MTV}).

 In order to prove that $u_0$ is a global minimizer of the Alt-Caffarelli functional in $\RN$, we fix $R>0$ and $w\in H^1_{loc}(\RN)$ such that $w=u_0$ outside $B_R(0)$. For every $\eps\in ]0,R[$ let $\eta_{\eps}\in C^\infty_c(B_R(0))$ be such that $0\leq \eta_{\eps} \leq 1$, $\eta_{\eps}\equiv 1$ in $B_{R-\eps}(0)$, and define $$w_{r}:=w+(1-\eta_\eps)((u_*)_{x_0,r}-u_0),$$
for $r>0$ sufficiently small such that $\frac{r_0}{r}>R$. By construction $w_{r}=(u_*)_{x_0,r}$ in $B_{\frac{r_0}{r}}(0)\setminus \overline{B_R(0)}$ and considering the blow-down $(w_r)^{r,x_0}$ of $w_r$ centred at $x_0$, i.e. the rescaled function $$(w_r)^{r,x_0}(x)=rw_r\left(\frac{x-x_0}{r}\right),$$ defined for $x\in B_{r_0}(x_0)$, we have that $(w_r)^{r,x_0}=u_*$ in $B_{r_0}(x_0)\setminus\overline{B_{rR}(x_0)}$. Hence, extending $(w_r)^{r,x_0}$ as $u_*$ in $\C_\omega\setminus \overline{B_{r_0}(x_0)}$, we have that $(w_r)^{r,x_0}\in\mathcal{F}$, where $\mathcal{F}=\mathcal{F}(u_*,x_0,r_0)$ is the class of functions defined in \eqref{defsetF}, in addition, we have
$$||\Omega_{(w_r)^{r,x_0}}|-|\Omega_{u_*}||\leq C_2 r^N,$$
for some positive constant $C_2$ independent of $r$. Then, setting $h_r:=\left||\Omega_{(w_r)^{r,x_0}}|-|\Omega_{u_*}|\right|$, it holds that $h_r\to 0^+$, as $r\to 0^+$, and by \textit{Step 4} we infer that
\beq\label{eq1:Step9Pteo13}
J(u_*) + \mu(h_r) \left|\Omega_{u_*}\right|\leq J((w_r)^{r,x_0}) + \mu(h_r) \left|\Omega_{(w_r)^{r,x_0}}\right|,
\eeq
for all sufficiently small $r>0$, for some $\mu(h_r)>0$ such that $\mu(h_r)\to \Lambda$, as $r\to 0^+$. From \eqref{eq1:Step9Pteo13}, by definition of $J$, taking into account that $(w_r)^{r,x_0}=u_*$  in $\C_\omega\setminus \overline{B_{rR}(x_0)}$ we deduce that
\begin{equation}\label{eq2:Step9Pteo13}
\begin{split}
\frac{1}{2}\int_{B_R(0)}|\nabla (u_*)_{x_0,r}|^2\, dx &- r\int_{B_R(0)} (u_*)_{x_0,r}\, dx + \mu(h_r) \left|\Omega_{(u_*)_{x_0,r}}\cap B_R(0)\right|\\
&\leq\displaystyle \frac{1}{2}\int_{B_R(0)}|\nabla w_r|^2\, dx - r\int_{B_R(0)} w_{r}\, dx + \mu(h_r) |\Omega_{w_r}\cap B_R(0)|.
\end{split}
\end{equation}
We now notice that, since $(u_*)_{x_0,r}\to u_0$ strongly in $H^1(B_R(0))$, as $r\to 0^+$, then, by definition of $w_r$, it holds that $w_r\to w$ strongly in $H^1(B_R(0))$, as $r\to 0^+$. Moreover, in view of the inclusion $\Omega_{w_r}\cap B_R(0)\subset \{x\in B_R(0);\ w(x)>0,\, \eta_\eps(x)=1\}\cup\{x\in B_R(0);\  0\leq\eta_\eps(x)<1\}$, we have the bound $$|\Omega_{w_r}\cap B_R(0)|\leq |\Omega_{w}\cap\{x\in B_R(0);\, \eta_\eps(x)=1\}|+|\{x\in B_R(0);\  0\leq\eta_\eps(x)<1\},$$
and since $\chi_{\Omega_{(u_*)_{x_0,r}}}\to\chi_{\Omega_{u_{0}}}$ in $L^1(B_R(0))$ as $r\to 0^+$, passing to the limit in \eqref{eq2:Step9Pteo13}, we get that
$$%\begin{equation}\label{eq3:Step9Pteo13}
\frac{1}{2}\int_{B_R(0)}|\nabla u_{0}|^2\, dx +\Lambda |\Omega_{u_0}\cap B_R(0)|\\
\leq \frac{1}{2}\int_{B_R(0)}|\nabla w|^2\, dx +\Lambda |\Omega_{w}\cap B_R(0)|+ C_3\eps
$$%\end{equation}
for some positive constant $C_3$ independent of $\eps$, because by construction $\eta_\eps\equiv 1$ in $B_{R-\eps}(0)$. Finally, since $\eps>0$ is arbitrary, we obtain the desired inequality \eqref{eq:thesisStep9ProofTh13}. To conclude, it remains to show that $u_0$ is positively 1-homogeneous and non-trivial. By definition we check that, for all $s\in]0,1]$, it holds
 \begin{equation}\label{eq4:Step9Pteo13}
W_\Lambda([(u_*)_{x_0,r}]_{0,s}) = W_\Lambda((u_*)_{x_0, sr}).
\eeq
Now, in view of \textit{Step 8} we infer that there exists $l^+\in \R$ such that, for any $s\in ]0,1]$, $W_\Lambda((u_*)_{x_0, sr})\to l^+$ as $r\to 0^+$ and thus, from \eqref{eq4:Step9Pteo13} we deduce that 
 \begin{equation}\label{eq5:Step9Pteo13}
\lim_{r\to 0^+}W_\Lambda([(u_*)_{x_0,r}]_{0,s})=l^+,~~\forall s\in \left]0,1\right]\,.
\eeq
On the other hand, from the strong convergence of $(u_*)_{x_0,r}$ to $u_0$, in $H^1(B_R(0))$, we also infer that $$\lim_{r\to 0^+}W_\Lambda([(u_*)_{x_0,r}]_{0,s})=W_\Lambda((u_0)_{0, s}),$$
which, together with \eqref{eq5:Step9Pteo13} implies that $[W_\Lambda((u_0)_{0, s})]^\prime(s)=0$ for all $s\in]0,1[$, and from the results of \cite{W1} (see \cite[Remark 4.24]{LS}) this condition implies that $u_0$ is (positively) 1-homogeneous. Finally, from \textit{Step 6}, by a classical argument (see e.g. \cite{MTV,V}), it follows that $u_0$ is non trivial. The proof of \textit{Step 9} is complete.\vspace{6pt}

In order to complete the proof of Theorem \ref{teoreg}, for future convenience let us introduce, similarly to \cite[Definition 6.1]{V}, the following definition:
\begin{definition}\label{def:regsingparts}
 The regular part of the relative boundary $\Gamma_{\Omega^*}= \partial\Omega_{u_*}\cap \C_\omega$ is the set of points $x_0\in \Gamma_{\Omega^*}$ for which the blow-up limit $u_0$, obtained in \textit{Step 9}, is of the form $u_0(x)= \sqrt{2\Lambda} (x\cdot \nu)^+$ for some vector $\nu \in \mathbb{S}^{N-1}$. In that case, $x_0$ is said to be a regular point of $\Gamma_{\Omega^*}$, and we denote by $\mathrm{Reg}(\Gamma_{\Omega^*})$ the regular part of $\Gamma_{\Omega^*}$, i.e., the subset of $\Gamma_{\Omega^*}$ consisting of all its regular points. We define the singular part of $\Gamma_{\Omega^*}$, and denote it by $\mathrm{Sing}(\Gamma_{\Omega^*})$, as the complement in $\Gamma_{\Omega^*}$ of the regular part, namely, $\mathrm{Sing}(\Gamma_{\Omega^*}):=\Gamma_{\Omega^*}\setminus\mathrm{Reg}(\Gamma_{\Omega^*})$.
 \end{definition}

\noindent\textbf{Conclusion of the proof of Theorem \ref{teoreg}}: $\mathrm{Reg}(\Gamma_{\Omega^*})$ is locally a smooth manifold and there exists $d^*\in\{5,6,7\}$ such that:
\begin{itemize}
 \item[(i)] if $N<d^*$, then $\mathrm{Sing}(\Gamma_{\Omega^*})$ is empty;
 \item[(ii)] if $N=d^*$, then $\mathrm{Sing}(\Gamma_{\Omega^*})$ is locally finite;
 \item[(iii)] if $N>d^*$, then $\mathrm{dim}_{\H}(\mathrm{Sing}(\Gamma_{\Omega^*}))\leq N-d^*$.
\end{itemize}
\vspace{6pt}

We begin by proving that $\mathrm{Reg}(\Gamma_{\Omega^*})$ is locally a smooth manifold. Let $x_0 \in \mathrm{Reg}(\Gamma_{\Omega^*})$ and let $r_0>0$ be such that $B_{r_0}(x_0)\subset\subset \C_\omega$. In view of \textit{Step 9}, and arguing as in the proof of \cite[Lemma 5.2]{MTV}, it follows that $u_*$ satisfies the condition $|\nabla u_*|=\sqrt{2\Lambda}$ in the viscosity sense on $\Gamma_{\Omega^*}$ (for the definition, see e.g., \cite[Definition 5.1]{MTV} or \cite[Proposition 4.29]{LS}). Since $u_*$ also satisfies the equation $-\Delta u_*=1$ in $\Omega_{u_*}$ in the weak sense (and thus also in the classical sense), it turns out that $u_*$ is a viscosity solution (in the sense of \cite[Definition 2.2]{DS}) of
$$
\begin{cases}
-\Delta u =1& \hbox{in $\{u>0\}\cap B_{r_0}(x_0)$},\\
|\nabla u|=\sqrt{2\Lambda} & \hbox{on $\partial\{u>0\}\cap B_{r_0}(x_0)$}.
\end{cases}
$$
Considering the rescaled function $\tilde u_r:=(2\Lambda)^{-1/2} (u_*)_{x_0,r}:B_{r_0/r}(0)\to\R$, for all $0<r\leq r_0$, we have that $\tilde u_r$ is, in particular, a viscosity solution to
$$%\beq\label{eq1:ConclusionPrTeo13}
\begin{cases}
-\Delta u =r& \hbox{in $\{u>0\}\cap B_{1}(0)$},\\
|\nabla u|=1 & \hbox{on $\partial\{u>0\}\cap B_{1}(0)$}.
\end{cases}
$$%\eeq
Moreover, by definition of regular point, up to a rigid motion, we can assume without loss of generality that $\nu=\mathbf{e}_{N}$, so that $\tilde u_r\to u_0=x_N^+$ uniformly on compact subsets of $\RN$, as $r\to 0^+$. Let $\bar\eps>0$ be the universal constant given by \cite[Theorem 1.1]{DS}, then, there exists $0<\bar r\leq r_0$ sufficiently small such that $|\tilde u_r(x) - x_N^+|\leq \bar\eps$ for all $x\in B_1(0)$, for any $0<r<\bar r$. This implies that $\tilde u_r$ is $\bar\eps$-flat in $B_1(0)$, namely $(x_N-\bar\eps)^+\leq \tilde u_r(x)\leq (x_N+\bar\eps)^+$, for all $x\in B_1(0)$, for any $0<r<\bar r$. Hence, fixing $0<r<\min\{\bar r, \bar\eps\}$ we get that $\tilde u_r$ satisfies the conditions \cite[(1.2)-(1.3)]{DS} (with $a_{ij}\equiv 1$, $g\equiv 1$, $f\equiv r$ in $B_1$) and thus by \cite[Theorem 1.1]{DS} we obtain that $\partial\{\tilde u_r>0\}\cap B_{1/2}(0)$ is of class $C^{1,\alpha}$, for some $\alpha\in(0,1)$, which in turn implies that $\partial\{u_*>0\}\cap B_{r/2}(x_0)$ is of class $C^{1,\alpha}$. Then, by a classical argument (see \cite{KN}) we conclude that $\partial\{u_*>0\}\cap B_{r/2}(x_0)$ is of class $C^{\infty}$ and this proves that $\mathrm{Reg}(\Gamma_{\Omega^*})$ is a smooth manifold.

Finally, we prove that, for any open bounded subset $D\subset\subset \C_\omega$ such that $D\cap \Gamma_{\Omega_*}\neq\emptyset$, $u_*$ satisfies the hypotheses \textit{(a)--(d)} of \cite[Proposition 10.13]{V} in $D$, and by the arbitrariness of $D$ we conclude. 

To this end, fixing $D\subset\subset \C_\omega$, since $\overline{D}$ is a compact subset of $\C_\omega$, arguing as above, we find $\bar\eps>0$ (the universal constant given by  \cite[Theorem 1.1]{DS}) and $R>0$, independent of $x_0$, such that if $x_0\in D\cap \Gamma_{\Omega_*}$, $r\in (0,R)$ are such that $B_r(x_0)\subset D$ and $\|u_*(x) - \sqrt{2\Lambda} ((x-x_0)\cdot \nu)^+\|_{L^\infty(B_r(x_0))}\leq \bar\eps r$, for some $\nu \in \mathbb{S}^{N-1}$, then $\Gamma_{\Omega_*}\cap B_{r/2}(x_0)\subset \mathrm{Reg}(\Gamma_{\Omega^*})$. This means that $u_*$ satisfies the  property \emph{(a) $\eps$-regularity}. Moreover, up to taking a finite subcover of $\overline{D}$, made by open balls strictly contained in $\C_\omega$, from \textit{Step 6} we infer that there exist constants $\kappa>0$ and $\bar r>0$ such that if $x_0\in \Gamma_{\Omega^*}\cap D$ and $r\in (0,\bar r)$ such that $B_r(x_0)\subset D$ then $\|u\|_{L^\infty(B_r(x_0))}\geq \kappa r$. In other words $u_*$ satisfies the  property \emph{(b) Non-degeneracy}. Finally, by \textit{Step 9}, it also follows that $u_*$ satisfies the hypotheses \emph{(c) Convergence of the blow-up sequences} and \emph{(d) Homogeneity and minimality of the blow-up limit}. Thus, we can apply \cite[Proposition 10.13]{V} and the proof is complete.
\end{proof}

We conclude this section with the following remarkable topological property.
\begin{theorem}\label{connected}
If $\Omega^*$ is a minimizer for $\E(\cdot;\C_\omega)$ in $\A_{\omega,c}$, for some $c>0$, then $\Omega^*$ is connected.
\end{theorem}
\begin{proof}
Let $\Omega^*\in \A_{\omega,c}$ be a minimizer for $\E(\cdot;\C_\omega)$. For simplicity in notation, we will omit the superscript $*$ and simply write $\Omega$ in place of $\Omega^*$.  By Proposition \ref{prop1:sect5} and Proposition \ref{prop2:sect5}, we know that $\Omega$ is a bounded open subset of $\C_\omega$. Assume by contradiction that $\Omega$ is not connected. Then there exist nonempty open sets $\Omega_{1}$ and $\Omega_{2}$ such that $\Omega=\Omega_{1}\cup\Omega_{2}$ and $\Omega_{1}\cap\Omega_{2}=\varnothing$. Setting  $c_i:=|\Omega_i|>0$, for $i=1,2$ we have
\beq\label{eq-1:teo1sect5}
0< c_1<c,\ \ 0<c_2<c, \ \ c=c_1+c_2.
\eeq
Moreover, as $\partial\Omega=\partial\Omega_1\cup\partial\Omega_2$ we have $\Gamma_{\Omega}=\Gamma_{\Omega_1} \cup \Gamma_{\Omega_2}$ where $\Gamma_{\Omega_{i}}=\partial\Omega_{i}\cap\C_{\omega}$. Since $\Omega$ is a minimizer for $\E(\cdot;\C_\omega)$ in $\A_{\omega,c}$, then by Theorem \ref{teoreg} we infer that $\H^{N-2}(\mathrm{Sing}(\Gamma_{\Omega}))=0$, and thus, by \cite[Theorem 4.16]{EG} we get that 
\begin{equation}
\label{null-capacity}
\mathrm{cap}(\mathrm{Sing}(\Gamma_{\Omega}))=0,
\end{equation} 
where $\mathrm{Sing}(\Gamma_{\Omega})$ is the singular part of  $\Gamma_{\Omega}$ (see Definition \ref{def:regsingparts}).
In addition, by standard elliptic regularity theory, taking into account of Theorem \ref{teoreg}, we deduce that the energy function $u_\Omega \in H_0^1(\Omega; \C_\omega)$ corresponding to $\Omega$ is smooth in $\Omega\cup \mathrm{Reg}(\Gamma_\Omega)$ and has null trace on $\mathrm{Reg}(\Gamma_\Omega)$ in the classical sense. For $i=1,2$ we set
$$
u_{i}:=\begin{cases}u_{\Omega}&\text{in $\Omega_{i}$}\\ 0&\text{in $\C_{\omega}\setminus\Omega_{i}$.}\end{cases}
$$
\textbf{Claim I}: $u_{i}\in H^{1}_{0}(\Omega_{i};\C_{\omega})\quad(i=1,2)$.\vspace{6pt}

This claim seems not to be so obvious because we only have partial information on the regularity of $\partial\Omega$, in  particular on $\partial\Omega\cap\partial\C_\omega$. In order to prove it we construct a suitable sequence $(u_{n,i})_{n}\subset H^{1}_{0}(\Omega_{i};\C_{\omega})$ such that $u_{n,i}\to u_{i}$ in $H^{1}(\C_{\omega})$ and each $u_{n,i}$ is continuous in $\C_{\omega}$ for every $n$. More precisely, thanks to \eqref{null-capacity}, arguing as in the proof of Lemma \ref{lem:tech2} and, up to mollification, we find a sequence $(v_{n})_n\subset C_c^\infty(\RN)$ such that 
\begin{equation}\label{eq:propvn-1}
\begin{split}\|v_{n}\|_{H^{1}(\R^N)}\to 0,~~ v_{n}\to 0 \text{~~a.e. in $\R^N$,}\\0\le v_{n}\le 1\text{ in $\R^N$,~~$v_{n}=1$ in $\Lambda_{n}$,}
\end{split}\end{equation}
where $\Lambda_{n}$ is a neighborhood of $\mathrm{Sing}(\Gamma_{\Omega})$. Let us remark that
\begin{equation}\label{Lambdan}
|\overline{\Lambda_{n}}|\to 0\quad\text{as }n\to+\infty
\end{equation}
because
$$
|\overline{\Lambda_{n}}|=\int_{\overline{\Lambda_{n}}}v_{n}^{2}\,dx\le\int_{\R^{N}}v_{n}^{2}\,dx\to 0\quad\text{as }n\to+\infty.
$$

For every $n\in\N$, the function 
$$
u_{n}:=(1-v_{n})u_{\Omega}
$$
belongs to $H^{1}(\C_{\omega})$ and is continuous in $\C_{\omega}$ because $u_{\Omega}$ is so in $(\Omega\cup\mathrm{Reg}(\Gamma_\Omega)\cup(\C_{\omega}\setminus\overline\Omega)$ and $v_{n}=1$ in a neighborhood of $\mathrm{Sing}(\Gamma_\Omega)$. In particular $u_{n}=0$ on $\Lambda_{n}\cup\mathrm{Reg}(\Gamma_\Omega)$. We consider the open set
$$\Omega_{n,i}:= \Omega_{i}\setminus\overline{\Lambda_{n}}$$
%\Omega_{i}\setminus\overline{\Lambda_{n}}\subset\Omega_{n,i}\subset\Omega_{i}\text{~~with~~}\partial\Omega_{n,i}\cap\C_{\omega}\text{~~smooth.}
%$$
%Since $\Gamma_{\Omega}^{\mathrm{reg}}$ is smooth, for every $n\in\N$ and for $i=1,2$, there exists an open set $\Omega_{n,i}$ such that
%$$
%\Omega_{i}\setminus\overline{\Lambda_{n}}\subset\Omega_{n,i}\subset\Omega_{i}\text{~~with~~}\partial\Omega_{n,i}\cap\C_{\omega}\text{~~smooth.}
%$$
and we define
$$
u_{n,i}:=\begin{cases}u_{n}&\text{in $\Omega_{n,i}$}\\ 0&\text{in $\C_{\omega}\setminus\Omega_{n,i}$.}\end{cases}
$$
Let us point out that $u_{n,i}$ is continuous in $\C_{\omega}$ because it is so in $\Omega_{n,i}$, in $\C_{\omega}\setminus\Omega_{n,i}$, and $u_{n}$ is continuous in $\C_{\omega}$ with $u_{n}=0$ on $\partial\Omega_{n,i}\cap\C_{\omega}$. Let us check that
\begin{itemize}
\item[(a)] $u_{n,i}\in H^{1}(\C_{\omega})$;
\item[(b)] $u_{n,i}\in H^{1}_{0}(\Omega_{i};\C_{\omega})$;
\item[(c)] $u_{n,i}\to u_{i}$ as $n\to+\infty$, in $H^{1}(\C_{\omega})$.
\end{itemize}
Let us prove (a): by definition of $u_{n,i}$ and since $u_{n}\in L^{2}(\C_{\omega})$, also $u_{n,i}\in L^{2}(\C_{\omega})$. Fixing $j\in\{1,...,N\}$, we set
$$
g_{n,i,j}:=\begin{cases}\frac{\partial u_{n}}{\partial x_{j}}&\text{in $\Omega_{n,i}$}\\ 0&\text{in $\C_{\omega}\setminus\Omega_{n,i}$.}\end{cases}
$$
Since $u_{n}\in H^{1}(\C_{\omega})$, we have that $g_{n,i,j}\in L^{2}(\C_{\omega})$. Let us check that $g_{n,i,j}$ is the weak derivative (with respect to $x_{j}$) of $u_{n,i}$.

Indeed, fixing $\varphi\in C^{\infty}_{c}(\C_{\omega})$, since $\mathrm{dist}(\mathrm{supp}\,\varphi,\partial\C_{\omega})>0$, there exists a domain $\omega'\subset\R^{N-1}$ such that $\overline{\omega'}\subset\omega$ and $\mathrm{supp}\,\varphi\subset\C_{\omega'}:=\omega'\times\R$. Then, letting $\hat x_{j}=(x_{1},...,x_{j-1},x_{j+1},...,x_{N})$, by Fubini's theorem and since for a.e. $\hat x_j$ the mapping $x_j\mapsto u_n (x_1,\ldots, x_j,\ldots,x_N)$ is absolutely continuous in $\{x_j\in\R;\, (x_1,\ldots, x_j,\ldots,x_N)\in \C_{\omega'}\cap\Omega_{n,i}\}$, which is an open set of $\R$ (as $\C_{\omega'}\cap\Omega_{n,i}$ is an open set of $\R^N$), integrating by parts with respect to $x_j$, we infer that
$$
\int_{\C_{\omega}}\hspace{-3pt}\varphi\,g_{n,i,j}\,dx%=\int_{\C_{\omega'}\cap\Omega_{n,i}}\hspace{-9pt}\varphi\,\frac{\partial u_{n}}{\partial x_{j}}\,dx
=\int_{\C_{\omega'}\cap\Omega_{n,i}}\hspace{-9pt}\varphi\,\frac{\partial u_{n}}{\partial x_{j}}\,dx_{j}\,d\hat x_{j}=-\int_{\C_{\omega'}\cap\Omega_{n,i}}\hspace{-9pt}u_{n}\,\frac{\partial\varphi }{\partial x_{j}}\,dx_{j}\,d\hat x_{j}=-\int_{\C_{\omega}}\hspace{-3pt}u_{n,i}\,\frac{\partial\varphi }{\partial x_{j}}\,dx
$$
because $\partial(\C_{\omega'}\cap\Omega_{n,i})\subset\partial\C_{\omega'}\cup(\partial\Omega_{n,i}\cap\C_{\omega})$ and $\varphi u_{n}=0$ on $\partial\C_{\omega'}\cup(\partial\Omega_{n,i}\cap\C_{\omega})$. Hence $g_{n,i,j}$ is the weak derivative (with respect to $x_{j}$) of $u_{n,i}$ and (a) holds true. 
\medskip

Then, as $u_{n,i}\in H^{1}(\C_{\omega})$ and since by construction $u_{n,i}$ is continuous in $\C_\omega$ and $u_{n,i}=0$ in $\C_{\omega}\setminus\Omega_{n,i}$, and $\C_{\omega}\setminus\Omega_{i} \subset \C_{\omega}\setminus\Omega_{n,i}$,  we readily deduce (b). 
\medskip

Let us prove (c). First, we have that
\begin{equation*}
\begin{split}
\int_{\C_{\omega}}|u_{n,i}-u_{i}|^{2} \, dx&=\int_{\Omega_{i}}|u_{n,i}-u_{\Omega}|^{2} \, dx =\int_{\Omega_{n,i}}|(1-v_{n})u_{\Omega}-u_{\Omega}|^{2} \, dx+\int_{\Omega_{i}\setminus\Omega_{n,i}}u_{\Omega}^{2}  \, dx\\
&=\int_{\Omega_{n,i}}v_{n}^{2}u_{\Omega}^{2} \, dx+\int_{\overline{\Lambda_n}}u_{\Omega}^{2}  \, dx.
\end{split}
\end{equation*}
Since $u_\Omega$ is bounded in $\C_\omega$ (see Proposition \ref{prop1:sect2}), and thanks to \eqref{eq:propvn-1}, we get that
$$
\int_{\Omega_{n,i}}v_{n}^{2}u_{\Omega}^{2} \, dx\le\|u_{\Omega}\|_{\infty}^{2}\int_{\R^{N}}v_{n}^{2}\, dx\to 0\quad\text{as~~}n\to\infty.
$$
Moreover, by \eqref{Lambdan} and since $u_{\Omega}\in L^{2}(\C_{\omega})$, we have
$$
\int_{\overline{\Lambda_{n}}}u_{\Omega}^{2}\, dx\to 0\quad\text{as~~}n\to\infty.
$$
We now handle the gradient term. To this end we write
$$
\int_{\C_{\omega}}|\nabla(u_{n,i}-u_{i})|^{2} \, dx=\int_{\Omega_{i}}|\nabla(u_{n,i}-u_{\Omega})|^{2} \, dx=\underbrace{\int_{\Omega_{n,i}}|\nabla(v_{n}u_{\Omega})|^{2} \, dx}_{\mathbf{I_1}}+\underbrace{\int_{\Omega_{i}\setminus\Omega_{n,i}}|\nabla u_{\Omega}|^{2} \, dx}_{\mathbf{I_2}}.
$$
For $\mathbf{I_1}$ we have
$$
\int_{\Omega_{n,i}}|\nabla(v_{n}u_{\Omega})|^{2} \, dx\le 2\int_{\Omega_{n,i}}|\nabla v_{n}|^{2}\,u_{\Omega}^{2} \, dx+2\int_{\Omega_{n,i}}v_{n}^{2}\,|\nabla u_{\Omega}|^{2} \, dx.
$$
Now, from \eqref{eq:propvn-1} we readily deduce that
$$
\int_{\Omega_{n,i}}|\nabla v_{n}|^{2}\,u_{\Omega}^{2} \, dx\le\|u_{\Omega}\|_{\infty}^{2}\int_{\R^{N}}|\nabla v_{n}|^{2} \, dx\to 0,
$$
as $n\to +\infty$. Moreover, since $\Omega_{n,i}\subset\Omega_{i}$, $v_n\to 0$ a.e. in $\R^N$, as $n\to +\infty$, and $0\leq v_n\leq1$, then, by the dominated convergence theorem, applied in $\Omega_i$, we infer that 
$$
\int_{\Omega_{n,i}}v_{n}^{2}\,|\nabla u_{\Omega}|^{2} \, dx\le\int_{\Omega_i}v_n^2|\nabla u_{\Omega}|^{2} \, dx\to 0.
$$
%by \eqref{Kappan}. 
Finally, for $\mathbf{I_2}$, exploiting \eqref{Lambdan} and since $u_{\Omega}\in H^{1}(\C_{\omega})$, as $n\to +\infty$, we conclude that
$$
\int_{\Omega_{i}\setminus\Omega_{n,i}}|\nabla u_{\Omega}|^{2}\, dx=\int_{\overline{\Lambda_{n}}}|\nabla u_{\Omega}|^{2}\, dx\to 0.
$$
Hence (c) holds true, and \textit{Claim I} follows immediately from (a), (b) and (c).
\bigskip

Thanks to $\textit{Claim I}$, and by definition of torsional energy, we have $ \E(\Omega_{i};\C_{\omega}) \leq J_{\Omega_{i}}(u_{i})$, for $i=1,2$. Then, from this and by construction, we infer that
$$
\E(\Omega_{1};\C_{\omega})+\E(\Omega_{2};\C_{\omega})\le J_{\Omega_{1}}(u_{1})+ J_{\Omega_{2}}(u_{2})=J_{\Omega_{1}\cup\Omega_{2}}(u_\Omega)=J_{\Omega}(u_\Omega)=\E(\Omega;\C_{\omega}),
$$
which, together with \eqref{subadditivity-energy}, implies that

\beq\label{eq0:teo1sect5}
\E(\Omega_1; \C_\omega) + \E(\Omega_2; \C_\omega)= \E(\Omega; \C_\omega)=\O_c(\C_\omega).
\eeq
We now prove the following.
\medskip

\noindent\textbf{Claim II}: $\Omega_i$ is a minimizer for $\E(\cdot;\C_{\omega})$ in $\A_{\omega,c_{i}}$, for $i=1,2$.
\medskip

Indeed, on the contrary, assume for instance that $\Omega_1$ is not a minimizer of $\E(\cdot;\C_{\omega})$ in $\A_{\omega,c_{1}}$, then, denoting by $\widetilde\Omega_1$ a minimizer (which always exists in view of Theorem \ref{teo:mainteoexistmin}) we have
\beq\label{eq00:teo1sect5}
\E(\widetilde\Omega_1; \C_\omega)< \E(\Omega_1; \C_\omega).
\eeq
By Proposition \ref{prop1:sect5} and Corollary \ref{cor1:sect5} we have that $\widetilde\Omega_1$ is a bounded open set and $|\widetilde\Omega_1|=c_1$.
Since $\widetilde\Omega_1$ and $\Omega_2$ are bounded, then, up to a translation along the axial $x_N$-direction of $\widetilde\Omega_1$ we can assume without loss of generality that $\overline{\widetilde\Omega_1}\cap\overline{\Omega_2}=\varnothing$. Then, applying Lemma \ref{additivity} and exploiting \eqref{eq00:teo1sect5}, we get that
\beq\label{eq000:teo1sect5}
\E(\widetilde\Omega_1 \cup\Omega_2; \C_\omega)=\E(\widetilde\Omega_1; \C_\omega) + \E(\Omega_2; \C_\omega) < \E(\Omega_1; \C_\omega) + \E(\Omega_2; \C_\omega).
\eeq
Finally, since $|\widetilde\Omega_1 \cup\Omega_2|=c$ then from \eqref{eq0:teo1sect5} and \eqref{eq000:teo1sect5} we contradict the minimality of $\Omega$. The proof of \textit{Claim II} is complete.
\medskip

Next we prove the following.
\medskip

\noindent\textbf{Claim III}: For $i=1,2$ it holds that
$$
\frac{\E(\Omega_i; \C_\omega)}{|\Omega_i|}>  \frac{\E(\Omega; \C_\omega)}{|\Omega|}.
$$

We first notice that by Lemma \ref{lem1:sect5} we have $\frac{\E(\Omega_i; \C_\omega)}{|\Omega_i|}\geq  \frac{\E(\Omega; \C_\omega)}{|\Omega|}$ for $i=1,2$. To prove the strict inequality, for $i=1,2$ we set $t_i:=\frac{|\Omega|}{|\Omega_i|}$. By \eqref{eq-1:teo1sect5} we have that $t_i>1$, for $i=1,2$. Setting
$$ 
u_i:= u_{\Omega_i}\circ F_{t_i}^{-1},
$$ 
where $u_{\Omega_i}$ is the energy function of $\Omega_i$ and $F_{t_i}$ is the diffeomorphism defined in \eqref{eq:defFt} (with $t=t_i$), then, arguing as in \eqref{eq3:prop1:sect5}--\eqref{eq6:prop1:sect5} we infer that
\beq\label{eq1a:teo1sect5}
\begin{array}{lll}
\displaystyle t_i \E(\Omega_i; \C_\omega)
&=&\displaystyle \int_{F_{t_i}(\Omega_i)}\left[\frac{1}{2}\left( \left|\nabla^\prime  u_i \right|^2 + t_i^2 \left(\frac{\partial   u_i }{\partial x_N}\right)^2\right) -   u_i \right]\, dx\\[15pt]
&=&\displaystyle \int_{F_{t_i}(\Omega_i)}\left(\frac{1}{2} \left|\nabla  u_i \right|^2  -   u_i \right)\, dx +  \frac{t_i^2-1}{2}\int_{F_{t_i}(\Omega_i)}  \left(\frac{\partial   u_i }{\partial x_N}\right)^2\, dx\\[15pt]
&=& \displaystyle J_{F_{t_i}(\Omega_i)} ( u_i) +  \frac{t_i^2-1}{2}\int_{F_{t_i}(\Omega_i)}  \left(\frac{\partial   u_i }{\partial x_N}\right)^2\, dx\\[9pt]
&\geq&\displaystyle  \E(F_{t_i}(\Omega_i); \C_\omega) +  \frac{t_i^2-1}{2}\int_{F_{t_i}(\Omega_i)}  \left(\frac{\partial   u_i }{\partial x_N}\right)^2\, dx.
\end{array}
\eeq
Hence, as $|F_{t_i}(\Omega_i)|=|\Omega|$ (see \eqref{eq1:lem2sect5}), from \eqref{eq1a:teo1sect5} we obtain
$$
t_i \E(\Omega_i; \C_\omega) \geq \E(\Omega; \C_\omega) +  \frac{t_i^2-1}{2}\int_{F_{t_i}(\Omega_i)}  \left(\frac{\partial   u_i }{\partial x_N}\right)^2\, dx.
$$
To conclude, as $t_i>1$, it is sufficient to show that 
\begin{equation}\label{positive}
\int_{F_{t_i}(\Omega_i)}  \left(\frac{\partial   u_i }{\partial x_N}\right)^2\, dx=\frac1{t_{i}}\int_{\Omega_{i}}\left(\frac{\partial   u_{\Omega_{i}} }{\partial x_N}\right)^2\, dx>0.
\end{equation}
Indeed, since $\Omega_i$ is an open subset of $\C_\omega$, $u_{\Omega_i}\in H_0^1(\Omega_i; \C_\omega)$ is a weak solution to \eqref{eq:mixbvprobquasiopen} (with $\Omega=\Omega_i$), then by standard elliptic regularity theory $u_{\Omega_i} \in C^2(\Omega_i)$ and, by the strong maximum principle, 
\begin{equation}\label{ui-positive}
u_{\Omega_{i}}>0\text{~~in~~}\Omega_{i}.
\end{equation} 
Now, assume by contradiction that \eqref{positive} is false. Then
\begin{equation}\label{contr}
\frac{\partial   u_{\Omega_{i}} }{\partial x_N}(x)=0~~\forall x\in\Omega_{i}.
\end{equation}Since $\Omega_{i}$ is bounded, for every $x=(x',x_{N})\in\Omega_{i}$ there exists $\tau(x')>1$ such that 
\begin{equation}
\label{tildex}
\widetilde{x}:=(x',\tau(x')x_{N})\in\partial\Omega_{i}\quad\text{and}\quad (x',tx_{N})\in\Omega_{i}~~\forall t\in\left[1,\tau(x')\right[.
\end{equation} 
More precisely, $\widetilde{x}\in\Gamma_i:=\partial\Omega_i \cap \C_\omega$. Since, by \textit{Claim II}, $\Omega_i$ is a minimizer of $\E(\cdot;\C_\omega)$ in $\A_{\omega,c_{i}}$, by Theorem \ref{teoreg}, we infer that $\H^{N-1}(\mathrm{Sing}(\Gamma_i))=0$. Therefore we can always find $x\in\Omega_{i}$ such that $\widetilde{x}\in \mathrm{Reg}(\Gamma_i)$, where $\widetilde{x}$ is defined in \eqref{tildex}. Then, by \eqref{contr}, $u_{i}(tx+(1-t)\widetilde{x})=0$ for every $t\in[0,1]$, in contrast with \eqref{ui-positive}. The proof of \textit{Claim III} is complete.
\medskip

We now conclude the proof of the theorem. In view of \textit{Claim III}, and exploiting \eqref{eq-1:teo1sect5}--\eqref{eq0:teo1sect5}, we can write
\begin{gather*}
\frac{\E(\Omega_1; \C_\omega)}{c_1}>  \frac{\E(\Omega_1; \C_\omega)+\E(\Omega_2; \C_\omega)}{c_1+c_2},\\
\frac{\E(\Omega_2; \C_\omega)}{c_2}> \frac{\E(\Omega_1; \C_\omega)+\E(\Omega_2; \C_\omega)}{c_1+c_2}.\\
\end{gather*}
Then, isolating $\E(\Omega_1; \C_\omega)$ and $\E(\Omega_2; \C_\omega)$ in each inequality, by elementary algebraic computations we readily infer that
\begin{gather}
\E(\Omega_1; \C_\omega)\,\frac{c_2}{c_1}>  \E(\Omega_2; \C_\omega),\label{eq0:cstep1teo1sect5}\\
\E(\Omega_2; \C_\omega)> \frac{c_2}{c_1}\,\E(\Omega_1; \C_\omega).\label{eq1:cstep1teo1sect5}
\end{gather}
Hence, combining \eqref{eq0:cstep1teo1sect5}--\eqref{eq1:cstep1teo1sect5} we obtain
$$\E(\Omega_1; \C_\omega)>\E(\Omega_1; \C_\omega),$$
which is clearly a contradiction. The proof is complete.
\end{proof}

We conclude this section with the proof of Theorem \ref{mainteo1}.
\begin{proof}[Proof of Theorem \ref{mainteo1}]
Let $\omega$ be any bounded domain in $\R^{N-1}$ with Lipschitz boundary and let $c>0$. The existence of a minimizer $\Omega^*$ of $\E(\cdot,\C_\omega)$ in $\A_{\omega,c}$ is given by Theorem \ref{teo:mainteoexistmin}. By Proposition \ref{prop1:sect5}, Proposition \ref{prop2:sect5}, and Corollary \ref{cor1:sect5}, we know that $\Omega^*$ is a bounded open set and $|\Omega^*|=c$. Finally, in view Theorem \ref{connected}, $\Omega^*$ is connected. The proof is complete.
\end{proof}

\begin{remark}\label{R:CF}
If $\Omega^{*}$ is a minimizer of $\E(\cdot,\C_\omega)$ in $\A_{\omega,c}$, for some $c>0$, and its relative boundary $\Gamma_{\Omega^{*}}$ is almost nowhere parallel to the $x_{N}$-axis inside $\C_{\omega}$ in the sense of \cite[{(2.11)}]{CF}, then $\Omega^{*}$ is convex in the $x_{N}$-direction and symmetric with respect to some horizontal hyperplane $x_{N}=a$ for some $a\in\R$. This is obtained by applying \cite[Theorem 2.2]{CF}.
\end{remark}

\section{Related overdetermined problems and estimates on the free boundary}\label{S:odtpb}
We begin this section by showing that on the regular part of the relative boundary of a minimizer $\Omega^*$ the normal derivative of the corresponding energy function $u_{\Omega^*}$ is constant and we determine a bound on its value. 

%Firstly, we clarify the notion of regular part of the boundary. 

%\begin{definition}\label{def:regular-part}
%Let $\Omega\subset\R^{N}$ be an open set. The \emph{regular part of $\partial\Omega$}, denoted $\partial\Omega^{\mathrm{reg}}$, is the set of points $x\in\partial\Omega$ such that there exist $r>0$ and $\psi\in C^{1}(B_{r}(x))$ with the property that $B_{r}(x)\cap\Omega=\{\psi<0\}$ and $B_{r}(x)\cap\partial\Omega=\{\psi=0\}$. In this case, for every $x\in \partial\Omega^{\mathrm{reg}}$, the versor $\nu(x):=\frac{\nabla\psi(x)}{|\nabla\psi(x)|}$ is the outer unit normal to $\Omega$ at $x$. %Moreover, we say that the set $\Omega$ has \emph{almost $C^{k,\alpha}$-boundary} (or \emph{almost $C^{k}$-boundary}) if $\mathcal{H}^{N-1}(\partial\Omega\setminus\partial\Omega^{\mathrm{reg}})=0$ and $\partial\Omega^{\mathrm{reg}}$ is of class $C^{k,\alpha}$ (or $C^{k}$), namely the functions $\psi$ are so. 
%\end{definition}

\begin{proposition}\label{lem2:sect5}
Let $c>0$, let $\omega\subset\R^{N-1}$ be a Lipschitz bounded domain.
Let $\Omega^*$ be a minimizer for $\E(\cdot;\C_\omega)$ in $\A_{\omega,c}$, let $\Gamma_{\Omega^*}=\partial\Omega^*\cap\C_\omega$ be its relative boundary and let $u_{\Omega^*}\in H_0^1(\Omega^*; \C_\omega)$ be the energy function of $\Omega^*$. Then there exists a positive constant $C_0$ such that
$$%\beq\label{eq:dercostregpart}
\left(\frac{\partial u_{\Omega^*}}{\partial \nu}\right)^2=|\nabla u_{\Omega^*}|^2\equiv C_0 \  \ \hbox{on $\mathrm{Reg}(\Gamma_{\Omega^*})$,}
$$%\eeq
where  $\mathrm{Reg}(\Gamma_{\Omega^*})$  is the regular part of  $\Gamma_{\Omega^*}$ (see Definition \ref{def:regsingparts}). Moreover it holds that
\beq\label{eq:boundlemmasect5}
C_0\geq \frac{2|\O_c(\C_\omega)|}{c}\,.
\eeq
\end{proposition}
\begin{proof}
Let $\Omega^*$ and $\Gamma_{\Omega^*}$ be as in the statement. To ease the notation we write $\Omega$, $\Gamma$ instead of $\Omega^*$, $\Gamma_{\Omega^*}$, respectively. By the proof of Theorem \ref{teoreg} we know that $\mathrm{Reg}(\Gamma)$ is locally a smooth manifold.
%
%we know that there exists a critical dimension $d^*\in\{5,6,7\}$, such that: $\Gamma$ is a smooth manifold if $N<d^*$; $\Gamma$ can have countable isolated singularities if $N=d^*$, and $\Gamma$ can have a singular set of dimension $N-d^*$, if $N>d^*$.
%
%In particular, for any dimension $N\geq 2$ the regular part of $\Gamma$, denoted $\Gamma^{\textrm{reg}}$, is a nonempty, relative open subset of $\Gamma$.
 Let $x_0\in\mathrm{Reg}(\Gamma)$, and let $B_r(x_0)$ be a small ball such that ${B_{2r}(x_0)}\subset \C_\omega$, $\Gamma\cap{B_{2r}(x_0)} \subset \mathrm{Reg}(\Gamma)$ and the outer normal $\nu$ to $\Gamma\cap{B_{2r}(x_0)}$ admits an extension to a smooth vector field $\overline\nu$ defined in $\overline{B_r(x_0)}$. We also observe that by standard elliptic regularity theory the energy function $u_\Omega$ is smooth, as well as its derivatives, on $\Omega\cup(\Gamma\cap{B_{2r}(x_0)})$, in particular, $\left(\frac{\partial u_{\Omega}}{\partial \nu}\right)^2\big\vert_{\Gamma\cap\overline{B_{r}(x_0)}}=|\nabla u_\Omega|^2\big\vert_{\Gamma\cap\overline{B_{r}(x_0)}}$ extends to a smooth function defined in $\overline{B_r(x_0)}$, denoted by $\overline{|\nabla u_\Omega|^2}$.
\medskip

\noindent\textbf{Claim I:} $|\nabla u_\Omega|^2$ is constant on $\Gamma\cap \overline{B_r(x_0)}$.
\medskip

\noindent
Assume by contradiction that the claim is false. Then we find a relative open set $S\subset \subset \Gamma\cap \overline{B_r(x_0)}$ such that $|\nabla u_\Omega|^2$ is not constant on $S$. Let $\varphi \in C_c^\infty(B_r(x_0))$ be a cut-off function such that $\varphi\equiv 1$ in a neighborhood of $\overline{S}$, let us set 
\beq\label{eq:defCvarphi}
C_{\varphi}:=\frac{\int_{\Gamma\cap \overline{B_r(x_0)}}  \varphi |\nabla u_\Omega|^2\, d\sigma}{\int_{\Gamma\cap \overline{B_r(x_0)}}  \varphi\, d\sigma}
\eeq
%where $d\sigma$ denotes the $(N-1)$-dimensional surface area element, 
and consider the smooth vector field $V\colon\R^N\to\R^N$, defined by 
\beq\label{eq:defV}
V=\begin{cases}\varphi (\overline{|\nabla u_\Omega|^2}-C_\varphi)\overline{\nu} & \hbox{in $B_r(x_0)$,}\\[2pt] 0 & \hbox{in $\R^{N}\setminus B_r(x_0)$.}\end{cases}
\eeq
Let $\xi\colon\R\times\R^{N}\to \R^{N}$ be the associated flow, i.e., the solution to
$$
\begin{cases}\dfrac{\partial\xi}{\partial t}(t,x)=V(\xi(t,x))\\[4pt] \xi(0,x)=x.\end{cases}
$$
Since $V(x)=0$ for $x\in\R^{N}\setminus B_{r}(x_0)$, we have that $\xi(t,x)=x$ for every $(t,x)\in\R\times(\R^{N}\setminus B_{r}(x_0))$. In particular, since $B_{r}(x_0)\subset \C_\omega$, $\xi(t,x)=x$ for every $(t,x)\in\R\times\partial\C_{\omega}$ and $\xi(t,x)\in\C_{\omega}$ for every $(t,x)\in\R\times\C_{\omega}$. Moreover, for every $t\in\R$ the set $\Omega_t:=\xi(t,\Omega)$ is an open subset of $\C_{\omega}$. 
From \cite[Sect. 5.9.3, (5.103)]{HP} (see also \cite[Proposition 7.4]{IPW} and \cite[Proposition 4.3]{PT2}) we know that the function $t\mapsto \mathcal{E}(\Omega_t;\C_\omega)= -\frac{1}{2}\int_{\Omega_t}|\nabla u_{\Omega_t}|^2 \ dx$ is differentiable at $t=0$ and since $V$ has compact support in a neighborhood of $x_0\in  \mathrm{Reg}(\Gamma)$ we have
\begin{equation}\label{eq:diffTorsionalOmegat} 
\begin{array}{lll}
\displaystyle \frac{d}{dt}\left.\left( \mathcal{E}(\Omega_t;\C_\omega)\right)\right|_{t=0} &=&\displaystyle - \frac{1}{2}\int_{\Gamma\cap \overline{B_r(x_0)}}  |\nabla u_\Omega|^2 \langle V,\nu\rangle \ d\sigma.
\end{array}
\end{equation}
Similarly, for the volume (see \cite[Sect. 5.9.3]{HP}) we have
\begin{equation}\label{eq:diffVolOmegat}
 \frac{d}{dt}\left.\left(|\Omega_t|\right)\right|_{t=0} =\int_{\Gamma\cap \overline{B_r(x_0)}}\langle V,\nu\rangle \ d\sigma.%=\int_{\Gamma\cap \overline{B_r(x_0)}} \psi \ d\sigma.
\end{equation}
However, in general, the deformed sets $\Omega_t=\xi(t,\Omega)$ do not preserve the volume. Nevertheless, arguing as in the proof of \cite[Theorem 2.2]{SZ}, we can prove that there exist $t_1>0$ and a volume-preserving deformation of $\Omega$, obtained by suitably modifying $(\Omega_t)_{t\in (-t_1, t_1)}$. For the reader's convenience we provide some details of the proof.

Up to choosing a smaller radius $r>0$ at the beginning of the proof we can assume without loss of generality that $\Gamma\cap \overline{B_r(x_0)}$ is the graph of a smooth function and that the same property holds for 
$$
\Gamma_t\cap \overline{B_r(x_0)}:=\partial\Omega_t\cap \overline{B_r(x_0)}\cap\C_\omega,
$$ 
for $|t|$ sufficiently small.
More precisely, there exist $t_{1}>0$, an open set $D^\prime\subset \R^{N-1}$, and a smooth function $w\colon\overline{D^\prime}\times (-t_1,t_1)\to \R$ such that for all $t\in\left]-t_{1},t_{1}\right[$ 
\begin{eqnarray*}
\Gamma_t\cap \overline{B_r(x_0)}&=&\{x=(x^\prime, x_N)\in \RN\,;~x_N=w(x^\prime,t), \, x^\prime\in \overline{D^\prime}\},\\
\Omega_{t}&=&\left(\Omega\setminus B_{r}(x_{0})\right)\cup\{x=(x^\prime, x_N)\in B_{r}(x_{0})\,;~ x_N<w(x^\prime,t), \, x^\prime\in \overline{D^\prime}\}.
\end{eqnarray*}
Notice that $\Omega_0=\Omega$ and $\Gamma_0\cap \overline{B_r(x_0)}=\Gamma\cap \overline{B_r(x_0)}$. Now we perturb $\Omega_{t}$ in a neighborhood of $x_{0}$ in order to preserve the volume. To this aim, we fix a mapping $g_{0}\in C^{\infty}_{c}(D')$ such that $g\ge 0$ and 
\begin{equation}\label{defg}
\int_{D'}g_{0}(x')\,dx'=1.
\end{equation}
We define
\begin{equation}\label{defh}
h(t):=\begin{cases}\displaystyle \frac{|\Omega|-|\Omega_t|}{t^2} & \hbox{if $t\in \left]-t_1, t_1\right[\setminus\{0\}$,}\\[6pt] \displaystyle  -\frac{1}{2}\frac{d^2}{dt^2}(|\Omega_t|)\big\vert_{t=0} & \hbox{if $t=0$},\end{cases}
\end{equation}
\begin{equation}\label{defgh}
g(x',t):=g_{0}(x')h(t)\quad\quad\forall (x',t)\in D'\times\left]-t_{1},t_{1}\right[,
\end{equation}
and
$$
\widetilde\Omega_{t}:=\left(\Omega\setminus B_{r}(x_{0})\right)\cup\left\{x=(x^\prime, x_N)\in B_{r}(x_{0})\,;~ x_N<w(x^\prime,t)+t^{2}g(x',t), \, x^\prime\in \overline{D^\prime}\,\right\}\quad\forall t\in\left]-t_{1},t_{1}\right[.
$$
Then, by \eqref{defg}--\eqref{defgh}, for every $t\in\left]-t_{1},t_{1}\right[$ we have
\begin{equation}\label{volconst}
\big|\widetilde\Omega_{t}\big|=\left|\Omega_{t}\right|+\left|\left\{x=(x',x_{N})\in\R^{N}\,;~x'\in D',~0\le x_{n}\le t^{2}g(x',t)\right\}\right|=\left|\Omega_{t}\right|+t^{2}\int_{D'}g(x',t)\,dx'=\left|\Omega\right|.
\end{equation}
Since $\Omega$, by assumption, is a minimizer of the functional $\E(\cdot;\C_\omega)$ constrained to the family of the sets of volume $c$, by the Lagrange multiplier Theorem, there exists a constant $\lambda\in\R$ such that
$$
\frac{d}{dt}\left.\left( \mathcal{E}(\widetilde\Omega_t;\C_\omega)\right)\right|_{t=0}=\lambda  \frac{d}{dt}\left.\left(|\widetilde\Omega_t|\right)\right|_{t=0}. 
$$ 
Since $\widetilde\Omega_{t}$ differs from $\Omega_{t}$ by a term of order $t^{2}$, we have that
$$
\frac{d}{dt}\left.\left( \mathcal{E}(\widetilde\Omega_t;\C_\omega)\right)\right|_{t=0}=\frac{d}{dt}\left.\left( \mathcal{E}(\Omega_t;\C_\omega)\right)\right|_{t=0}\quad\text{and}\quad \frac{d}{dt}\left.\left(\big|\widetilde\Omega_t\big|\right)\right|_{t=0}=\frac{d}{dt}\left.\left(|\Omega_t|\right)\right|_{t=0}.
$$ 
Then, by \eqref{eq:diffTorsionalOmegat}, \eqref{eq:diffVolOmegat}, and \eqref{volconst}, we obtain that 
$$
\int_{\Gamma\cap \overline{B_r(x_0)}}  |\nabla u_\Omega|^2 \langle V,\nu\rangle \ d\sigma=0\quad\text{and}\quad \int_{\Gamma\cap \overline{B_r(x_0)}}  \langle V,\nu\rangle \ d\sigma=0,
$$
from which we infer that for all $C\in\R$
\begin{equation}\label{eq:diffTorsionalOmegat3}
\int_{\Gamma\cap \overline{B_r(x_0)}}  \left(|\nabla u_\Omega|^2 -C\right) \langle V,\nu\rangle \ d\sigma=0.
\end{equation}
Plugging the expression of $V$, given by \eqref{eq:defV}, into \eqref{eq:diffTorsionalOmegat3} and choosing $C=C_\varphi$, where $C_\varphi$ is the constant given by \eqref{eq:defCvarphi}, we obtain
$$ \int_{\Gamma\cap \overline{B_r(x_0)}} \varphi\left(|\nabla u_\Omega|^2-C_\varphi\right)^2 \, d\sigma=0.$$
In particular, since by construction it holds $\varphi\geq 0$ in $B_r(x_0)$, and $\varphi\equiv 1$ on $\overline{S}$, we get that
$$ \int_{\overline{S}} \left(|\nabla u_\Omega|^2-C_\varphi\right)^2 \, d\sigma=0$$
which implies that $|\nabla u_\Omega|^2$ is constant on $\overline{S}$, a contradiction. The proof of \textit{Claim I} is complete.
\medskip

We now conclude the proof of the first part of the Proposition. Let $x_1\in \mathrm{Reg}(\Gamma)$, $x_2\in \mathrm{Reg}(\Gamma)$, with $x_1\neq x_2$, and let $B_{r_1}(x_1)$, $B_{r_2}(x_2)$ be two disjoint balls, where, for $i=1,2$, the radius $r_i>0$ is sufficiently small so that the conditions stated at the beginning of the proof are satisfied, and $\Gamma\cap \overline{B_{r_i}(x_i)}$ is the graph of a smooth function. Let $V:\R^N\to\R^N$ a smooth vector field with compact support in $B_{r_1}(x_1) \cup B_{r_2}(x_2)$ and such that
\beq\label{eqGammaBrV}
\int_{\Gamma\cap \overline{B_{r_1}(x_1)}}\langle V,\nu\rangle \, d\sigma=-\int_{\Gamma\cap \overline{B_{r_2}(x_2)}}\langle V,\nu\rangle \, d\sigma\neq 0.
\eeq
A possible choice is
$$
V=\begin{cases}\varphi_1\overline{\nu_1} & \hbox{in $B_{r_1}(x_1)$,}\\
\varphi_2\overline{\nu_2} & \hbox{in $B_{r_2}(x_2)$,}\\
0 & \hbox{in $B_{r_1}^\complement(x_1)\cap B_{r_2}^\complement(x_2) $,}\end{cases}
$$
where, for $i=1,2$, $\overline{\nu_i}$ is the extension to a smooth vector field defined in $\overline{B_{r_i}(x_i)}$ of the outer normal $\nu\big\vert_{\Gamma\cap {B_{2r_i}(x_i)}}$, and $\varphi_i\in C^\infty_c(B_{r_i}(x_i))$ satisfy
$$ \int_{\Gamma\cap \overline{B_{r_1}(x_1)}} \varphi_1 \, d\sigma=-\int_{\Gamma\cap \overline{B_{r_2}(x_2)}} \varphi_2 \, d\sigma\neq 0.$$
Clearly, from \eqref{eqGammaBrV} we have
$$\int_{(\Gamma\cap \overline{B_{r_1}(x_1)})\cup (\Gamma\cap \overline{B_{r_2}(x_2)})}  \langle V,\nu\rangle \ d\sigma=0.$$
Then, considering the deformation induced by the flow associated to $V$, and arguing as in the proof of \textit{Claim I}, with slightly adjustments, we can construct a volume-preserving deformation starting from $V$.
Hence, as in the previous case, by minimality of $\Omega$ we have
\beq\label{eq:vcminimality}
- \frac{1}{2}\int_{(\Gamma\cap \overline{B_{r_1}(x_1)})\cup (\Gamma\cap \overline{B_{r_2}(x_2)})}  |\nabla u_\Omega|^2 \langle V,\nu\rangle \ d\sigma=0.
\eeq
Now, by \textit{Claim I}, we have that $|\nabla u_\Omega|^2\equiv C_i$ is constant on $\Gamma\cap\overline{B_{r_i}(x_i)}$, for $i=1,2$, and thus from \eqref{eqGammaBrV} and \eqref{eq:vcminimality} we infer that
$$ (C_1-C_2) \int_{\Gamma\cap \overline{B_{r_1}(x_1)}}  \langle V,\nu\rangle \ d\sigma = 0.$$
Hence, as the second factor is different from zero by \eqref{eqGammaBrV}, the only possibility is $C_1=C_2$ and from the arbitrariness of $x_1$, $x_2$ we conclude that $|\nabla u_\Omega|^2$ is constant on $\mathrm{Reg}(\Gamma)$. The proof of the first part of the proposition is complete.
\medskip

Let us prove now the bound \eqref{eq:boundlemmasect5}. Let $x_0\in \mathrm{Reg}(\Gamma)$ and let $B_r(x_0)$ a small ball satisfying the conditions stated at the beginning of the proof. Choose a smooth vector field $V$ with compact support in $B_r(x_0)$ and such that $\langle V,\nu\rangle\leq 0$ on $\Gamma\cap B_r(x_0)$ and $\langle V,\nu\rangle<0$ in $B_\sigma(y_0)$ for some $B_\sigma(y_0)\subset B_r(x_0)$, and consider the induced flow. In particular, from \eqref{eq:diffVolOmegat}, and since by construction $\Omega_0=\Omega$, then we deduce that there exists $t_1\in ]0,t_0[$ sufficiently small such that 
\beq\label{eq:volOmegat}
|\Omega_t|\leq |\Omega|,
\eeq
for all $0\leq t<t_1$.
In view of Lemma \ref{lem1:sect5}, $\Omega$ is also a minimizer for $\frac{\mathcal{E}(\widetilde\Omega;\C_\omega)}{|\widetilde\Omega|}$ in the class of subsets $\widetilde\Omega$ of $\A_{\omega, c}$ satisfying $0<|\widetilde\Omega|\leq |\Omega|$. Hence, from \eqref{eq:volOmegat}, exploiting \eqref{eq:diffTorsionalOmegat}, \eqref{eq:diffVolOmegat}, and taking into account that  $ |\nabla u_\Omega|^2\equiv C_0$ on $\mathrm{Reg}(\Gamma)$ for some constant $C_0>0$, $|\Omega|=c$, and $\mathcal{E}(\Omega;\C_\omega)= - |\mathcal{O}_c(\C_\omega)|$, we infer that
$$\frac{d}{dt}\left.\left(\frac{\mathcal{E}(\Omega_t;\C_\omega)}{|\Omega_t|} \right)\right|_{t=0^+}=- \frac{1}{2c}\int_{\Gamma\cap \overline{B_r(x_0)}} \left(C_0-\frac{2|\mathcal{O}_c(\C_\omega)|}{c} \right) \langle V,\nu\rangle \, d\sigma\geq0.$$
Hence, as $\int_{\Gamma\cap \overline{B_r(x_0)}} \langle V,\nu\rangle\, d\sigma< 0$ we deduce that the only possibility is \eqref{eq:boundlemmasect5} 
%% $$C_0\geq\frac{2}{c} |\mathcal{O}_c(\C_\omega)|,$$
and the proof is complete.
\end{proof}

\begin{remark}
Since $u_{\Omega^*}$ is a viscosity solution (in the sense of \cite[Definition 2.2]{DS}) to 
$$
\begin{cases}-\Delta u=1&\text{ in $\Omega^*$}\\ |\nabla u|=\sqrt{2\Lambda}&\text{ on $\Gamma_{\Omega^*}$}\end{cases}
$$
where $\Lambda$ is the Lagrange multiplier obtained in \textit{Step 2} of the proof of Theorem \ref{teoreg}, we have that $\Lambda=\frac{C_0}2$, where $C_0$ is the constant of Proposition \ref{lem2:sect5}, and thus by \eqref{eq:boundlemmasect5} we get
$$\Lambda \geq \frac{|\O_c(\C_\omega)|}{c}.$$
\end{remark}

The following result yields, under suitable assumptions, an estimate on the measure of the relative boundary of the minimizing domain for $\O_c(\C_\omega)$. In particular, this estimate provides an upper bound which is independent of $c$. This seems particularly interesting for large values of $c$. %For convenience, let us introduce the following definition.%, adapting a similar one stated in \cite[Sect. 9.3]{MA}). 

\begin{proposition}\label{P:bdr-estimate-1}
Let $\omega\subset\R^{N-1}$ be a bounded domain of class $C^{2,\alpha}$ and let $c>0$. Let  $\Omega^*$ be a minimizer of $\E(\cdot;\C_\omega)$ in $\A_{\omega,c}$. % and set
%$$
%\begin{array}{lll}
%&\Gamma_{\Omega^*}:=\partial\Omega^*\cap\C_\omega\,,~~&\Gamma_{\Omega^*}^{\mathrm{reg}}:=\Gamma_{\Omega^*}\cap(\partial\Omega^{*})^{\mathrm{reg}}\vspace{6pt}\\
%&\Gamma_{1,\Omega^*}:=\partial\Omega^*\cap\partial\C_\omega\,,~~&\Gamma_{1,\Omega^*}^{\mathrm{reg}}:=\mathrm{int}_{\partial\Omega^{*}}(\Gamma_{1,\Omega^*}\cap(\partial\Omega^{*})^{\mathrm{reg}}),
%\end{array}
%$$
%where $(\partial\Omega^{*})^{\mathrm{reg}}$ denotes the regular part of $\partial\Omega^{*}$ and $\mathrm{int}_{\partial\Omega^{*}}(E)$ is the interior of $E\subset\partial\Omega^{*}$ in the topology of $\partial\Omega^{*}$. 
Assume that the energy function $u_{\Omega^*}$ belongs to $W^{1,\infty}(\Omega^*)$ and that 
\begin{equation}
\label{assumpt-reg}
\mathcal{H}^{N-1}\left(\Gamma_{1,\Omega^*}\setminus\mathrm{int}_{\partial\mathcal{C}_\omega}(\Gamma_{1,\Omega^*})\right)=0
%\partial\Omega^{*}\setminus(\Gamma_{\Omega^*}^{\mathrm{reg}}\cup\Gamma_{1,\Omega^*}^{\mathrm{reg}}))=0
\end{equation}
where $\Gamma_{1,\Omega^*}:=\partial\Omega^*\cap\partial\C_\omega$ and $\mathrm{int}_{\partial\mathcal{C}_\omega}(\Gamma_{1,\Omega^*})$ is the interior of $\Gamma_{1,\Omega^*}$ in the topology of $\partial\mathcal{C}_\omega$. 
%$\Omega^*$ has almost $C^{1}$-boundary and there exists a relative (to $\partial\Omega^*$) open subset $\mathcal{U}_{1,\Omega^*}$, possibly empty, such that  $\mathcal{U}_{1,\Omega^*}\subset \Gamma_{1, \Omega^*}\cap(\partial\Omega^*)^{\mathrm{reg}}$ and $\H^{N-1}(\Gamma_{1, \Omega^*}\setminus\mathcal{U}_{1,\Omega^*})=0$. Moreover, assume 
Then, it holds that 
\beq\label{eqtesi:stimaGammasect5}
\H^{N-1}(\Gamma_{\Omega^*})\leq   2\sqrt{3}\,\H^{N-1}(\omega).
\eeq
\end{proposition}
\begin{proof}
Let us fix $c>0$ and, for $h>0$, let $\Omega_{\omega,h}\subset \C_\omega$ be the bounded cylinder of height $h$, that is $\Omega_{\omega,h}:=\omega\times\left]-\frac{h}{2},\frac{h}{2}\right[$. It is elementary to check that the energy function of $\Omega_{\omega,h}$ is given by
\beq\label{eq:torsionfunctbcylinder}
u_{\Omega_{\omega,h}}(x^\prime,x_N)=\begin{cases}\frac12\big(\frac{h^2}{4}-x_N^2\big)&\hbox{if $|x_N|<\frac{h}{2}$}\\0&\hbox{if $|x_N|\geq \frac{h}{2}$}\end{cases}.
\eeq
Then, by definition and a straightforward computation, we get that
\beq\label{eq1:lem3sect5}
\E(\Omega_{\omega,h}; \C_{\omega})=-\frac{1}{2} \int_{\Omega_{\omega,h}} u_{\Omega_{\omega,h}}\, dx= -\frac{1}{4}\H^{N-1}(\omega) \int_{-\frac{h}{2}}^{\frac{h}{2}}\left(\frac{h^2}{4}-x_N^2\right) \, dx_N=-\frac{1}{24} \H^{N-1}(\omega) h^3.
\eeq
When $h=\frac{c}{\H^{N-1}(\omega)}$, we have that $\Omega_{\omega,h}\in\A_{\omega,c}$ and then
$$
\O_c(\C_\omega)\leq \E( \Omega_{\omega,h}; \C_\omega)=-\frac{1}{24}\frac{c^3}{(\H^{N-1}(\omega))^2}<0.
$$
Hence, by \eqref{eq1:lem3sect5}, 
\beq\label{eq2:lem3sect5}
|\O_c(\C_\omega)|\geq \frac{1}{24}\frac{c^3}{(\H^{N-1}(\omega))^2}.
\eeq
Now let us find an upper bound for $|\O_c(\C_\omega)|$. To this aim, let $\Omega^* \subset\C_\omega$ a minimizer of $\E(\cdot;\C_\omega)$ in $\A_{\omega,c}$ satisfying \eqref{assumpt-reg} and let $u_{\Omega^*}$ be its energy function which, by assumption, belongs to $W^{1,\infty}(\Omega^{*})$. As in the proof of Proposition \ref{lem2:sect5} we simplify the notation and write $\Omega$, $\Gamma$, $\Gamma_1$, %$\Gamma^{\mathrm{reg}}$, $\Gamma_1^{\mathrm{reg}}$, 
$u_{\Omega}$ instead of $\Omega^*$, $\Gamma_{\Omega^*}$, $\Gamma_{1,\Omega^*}$, %$\Gamma_{\Omega^*}^{\mathrm{reg}}$, $\Gamma_{1,\Omega^*}^{\mathrm{reg}}$, 
$u_{\Omega^{*}}$, respectively. By Proposition \ref{prop2:sect5} and Corollary \ref{cor1:sect5} we know that $\Omega$ is an open set and $|\Omega|=c$. Moreover, as shown in the proof of Theorem \ref{teoreg}, the regular part of the relative boundary $\Gamma$, denoted $\mathrm{Reg}(\Gamma)$ and defined in Definition \ref{def:regsingparts}, is locally a smooth manifold, is a relative open subset of $\partial\Omega$, and 
\begin{equation}\label{singular-Gamma-negligible}
\mathcal{H}^{N-1}(\Gamma\setminus\mathrm{Reg}(\Gamma))=0\,.
\end{equation} 
In addition, $\mathrm{int}_{\partial\mathcal{C}_\omega}(\Gamma_{1})$ is a $C^{2,\alpha}$ manifold (because $\partial\omega$ is assumed to be so) and it is a relative open subset of $\partial\Omega$, too. We set
\begin{equation}\label{regular-partial-omega}
\mathrm{Reg}(\partial\Omega):=\mathrm{Reg}(\Gamma)\cup\mathrm{int}_{\partial\mathcal{C}_\omega}(\Gamma_{1})
\end{equation}
and we point out that, by the assumption \eqref{assumpt-reg} and by \eqref{singular-Gamma-negligible},
\begin{equation}\label{sing-partial-omega}
\mathcal{H}^{N-1}(\partial\Omega\setminus\mathrm{Reg}(\partial\Omega))=0\,.
\end{equation}
By standard elliptic regularity theory, the energy function $u_\Omega$ belongs to $C^2(\Omega\cup\mathrm{Reg}(\partial\Omega))$ and solves
\begin{equation}\label{neumann}
\begin{cases}-\Delta u = 1&\text{~~in~~$\Omega$}\vspace{3pt}\\
\dfrac{\partial u}{\partial \nu}= 0&\text{~~on~~$\mathrm{int}_{\partial\mathcal{C}_\omega}(\Gamma_{1})$}\vspace{3pt}\\
\dfrac{\partial u}{\partial \nu} = -\sqrt{C_0}&\text{~~on~~$\mathrm{Reg}(\Gamma)$}
\end{cases}
\end{equation}
for some positive constant $C_0$ satisfying \eqref{eq:boundlemmasect5}. 
%%\beq\label{eq0:lem3sect5}
%%\Lambda_*\geq \frac{|\O_c(\C_\omega)|}{c}\,.
%%\eeq
Since, by assumption, $u_\Omega\in W^{1,\infty}(\Omega)$, taking into account \eqref{sing-partial-omega}, we obtain that: 
\beq\label{eq:teodiv}
\int_\Omega \Delta u_\Omega\, dx=\int_{\mathrm{Reg}(\partial\Omega)}\frac{\partial u_\Omega}{\partial \nu}\, d\sigma. 
\eeq 
We postpone the proof of \eqref{eq:teodiv} and we complete the proof of \eqref{eqtesi:stimaGammasect5}. Taking into account \eqref{neumann} and using \eqref{regular-partial-omega} and \eqref{eq:teodiv}, we get that
$$
c=|\Omega|=\int_\Omega 1\, dx=-\int_\Omega \Delta u_\Omega\, dx=%-\int_{\mathrm{Reg}(\Gamma)}\frac{\partial u_\Omega}{\partial \nu}\, d\sigma=
\sqrt{C_0}\,\H^{N-1}(\mathrm{Reg}(\Gamma)).
$$
Then, by \eqref{singular-Gamma-negligible}, we conclude that
\beq\label{eq4:lem3sect5}
\H^{N-1}(\Gamma)=\frac{c}{\sqrt{C_0}}\,.
\eeq
Finally, from \eqref{eq:boundlemmasect5}, \eqref{eq2:lem3sect5}, and \eqref{eq4:lem3sect5}, we obtain \eqref{eqtesi:stimaGammasect5}.
\medskip

\textbf{Proof of \eqref{eq:teodiv}:} We suitably adapt the proof of \cite[Theorem 9.6]{MA}. 
We have that $\partial\Omega=\mathrm{Reg}(\partial\Omega)\cup M_{0}$ where $\mathrm{Reg}(\partial\Omega)$ is a manifold of class $C^{2,\alpha}$ ($C^{1}$ would be enough at this step), and $M_{0}:=\partial\Omega\setminus\mathrm{Reg}(\partial\Omega)$ is a closed subset of $\partial\Omega$, with $\H^{N-1}(M_0)=0$. Then, for any fixed $\eps\in(0,1)$, there exists a countable cover $\{B_{\eps_k}(y_k)\}_{k\in\N}$ of $M_0$, made by balls centered at $y_k\in M_0$ and such that
\beq\label{eq:boundepsk}
%% \sum_{k\in\N}2^{N-1} \eps_k^{N-1}< \eps.
\sum_{k\in\N} \eps_k^{N-1}\le \eps.
\eeq
Let us also consider the open cover of $\mathrm{Reg}(\partial\Omega)$ given by $\{B_{r}(x)\}_{x\in\mathrm{Reg}(\partial\Omega)}$, where $r=r(x)>0$. 
% is the number provided by Definition \ref{def:regular-part}. 
Since $\Omega$ is bounded (see Proposition \ref{prop1:sect5}) then $\partial\Omega$ and $\overline\Omega$ are compact sets. In particular, we find a finite subcover of $\partial\Omega$, namely
$$\partial\Omega\subset \bigcup_{k\in I}B_{\eps_k}(y_k) \cup  \bigcup_{h\in J}B_{r_h}(x_h),$$
where $I,J \subset \N$ are finite sets, $x_h\in\mathrm{Reg}(\partial\Omega)$, $r_h=r(x_h)$ for all $h\in J$. Moreover, by compactness of $\overline\Omega$, we find a finite set of balls $\{B_{s_l}(z_l)\}_{l\in L}$ such that $\overline{B_{s_l}(z_l)}\subset\Omega$, for all $l\in L$, and
$$\overline\Omega \subset \bigcup_{k\in I}B_{\eps_k}(y_k) \cup  \bigcup_{h\in J}B_{r_h}(x_h) \cup  \bigcup_{l\in L}B_{s_l}(z_l),$$
for some finite set $L\subset\N$.
Then, there exists a partition of unity subordinated to this open cover. In particular there exist $\eta_k \in C^1_c(B_{\eps_k}(y_k))$, $\zeta_h \in C^1_c(B_{r_h}(x_h))$, $\theta_l \in C^1_c(B_{s_l}(z_l))$ such that 
\begin{gather*}
0\leq \eta_k\leq 1~~\forall k\in I,~~ 0\leq \zeta_h\leq 1~~\forall h\in J,~~ 0\leq \theta_l\leq 1~~\forall l\in L,\vspace{6pt}\\
\sum_{k\in I}\eta_k +\sum_{h\in J} \zeta_h + \sum_{l\in L}\theta_l\equiv 1\text{~~on~~}\overline\Omega.
\end{gather*} 
Moreover, we can choose $\eta_k$ so that
\beq\label{eq:condnablaetak}
|\nabla \eta_k| \leq \frac{C}{\eps_k}\text{~~in~~}B_{\eps_k}(y_k),~~\forall k\in I,
\eeq
for some constant $C>0$ independent of $k$. We also notice that, since $\overline{B_{s_l}(z_l)}\subset\Omega$, for all $l\in L$, then
\beq\label{eq:reletakzetah}
\sum_{k\in I}\eta_k +\sum_{h\in J} \zeta_h \equiv 1\ \ \ \hbox{on $\partial\Omega$}.
\eeq
Now, by construction we have
\beq\label{eq:termsintegraldivteo}
\int_\Omega \mathrm{div}(\nabla u_\Omega)\, dx=  \sum_{k\in I} \int_\Omega \mathrm{div}(\eta_k\nabla u_\Omega )\, dx +\sum_{h\in J} \int_\Omega \mathrm{div}( \zeta_h\nabla u_\Omega)\, dx+\sum_{l\in L } \int_\Omega \mathrm{div}(\theta_l \nabla u_\Omega)\, dx.
\eeq
Since  $\overline{B_{s_l}(z_l)}\subset\Omega$, we have $u_\Omega\in C^2(\overline{B_{s_l}(z_l)})$, for all $l\in L$, and, as $\theta_l \in C^1_c(B_{s_l}(z_l))$, applying the standard divergence theorem, we infer that
\beq\label{eq:estimateterm3}
\sum_{l\in L } \int_\Omega \mathrm{div}(\theta_l \nabla u_\Omega)\, dx=0.
\eeq
Let us analyze the remaining terms in the right-hand side of \eqref{eq:termsintegraldivteo}. 
%To this end we first observe that, in view of Theorem \ref{teoreg} and since by assumption $\mathcal{U}_1\subset \Gamma_1\cap \partial\Omega^{\mathrm{reg}}$ it holds that $\Gamma^{\mathrm{reg}}\cup\mathcal{U}_1\subset  \partial\Omega^{\mathrm{reg}}$. Actually, as $\partial\Omega^{\mathrm{reg}}$ is a $C^1$-hypersurface, since $\Gamma^{\mathrm{reg}}$, $\U_1$ are two relative open subsets of $\partial\Omega$ and $\H^{N-1}(\partial\Omega\setminus \partial\Omega^{\mathrm{reg}})=0$, $\H^{N-1}(\Gamma\setminus\Gamma^{\mathrm{reg}})=0$, $\H^{N-1}(\Gamma_1\setminus\mathcal{U}_1)=0$, then, the only possibility is that $\partial\Omega^{\mathrm{reg}}=\Gamma^{\mathrm{reg}}\cup \U_1$. Hence, r
Recalling that $u_\Omega\in C^2(\Omega\cup\mathrm{Reg}(\partial\Omega))$, $\zeta_h\in C_c^1(B_{r_h}(x_h))$, and the definition of $B_{r_h}(x_h)$, it turns out that $\Omega\cap B_{r_h}(x_h)$ is a Lipschitz domain, for all $h\in J$, and applying the divergence theorem, we get that
\beq\label{eq:termsintegraldivteo2}
\sum_{h\in J} \int_\Omega \mathrm{div}( \zeta_h\nabla u_\Omega)\, dx= \sum_{h\in J} \int_{\mathrm{Reg}(\partial\Omega)\cap B_{r_h}(x_h)} \zeta_h\,\frac{\partial u_\Omega}{\partial \nu}\, d\sigma.
\eeq
Next, for the first term in the right-hand side of \eqref{eq:termsintegraldivteo}, from \eqref{neumann}, \eqref{eq:boundepsk},\eqref{eq:condnablaetak}, since $\eta_k \in C^1_c(B_{\eps_k}(y_k))$ and by assumption $u_\Omega \in W^{1,\infty}(\Omega)$, we deduce that
\beq\label{eq:estimateterm1}
\begin{array}{lll}
\displaystyle \left|\sum_{k\in I} \int_\Omega \mathrm{div}(\eta_k\nabla u_\Omega )\, dx  \right|&\leq&\displaystyle \sum_{k\in I}  \int_{B_{\eps_k}(y_k)\cap\Omega} \left( |\Delta u_\Omega| + |\nabla u_\Omega| |\nabla \eta_k| \right)\, dx\\[3pt]
&\leq&\displaystyle  \sum_{k\in I} \left(1+\frac{C}{\eps_k} |\nabla u_\Omega|_{\infty}\right)|B_{\eps_k}(y_k)|\\[3pt]
&\leq& \displaystyle \sigma_N  \left(1+{C} |\nabla u_\Omega|_{\infty}\right) \sum_{k\in I} \eps_k^{N-1}\\[3pt]
&\leq& \displaystyle \sigma_N  \left(1+{C} |\nabla u_\Omega|_{\infty}\right) \eps%%\\[3pt]
%% &\leq& \displaystyle \widetilde C(N,|\nabla u_\Omega|_\infty) \eps,
\end{array}
\eeq
where $\sigma_N=|B_1(0)|$. 
%, $\widetilde C(N,|\nabla u_\Omega|_\infty):=\frac{\sigma_N}{2^{N-1}} \left(1+{C} |\nabla u_\Omega|_{\infty}\right)$.\\
Summing up, from \eqref{eq:termsintegraldivteo}--%, \eqref{eq:estimateterm3}, \eqref{eq:termsintegraldivteo2} and 
\eqref{eq:estimateterm1}, we have
\beq\label{eq:finalestimate1}
\left|\int_{\mathrm{Reg}(\partial\Omega)}\frac{\partial u_\Omega}{\partial \nu}\, d\sigma-\int_\Omega \mathrm{div}(\nabla u_\Omega )\, dx\right|\leq \left|\int_{\mathrm{Reg}(\partial\Omega)}  \bigg(1-\sum_{h\in J} \zeta_h\bigg)
\frac{\partial u_\Omega}{\partial \nu}\, d\sigma \right| +  \sigma_N  \left(1+{C} |\nabla u_\Omega|_{\infty}\right) \eps.
%\widetilde C(N,|\nabla u_\Omega|_\infty)\eps.
\eeq
Setting $\mathcal{V}_\eps:=\bigcup_{k\in I}B_{\eps_k}(y_k)$, by \eqref{eq:reletakzetah} we have $ \sum_{h\in J} \zeta_h\equiv 1$ on $\mathrm{Reg}(\partial\Omega)\setminus \mathcal{V}_\eps$, and exploiting that $u_\Omega \in C^1(\Omega\cup\mathrm{Reg}(\partial\Omega))\cap W^{1,\infty}(\Omega)$ we obtain that
\beq\label{eq:finalestimate2}
 \left|\int_{\mathrm{Reg}(\partial\Omega)}  \bigg(1-\sum_{h\in J} \zeta_h\bigg) \frac{\partial u_\Omega}{\partial \nu}\, d\sigma \right|\leq |\nabla u_\Omega|_\infty \mathcal{H}^{N-1}(\mathrm{Reg}(\partial\Omega)\cap \V_\eps).
\eeq
Finally, as $\V_\eps$ is bounded, it holds that $\mathcal{H}^{N-1}(\mathrm{Reg}(\partial\Omega)\cap \V_\eps)<+\infty$ for all $\eps\in (0,1)$. Therefore, since $M_0=\bigcap_{\eps\in(0,1)} \V_\eps$, it follows that $\lim_{\eps\to 0^+} \mathcal{H}^{N-1}(\mathrm{Reg}(\partial\Omega)\cap \V_\eps)=\mathcal{H}^{N-1}(\mathrm{Reg}(\partial\Omega)\cap M_0)=0$. Hence, taking the limit as $\eps\to 0^+$ in \eqref{eq:finalestimate1}, \eqref{eq:finalestimate2}, 
%recalling that $\partial\Omega^{\mathrm{reg}}=\Gamma^{\mathrm{reg}}\cup \U_1$, and since $\Gamma^{\mathrm{reg}}\cap\U_1=\varnothing$, 
we obtain \eqref{eq:teodiv}. The proof is complete.
\end{proof}

The next result states another upper bound on the measure of the relative boundary of the minimizing domain for $\E(\cdot; \C_\omega)$ which is meaningful in particular for small values of $c$.

\begin{proposition}\label{P:bdr-estimate-2}
Under the same assumptions of Proposition \ref{P:bdr-estimate-1} it holds that
\beq\label{eqtesi:stimaGammasect5-B}
\H^{N-1}(\Gamma_{\Omega^*})\le c^{1-\frac1N}\left[\sqrt{N(N+2)}\left(\frac{\sigma_{N}}2\right)^{\frac1N}+o(1)\right]\quad\text{as }c\to 0^+,
\eeq
where $\sigma_{N}=|B_{1}(0)|$.
\end{proposition}
\begin{proof}
We follow the same strategy and notation used to show Proposition \ref{P:bdr-estimate-1} apart from the choice of the comparison domain. Indeed, for small $c$, as a comparison domain we take the intersection of a ball centered at some fixed point $x_{0}\in\partial\C_{\omega}$ with the container itself.  
Fixing such a point $x_{0}\in\partial\C_{\omega}$, there exists $\bar{r}>0$, an isometry $T\colon\R^{N}\to\R^{N}$ and a Lipschitz function $f\colon\{x'\in\R^{N-1}\,;~|x'|<\bar{r}\}\to\R$ such that $f(0)=0$, $f$ is differentiable at $0$ with $\nabla f(0)=0$ and 
$$
B_{\bar r}(x_{0})\cap\C_{\omega}=T(\{x=(x',x_{N})\in\R^{N}\,;~|x|<\bar r,~x_{N}>f(x')\}). 
$$
In the first step, for every fixed $c>0$ small enough, we select the competitor domains. 
\medskip

\noindent
\textbf{Step 1.} There exist $c_{0}>0$, $r_{0}\in\left]0,\bar{r}\right]$ and a continuous, increasing function $r\colon\left]0,c_{0}\right[\to\left]0,r_{0}\right[$ (depending on $x_{0}$) such that
\begin{equation}\label{st1}\begin{split}
&|B_{r(c)}(x_{0})\cap\C_{\omega}|=c\quad\forall c\in \left]0,c_{0}\right[\\
&r(c)\sim\left(\frac{2c}{\sigma_{N}}\right)^{\frac1N}\quad\text{as $c\to 0^+$.}
\end{split}\end{equation}
\medskip

\noindent
%Moreover there exist an isometry $T\colon\R^{N}\to\R^{N}$ and a smooth function $f\colon B'_{r_{0}}(0)\subset\R^{N-1}\to\R$ such that $f(0)=0$, $\nabla f(0)=0$ and 
%\begin{equation}\label{st2}
%B_{r(c)}(x_{0})\cap\C_{\omega}=T(\{x=(x',x_{N})\in\R^{N}~:~|x|<r(c),~x_{N}>f(x')\})\quad\forall c\in (0,c_{0}). 
%\end{equation}
%%
%%\emph{Proof of Step 1.} 
For every $r\in\left]0,\bar{r}\right[$ set ${\Omega}_{r}:=\{x=(x',x_{N})\in\R^{N}\,;~|x|<r,~x_{N}>f(x')\}$ where $\bar{r}>0$ and $f$ are given as before. Then $B_{r}(x_{0})\cap\C_{\omega}=T(\Omega_{r})$ and, since $T$ is an isometry,
$$
|B_{r}(x_{0})\cap\C_{\omega}|=|\Omega_{r}|=r^{N}|\widetilde{\Omega}_{r}|
$$
where
\begin{equation}\label{omegatilder}
\widetilde{\Omega}_{r}:=\left\{x=(x',x_{N})\in\R^{N}\,;~|x|<1,~x_{N}>\frac{f(rx')}{r}\right\}.
\end{equation}
Since $f(0)=0$ and $\nabla f(0)=0$, the condition $x_{N}>\frac{f(rx')}{r}$ in the limit as $r\to 0^+$ becomes $x_{N}>0$. Hence
\begin{equation}\label{alphaN}
|\widetilde{\Omega}_{r}|\to\frac{\sigma_{N}}{2}\text{~~as~~}r\to 0^{+}.
\end{equation}
Let us introduce the function $g\colon[0,\bar{r})\to\R$ defined by
$$
g(r):=\begin{cases}0&\text{as $r=0$}\\ r\,|\widetilde{\Omega}_{r}|^{\frac1N}&\text{as $r\in\left]0,\bar{r}\right[$.}\end{cases}
$$
This function is continuous and, by \eqref{alphaN}, $g'(0)>0$. Hence, in a right neighborhood of $0$ it admits an inverse function $g^{-1}\colon\left[0,\varepsilon_{0}\right[\to \left[0,r_{0}\right[$ which turns out to be continuous and increasing. Finally, one can plainly check that the function $r(c):=g^{-1}\big(c^{\frac1N}\big)$ for $c\in\left]0,c_{0}\right[$ with $c_{0}:=\varepsilon_{0}^{N}$, satisfies \eqref{st1} and 
\begin{equation}\label{st2}
B_{r(c)}(x_{0})\cap\C_{\omega}=T(\{x=(x',x_{N})\in\R^{N}\,;~|x|<r(c),~x_{N}>f(x')\}). 
\end{equation}
%\medskip
%
%\noindent
\textbf{Step 2.} For every $r\in(0,r_{0})$ the function
$$
u(x):=\begin{cases}\frac1{2N}(r^{2}-|x-x_{0}|^{2})&\text{as $x\in B_{r}(x_{0})\cap\C_{\omega}$}\\ 0&\text{as $x\in \C_{\omega}\setminus B_{r}(x_{0})$}\end{cases}
$$
belongs to $H^{1}_{0}(B_{r}(x_{0})\cap\C_{\omega};\C_{\omega})$ and
\begin{equation}\label{st3}
J_{B_{r}(x_{0})\cap\C_{\omega}}(u)=-r^{N+2}\left[\frac{\sigma_{N}}{4N(N+2)}+o(1)\right]\quad\text{as }r\to 0^+.
\end{equation}
%\emph{Proof of Step 2.} 
The fact that $u\in H^{1}_{0}(B_{r}(x_{0})\cap\C_{\omega};\C_{\omega})$ is immediate. Moreover, since $B_{r}(x_{0})\cap\C_{\omega}=T(\Omega_{r})$, we have that
$$
J_{B_{r}(x_{0})\cap\C_{\omega}}(u)=r^{N+2}\int_{\widetilde\Omega_{r}}\left(\frac{|x|^{2}}{2N^{2}}-\frac{1-|x|^{2}}{2N}\right)\,dx=-r^{N+2}\left[\frac{\sigma_{N}}{4N(N+2)}+o(1)\right]\quad\text{as }r\to 0,
$$
with $\widetilde\Omega_{r}$ as in \eqref{omegatilder}. Then \eqref{st3} follows from a direct computation, using the fact the the limit domain of $\widetilde\Omega_{r}$ tends to the half unit ball $\{x\in\R^{N}\,;~x_{N}>0\}$, as $r\to 0$. 

Finally, let us complete the proof of \eqref{eqtesi:stimaGammasect5-B}. By definition of $\mathcal{O}_{c}(\C_{\omega})$ and of $\mathcal{E}(\cdot;\C_{\omega})$ and by \eqref{st1}, \eqref{st2}, and \eqref{st3}, we have that
\begin{equation}\label{st3b}
\mathcal{O}_{c}(\C_{\omega})\le \mathcal{E}(B_{r(c)}(x_{0})\cap\C_{\omega};\C_{\omega})\le J_{B_{r(c)}(x_{0})\cap\C_{\omega}}(u)=-c^{1+\frac2N}\left[\frac1{2N(N+2)}\left(\frac2{\sigma_{N}}\right)^{\frac2N}+o(1)\right]~~\text{as }c\to 0^+.
\end{equation}
By \eqref{eq:boundlemmasect5} and \eqref{eq4:lem3sect5} we obtain
$$
\left(\frac{c}{\H^{N-1}(\Gamma)}\right)^{2}\ge\frac{2|\mathcal{O}_{c}(\C_{\omega})|}{c}
$$
which, together with \eqref{st3b}, after some algebraic computation, yields \eqref{eqtesi:stimaGammasect5-B}. 
\end{proof}

\begin{remark}\label{rem:extraassumptions}
We point out that the extra assumptions on the regularity of $\omega$, $\Omega^*$ and on its corresponding torsion function $u_{\Omega^*}$ play a role in the proofs of Proposition \ref{P:bdr-estimate-1} and Proposition \ref{P:bdr-estimate-2} because we use standard elliptic regularity theorem and the divergence theorem to obtain \eqref{eq4:lem3sect5}. This means that, if $\omega$ is merely a bounded Lipschitz domain of $\R^{N-1}$ and if we know that \eqref{eq4:lem3sect5} holds true, then we can rule out all the extra assumptions on $\Omega^*$, $u_{\Omega^*}$, and the estimates of Proposition \ref{P:bdr-estimate-1} and Proposition \ref{P:bdr-estimate-2} holds with minor modifications in the proofs.
\end{remark}

In view of Propositions \ref{P:bdr-estimate-1} and \ref{P:bdr-estimate-2}, we can deduce further properties of the minimizers of $\E(\cdot;\C_{\omega})$ in $\A_{\omega,c}$. In particular, for $c>0$ small enough, they cannot be bounded cylinders of the form $\omega\times\left]a,b\right[$ with volume $c$ and, at the same time, their closure intersects the boundary of the container $\C_{\omega}$. Moreover, this last property holds true also when $c$ is sufficiently large. More precisely, we have:

\begin{corollary}\label{cor1:sect6}
Let $\omega\subset\R^{N-1}$ be a Lipschitz bounded domain. Then:
\begin{itemize}
\item[\emph{(i)}] 
there exists $c_0>0$ such that for every $0<c<c_0$ the bounded cylinder $\Omega_{\omega,h}:=\omega\times\left]-\frac{h}{2},\frac{h}{2}\right[$, with $h=\frac{c}{\H^{N-1}(\omega)}$ (or any translation of it in the $x_{N}$-direction), cannot be a minimizer of $\E(\cdot;\C_{\omega})$ in $\A_{\omega,c}$. 
\item[\emph{(ii)}] There exists $c_1>0$ such that, if $0<c<c_1$ and $\Omega^{*}$ is a minimizer of $\E(\cdot;\C_{\omega})$ in $\A_{\omega,c}$, then $\mathcal{H}^{N-1}(\partial\Omega^{*}\cap\partial\C_{\omega})>0$.

\item[\emph{(iii)}] 
If $\omega$ is of class $C^{2,\alpha}$ and
\begin{equation}\label{c>}
c>\left[\frac{2\sqrt{3}\,\H^{N-1}(\omega)}{N\sigma_{N}^{\frac1N}}\right]^{\frac N{N-1}}
\end{equation}
then for any minimizer $\Omega^{*}$ of $\E(\cdot;\C_{\omega})$ in $\A_{\omega,c}$, satisfying the hypotheses of Proposition \ref{P:bdr-estimate-2}, with a torsion function $u_{\Omega^{*}}\in W^{1,\infty}(\Omega^*)$, it holds that $\mathcal{H}^{N-1}(\partial\Omega^{*}\cap\partial\C_{\omega})>0$.
\end{itemize}
\end{corollary}

\begin{proof} (i) Assume, by contradiction that the thesis is false, then there exists a sequence $c_{n}\to 0^{+}$ such that for every $n$ the bounded cylinder $\C_{\omega,h_{n}}:=\omega\times\left]-\frac{h_{n}}{2},\frac{h_{n}}{2}\right[$ with $h_{n}=\frac{c_{n}}{\H^{N-1}(\omega)}$ is a minimizer of $\E(\cdot;\C_{\omega})$ in $\A_{\omega,c_{n}}$. We observe that since $\C_{\omega,h_{n}}$ is a Lipschitz domain, and  its torsion function (given by \eqref{eq:torsionfunctbcylinder}), restricted to $\C_{\omega,h_{n}}$, is smooth up to the boundary, then, using the standard divergence theorem and by elementary computations we infer that \eqref{eq4:lem3sect5} holds true, with $C_0=\left(\frac{h_n}{2}\right)^2$, $c=c_n$. Then, from Proposition \ref{P:bdr-estimate-2}, taking into account Remark \ref{rem:extraassumptions}, we readily obtain a contradiction, because $\H^{N-1}(\Gamma_{\C_{\omega,h_{n}}})=2\mathcal{H}^{N-1}(\omega)$, whereas $c_{n}\to 0^+$.

For (ii), we argue again by contradiction, assuming that there exists a sequence $c_{n}\to 0^{+}$ such that for every $n$ there exists a minimizer $\Omega^{*}_{n}$ of $\E(\cdot;\C_{\omega})$ in $\A_{\omega,c_{n}}$ with $\mathcal{H}^{N-1}(\partial\Omega^{*}_n\cap\partial\C_{\omega})=0$. Then $H_0^1(\Omega^*_{n}; \C_\omega)=H_0^1(\Omega^*_{n})$ (see \cite[Remark 4.3]{BV}). In particular 
$$
\E(\Omega^*_{n}; \C_\omega)=\E(\Omega^*_{n}; \R^N),
$$ 
where $\E(\Omega^*_{n}; \R^N)$ denotes the ``free'' torsional energy of $\Omega^*_{n}$, namely $\E(\Omega^*_{n}; \R^N)$ is the minimizer in $H_0^1(\Omega^*_{n})$ of the functional $J(v)=\frac{1}{2}\int_{\R^2} \big(|\nabla v|^2 - v\big) \ dx$. Then, arguing as in \cite[(6.32)-(6.33)]{IPW}, considering the Schwartz symmetrization of the energy function $u_{\Omega^*_{n}}$ and thanks to the P\'olya-Szeg\"o inequality, we infer that $\Omega^*_{n}=B_{r_{n}}(x_{n})\subset\C_{\omega}$ for some $x_{n}\in \R^{N}$ and $r_{n}>0$ such that $\sigma_{N}r_{n}^{N}=|\Omega^{*}_{n}|=c_{n}$ ($\sigma_{N}$ is the measure of the unit ball in $\R^{N}$). Moreover $u_{\Omega^{*}_{n}}(x)=\frac{r_{n}^{2}-|x-x_{n}|^{2}}{2N}$, for $x\in\Omega^{*}_{n}$, and 
\begin{equation}\label{bdr1}
\E(\Omega^*_{n}; \C_\omega)=J_{\Omega_{n}^{*}}(u_{\Omega^{*}_{n}})=-\frac{c_{n}^{1+\frac2N}}{2N(N+2)\sigma_{N}^{\frac 2N}}.
\end{equation}
The same computations made in the proof of Proposition \ref{P:bdr-estimate-2} (see \eqref{st3b}) lead to estimate
\begin{equation}\label{bdr2}
\O_{c_{n}}(\C_{\omega})\le -c_{n}^{1+\frac2N}\left[\frac1{2N(N+2)}\left(\frac2{\sigma_{N}}\right)^{\frac2N}+o(1)\right]~~\text{as }n\to\infty.
\end{equation}
Since $\E(\Omega^*_{n}; \C_\omega)=\O_{c_{n}}(\C_{\omega})$, from \eqref{bdr1}--\eqref{bdr2} a contradiction follows.
\medskip

\noindent
Let us prove (iii). Let $\Omega^{*}$ be a minimizer of $\E(\cdot;\C_{\omega})$ in $\A_{\omega,c}$ satisfying the hypotheses of Proposition \ref{P:bdr-estimate-2} with a torsion function $u_{\Omega^{*}}\in W^{1,\infty}(\Omega^*)$. If $\mathcal{H}^{N-1}(\partial\Omega^{*}\cap\partial\C_{\omega})=0$ then $\mathcal{H}^{N-1}(\partial\Omega^{*})=\mathcal{H}^{N-1}(\Gamma_{\Omega^{*}})$ and, by Corollary \ref{cor1:sect5}, Proposition \ref{P:bdr-estimate-2} and the isoperimetric inequality (see e.g. \cite[Chap. 14]{MA}), 
$$
c\le\left[\frac{\mathcal{H}^{N-1}(\partial\Omega^{*})}{N\sigma_{N}^{\frac1N}}\right]^{\frac N{N-1}}\le
\left[\frac{2\sqrt{3}\,\H^{N-1}(\omega)}{N\sigma_{N}^{\frac1N}}\right]^{\frac N{N-1}}
$$
in contradiction with \eqref{c>}. The proof is complete. 
\end{proof}

Finally, in the 2-dimensional case, we can provide a more precise description of the minimizers of $\E(\cdot; \C_\omega)$, especially those which are symmetric with respect to the $x_{1}$-axis and convex in the $x_2$-direction (they always exist in view of Theorem \ref{mainteo1} and Lemma \ref{lem1:sect2}). Indeed, we have:

\begin{proposition}\label{prop:charelboundaryN2}
Let $\omega=\left]0,a\right[$ with $a>0$, $\C_{\omega}=\omega\times\R$, and let $c>0$. Let $\Omega^*$ be a minimizer of $\E(\cdot;\C_\omega)$ in $\A_{\omega,c}$. Then $\mathcal{H}^{1}(\partial\Omega^{*}\cap\partial\C_{\omega})>0$. Moreover, if: \vspace{3pt}
\begin{itemize}
\item[\emph{(a)}]
$\Omega^{*}$ is symmetric with respect to the $x_{1}$-axis and convex in the $x_2$-direction,\vspace{6pt}
\item[\emph{(b)}]
$u_{\Omega^{*}}\in W^{1,\infty}(\Omega^{*})$,
\vspace{3pt}
\end{itemize}
then only one of the following possibilities holds true: 
\begin{itemize}
\item[\emph{(i)}] $\Omega^{*}$ is a half-disk centered at a point of $\partial\C_{\omega}\cap\{x_{2}=0\}$ and with radius $\sqrt\frac{2c}\pi<a$.\vspace{3pt}
\item[\emph{(ii)}] $\overline{\Gamma_{\Omega^{*}}}=\overline{\Gamma_{+}}\cup\overline{\Gamma_{-}}$ where $\overline{\Gamma_{+}}$ is a curve in the upper semi-cylinder $[0,a]\times\left]0,+\infty\right[$ joining two points $(0,\ell_{0})$ and $(a,\ell_{1})$ with $\ell_{0},\ell_{1}>0$ and $\overline{\Gamma_{-}}$ is the symmetric curve of $\overline{\Gamma_{+}}$ with respect to the $x_{1}$-axis. 
\end{itemize}
Finally, if $c>\frac{3a^2}{\pi}$, then only \emph{(ii)} can happen.
\end{proposition}

\begin{proof}
Let $\Omega^*\subset\C_{\omega}$ be a minimizer of $\E(\cdot;\C_\omega)$ in $\A_{\omega,c}$.
\medskip

\noindent
\textbf{Claim I:} $\mathcal{H}^{1}(\partial\Omega^*\cap\partial\C_\omega)>0$.
\medskip

\noindent
Otherwise $H_0^1(\Omega^*; \C_\omega)=H_0^1(\Omega^*)$ (see \cite[Remark 4.3]{BV}). In particular 
$$
\E(\Omega^*; \C_\omega)=\E(\Omega^*; \R^2):=\inf_{v\in H_0^1(\Omega^*)}\int_{\R^2} \left(\frac12|\nabla v|^2 - v\right) \, dx,
$$ 
and as seen in the proof of Corollary \ref{cor1:sect6}-(ii), we infer that $\Omega^*$ must be a disk of some radius $r\in\left]0,\frac a2\right[$ and $u_{\Omega^{*}}(x)=\frac{r^{2}-|x-x_{0}|^{2}}4$, where $x_{0}\in\R^{2}$ is the center of $\Omega^{*}$. By Corollary \ref{cor1:sect5}, $|\Omega^{*}|=c$. Hence $r=\sqrt\frac{c}\pi$ and $\E(\Omega^{*};\C_{\omega})=J_{\Omega^{*}}(u_{\Omega^{*}})=-\frac{c^{2}}{16\pi}$. Moreover, as $\Omega^{*}\subset\C_{\omega}$, it must be $r\le\frac a2$. Now, let us consider the half-disk $\Omega_{0}:=\{x\in\C_{\omega}\,;~|x|<\sqrt{2}r\}$ and its corresponding energy function $u_{\Omega_{0}}(x)=\frac{2r^{2}-|x|^{2}}4$, for $x\in\Omega_0$, then $\Omega_{0}\in\mathcal{A}_{\omega,c}$ and $\E(\Omega_{0};\C_{\omega})=J_{\Omega_{0}}(u_{\Omega_{0}})=-\frac{c^{2}}{8\pi}<\E(\Omega^{*};\C_{\omega})$, against the fact that $\E(\Omega^{*};\C_{\omega})=\O_{c}(\C_{\omega})$. The proof of \textit{Claim I} is complete.
\medskip

\noindent
Assume now that $\Omega^*$ is a minimizer of $\E(\cdot;\C_\omega)$ in $\A_{\omega,c}$ satisfying (a) and its corresponding energy function verifies (b). Since $\partial\Omega^*\cap\C_{\omega}=(\partial\Omega\cap\{x_{1}=0\})\cup(\partial\Omega\cap\{x_{1}=a\})$, by \textit{Claim I}, we can assume without loss of generality that 
\begin{equation}
\label{bordo-sinistro-pos}
\mathcal{H}^{1}(\partial\Omega^{*}\cap\{x_{1}=0\})>0. 
\end{equation}
Let $\omega^{*}$ be the orthogonal projection of $\Omega^{*}$ on the $x_{1}$-axis. Since $\Omega^{*}$ is connected (see Theorem \ref{connected}), $\omega^{*}$ is an open interval and there exists $a^{*}\in\left]0,a\right]$ such that 
$$
\omega^{*}=\left]0,a^{*}\right[.
$$
\textbf{Claim II:} $\overline{\Gamma_{\Omega^{*}}}\cap\{x_{1}=0\}=\{(0,\ell_{0}),(0,-\ell_{0})\}$ for some $\ell_{0}>0$ and $\overline{\Gamma_{\Omega^{*}}}\cap\{x_{1}=a\}$ either is empty or is the set $\{(0,\ell_{1}),(0,-\ell_{1})\}$ for some $\ell_{1}\ge 0$.
\medskip

\noindent
Since $\Omega^{*}$ satisfies (a), \eqref{bordo-sinistro-pos}, and $\Gamma_{\Omega^*}$ is smooth (by Theorem \ref{teoreg}), there exist $\varepsilon_{0}>0$ and $\gamma\colon\left]0,1\right]\to\R^{2}$ of class $C^{\infty}$, $\gamma(t)=(\gamma_1(t),\gamma_2(t))$, such that 
\begin{gather*}\Gamma_{\Omega^{*}}\cap\{x\in\R^{2}\,;~x_{1}\in\left]0,\varepsilon_{0}\right]\,,~x_{2}>0\}=\gamma(\left]0,1\right]),\\\gamma_{1}(t)\to 0~~\text{as }t\to 0^+~~\text{and}~~~~\gamma_{1}'(t)\ge 0~\forall t\in\left]0,1\right].
\end{gather*}
For every $\varepsilon\in\left]0,\varepsilon_{0}\right[$ let $\Omega_{\varepsilon}:={\Omega^{*}}\cap\{x\in\R^{2}\,;~x_{1}\in\left]\varepsilon,\varepsilon_{0}\right[\}$. Since $u_{\Omega^{*}}$ solves 
$$
\begin{cases}
-\Delta u_{\Omega^{*}}=1&\text{in~~$\Omega^{*}$}\\ u_{\Omega^{*}}=0&\text{on~~$\Gamma_{\Omega^{*}}$}\\
\frac{\partial u_{\Omega^{*}}}{\partial\nu}=-\sqrt{C_{0}}&\text{on~~$\Gamma_{\Omega^{*}}$}
\end{cases}
$$
for some $C_{0}>0$, by the divergence Theorem we have
$$
|\Omega_{\varepsilon}|=-\int_{\Omega_{\varepsilon}}\Delta u_{\Omega^{*}}\,dx
%=-\int_{\partial\Omega_{\varepsilon}}\frac{\partial u_{\Omega^{*}}}{\partial\nu}\,d\sigma
=\int_{\Omega^{*}\cap\{x_{1}=\varepsilon\}}\frac{\partial u_{\Omega^{*}}}{\partial x_{1}}\,ds-\int_{\Omega^{*}\cap\{x_{1}=\varepsilon_{0}\}}\frac{\partial u_{\Omega^{*}}}{\partial x_{1}}\,ds+2\sqrt{C_{0}}\mathcal{H}^{1}(\Gamma_{\Omega^{*}}\cap\{\varepsilon<x_{1}<\varepsilon_{0}\})
$$
and then, using (b), we infer that
$$
\mathcal{H}^{1}(\Gamma_{\Omega^{*}}\cap\{\varepsilon<x_{1}<\varepsilon_{0}\})\le \frac{|\Omega_{\varepsilon}|+\|\nabla u_{\Omega^{*}}\|_{\infty}\left[\mathcal{H}^{1}(\Omega^{*}\cap\{x_{1}=\varepsilon\})+\mathcal{H}^{1}(\Omega^{*}\cap\{x_{1}=\varepsilon_{0}\})\right]}{2\sqrt{C_{0}}}\le C,
$$
with $C<+\infty$ independent of $\varepsilon$, because $\Omega^{*}$ is bounded (see Proposition \ref{prop1:sect5}). Hence $$\mathcal{H}^{1}(\Gamma_{\Omega^{*}}\cap\{0<x_{1}<\varepsilon_{0}\})<+\infty.$$ This implies that there exists $\lim_{t\to 0^+}\gamma(t)=(0,\ell_{0})\in\R^{2}$ for some $\ell_{0}\ge 0$. Indeed, if not, then it holds that $\ell_{-}:=\liminf_{t\to 0^+}\gamma_{2}(t)<\limsup_{t\to 0^+}\gamma_{2}(t)=:\ell_{+}$. In this case, we have that $\Gamma_{\Omega^{*}}\cap\{\varepsilon<x_{1}<\varepsilon_{0}\}$ is a smooth curve which crosses infinitely many times the strip $\{\ell_{-}+\delta\le x_{2}\le \ell_{+}-\delta\}$ ($\delta>0$ small enough) and then $\mathcal{H}^{1}(\Gamma_{\Omega^{*}}\cap\{0<x_{1}<\varepsilon_{0}\})=+\infty$, a contradiction. The fact that $\ell_{0}>0$ follows from \eqref{bordo-sinistro-pos}. A similar argument holds true for $\overline{\Gamma_{\Omega^{*}}}\cap\{x_{1}=a\})$ if it is nonempty. Thus \textit{Claim II} is proved. 
\medskip

\noindent
\textbf{Claim III:} if $a^{*}<a$, then $\Omega^{*}$ is the right half-disk centered at $0$ and with radius $a^{*}$. Hence, case (i) occurs. 
\medskip

\noindent
We reflect the domain $\Omega^{*}$ with respect to the $x_{2}$-axis. More precisely, we introduce the sets 
\begin{equation}\label{doubling}
\Omega^{*}_{-}:=\{(x_{1},x_{2})\in\R^{2}\,;~(-x_{1},x_{2})\in\Omega\}\,,~~\Omega:=\mathrm{int}(\overline{\Omega^{*}}\cup\overline{\Omega^{*}_{-}}),
\end{equation}
and the function $u\colon\overline\Omega\to\R$ defined by $u(x_{1},x_{2})=:u_{\Omega^{*}}(|x_{1}|,x_{2})$ for every $(x_{1},x_{2})\in\overline\Omega$. We point out that, since $\ell_{0}>0$, the set $\Omega$ is a bounded domain contained in $\left]-a,a\right[\times\R$ and $u$ is a (weak) positive solution to the overdetermined problem
\beq\label{eq:overdetSect6P}
\begin{cases}
-\Delta u = 1 & \text{in $\Omega$}\\
 u = 0 & \text{on $\partial\Omega$}\\
\frac{\partial u}{\partial \nu} = -\sqrt{C_0} & \text{on $\partial\Omega$}
\end{cases}
\eeq
for some positive constant $C_0$. Hence by the Weinberger version of the Serrin symmetry theorem without assumptions on the smoothness of $\partial\Omega$ (see \cite[Theorem 1]{GL} and \cite{SE, WE}), $\Omega$ is a disk centered at some point lying on the line $x_{1}=0$ and the claim is proved. 
\medskip

\noindent
\textbf{Claim IV:} if $a^{*}=a$, then case (ii) occurs. 
\medskip

\noindent
Indeed, in this case $\partial\Omega^{*}$ intersects $\partial\C_{\omega}$ both on the line $x_{1}=0$ and on the line $x_{1}=a$. Then, by (a) 
$$
\Gamma_{\Omega^{*}}=\Gamma_{+}\cup\Gamma_{-},
$$ 
where $\Gamma_{+}$ is a smooth curve in $\left]0,a\right[\times\left]0,+\infty\right[$, and $\Gamma_{-}$ is the symmetric curve to $\Gamma_{+}$ with respect to the line $x_{2}=0$, namely $\Gamma_{-}:=\{(x_{1},x_{2})\in\R^{2}\,;~(-x_{1},x_{2})\in\Gamma_{+}\}$. Moreover, by \textit{Claim II}, there exists $\ell_{1}\ge 0$ such that  
$$
\overline{\Gamma_{+}}=\Gamma_{+}\cup\{(0,\ell_{0}),(a,\ell_{1})\}\text{~~and~~}\partial\Omega^{*}=\Gamma_{\Omega^{*}}\cup S_{0}\cup S_{1}
$$
where $S_{0}$ is the vertical segment joining $(0,\ell_{0})$ and $(0,-\ell_{0})$ and $S_{1}$ is the vertical segment joining $(a,\ell_{1})$ and $(a,-\ell_{1})$. Let us exclude that $\ell_{1}=0$. Indeed, if this occurs, then, arguing as in the proof of \textit{Claim III}, we conclude that $\Omega^{*}$ is the right half-disk centered at $0$ with radius $a$. Then, repeating the computations already made in \textit{Claim I}, we have
$$
\O_{c}(\C_{\omega})=\E(\Omega^{*};\C_{\omega})=-\frac{c^{2}}{8\pi}=-\frac{\pi a^{2}}{32}
$$
because 
\begin{equation}
\label{c=}
c=|\Omega^{*}|=\frac{\pi a^{2}}{2}.
\end{equation} 
Now, let us consider the open rectangle $\Omega_{\omega,h}=\omega\times\left]-\frac h2,\frac h2\right[$, with $h=\frac ca$. Then also $\Omega_{\omega,h}\in\A_{\omega,c}$, and its torsion function is $u_{\Omega_{\omega,h}}(x_{1},x_{2})=\frac12\big(\frac{h^{2}}4-x_{2}^{2}\big)$, for $(x_1,x_2)\in\Omega_{\omega,h}$, and 
$$
\E(\Omega_{\omega,h};\C_{\omega})=J_{\Omega_{\omega,h}}(u_{\Omega_{\omega,h}})=-\frac{c^{3}}{24 a^{2}}.
$$ 
Then, by \eqref{c=}, $\E(\Omega_{\omega,h};\C_{\omega})<\O_{c}(\C_{\omega})$, a contradiction. Thus $\ell_{1}>0$ and \textit{Claim IV} is completely proved.
\medskip

\noindent
\textbf{Claim V:} if $c>\frac{3a^2}{\pi}$, then only (ii) can happen.
\medskip

\noindent
Let us consider the right half-disk $\Omega_{0}$ centered at $0$ and with radius $r:=\min\big\{a,\sqrt\frac{2c}\pi\big\}$. By the computations made in \textit{Claim I}, we have that $\Omega_{0}\in\A_{\omega,c}$ and $\E(\Omega_{0};\C_{\omega})\ge -\frac{c^{2}}{8\pi}$. Now, let us consider the open rectangle $\Omega_{\omega,h}=\left]0,a\right[\times\left]-\frac h2,\frac h2\right[$ with $h=\frac ca$. Then, as already noted in the proof of \textit{Claim III}, also $\Omega_{\omega,h}\in\A_{\omega,c}$ and $\E(\Omega_{\omega,h};\C_{\omega})=J_{\Omega_{\omega,h}}(u_{\Omega_{\omega,h}})=-\frac{c^{3}}{24 a^{2}}$. If $c>\frac{3a^2}{\pi}$, then $\E(\Omega_{0};\C_{\omega})>\E(\Omega_{\omega,h};\C_{\omega})\ge\mathcal{O}_{c}(\C_{\omega})$. Then (i) cannot occur. Hence, also \textit{Claim V} is proved.
\end{proof}

\begin{remark}
In the setting of Proposition \ref{prop:charelboundaryN2}, let $\Omega^*$ be a minimizer of $\E(\cdot;\C_\omega)$ in $\A_{\omega,c}$ satisfying (a) and (b), and let $\Gamma_{+}$ be the upper component of the relative boundary $\Gamma_{\Omega^{*}}$, defined as in the previous proof. Assuming that case (ii) holds true, $\Gamma_{+}$ admits a parametrization by arc-length $\gamma\in C^\infty(\left]0,1\right[; \R^2)$ which admits a continuous extension on $\left[0,1\right]$, with $\gamma(0)=(0,\ell_{0})$ and $\gamma(1)=(0,\ell_{1})$. Assume that there exist $\lim_{t\to 0}\gamma'(t)=\tau_{0}$ and $\lim_{t\to 1}\gamma'(t)=\tau_{1}$. Then $\tau_{0}$ and $\tau_{1}$ are two versors in $\R^{2}$ and we can define the corresponding contact angles 
$$
\theta_{0}:=\arccos\langle\tau_{0}, e_{2}\rangle~~\text{and}~~\theta_{1}=\arccos\langle\tau_{1}, e_{2}\rangle
$$ 
where $e_{2}=(0,1)$. By symmetry, the same can be done also for the lower component $\Gamma_{-}$ of $\Gamma_{\Omega^{*}}$. Hence, $\Omega^*$ is a boundary edge domain of $\C_\omega$, in the sense of \cite[Definition 2.1]{LS}. Now, since $\mathrm{Reg}(\Gamma_{\Omega^{*}})=\Gamma_{\Omega^{*}}$ (because $N=2$), then by Proposition \ref{lem2:sect5} we get that $\frac{\partial u_{\Omega^*}}{\partial \nu}$ is constant on $\Gamma_{\Omega^*}$, and from the same proof of \cite[Corollary 2.6]{LS} we deduce that $\theta_{0}=\theta_{1}=\frac{\pi}{2}$, i.e., $\overline{\Gamma_{\Omega^{*}}}$ meets $\partial\C_{\omega}$ orthogonally. This implies that we can construct a domain $\Omega_{1}$ by repeating the set $\Omega$ defined in \eqref{doubling} in a periodic way in the $x_{1}$-direction and a corresponding function $u_{1}\colon\overline{\Omega_{1}}\to\R$ which solves an overdetermined problem like \eqref{eq:overdetSect6P} in $\Omega_{1}$. This is related to a result in \cite{FMW} (see also \cite{RSW} and \cite{PRS}). 
We also note that, in view of the results obtained in recent work \cite{FTV}, we believe that the additional assumptions on $\Gamma_{+}$ and $\Gamma_-$, as well as the hypothesis $u_{\Omega^*}\in W^{1,\infty}(\Omega^*)$, can be omitted. However, we do not address this aspect which is beyond the scope of the present paper.
\end{remark}

\section{The minimization problem in the half-cylinder}\label{S:half} 

 Let $\omega$ be a bounded domain of $\R^{N-1}$ with Lipschitz boundary and let $\mathcal{C}^{+}_\omega:=\omega\times\left]0,\infty\right[$ be the half-cylinder generated by $\omega$.

For any quasi-open set $\Omega\subset \C^{+}_\omega$ we consider the Sobolev space:
$$H^1_0(\Omega; \C^{+}_\omega):=\left\{u \in H^1(\C^{+}_\omega)\,;~ u=0 \ \ \hbox{q.e. on} \ \C^{+}_\omega\setminus\Omega  \right\},$$
where q.e. means quasi-everywhere, i.e. up to sets of zero capacity. 

The space $H^1_0(\Omega; \C^{+}_\omega)$ is the natural functional space to study the torsion problem
\begin{equation}\label{eq:mixbvprobquasiopen+}
\begin{cases}
-\Delta u = 1 & \text{in $\Omega$,}\\
 u = 0 & \text{on $\partial\Omega \cap \C^{+}_\omega$,}\\
 \frac{\partial u}{\partial \nu} = 0 & \text{on $\partial \C^{+}_\omega$.}
\end{cases}
\end{equation}
A (weak) solution of \eqref{eq:mixbvprobquasiopen+} is a critical point of the functional $J^{+}_\Omega\colon H^1_0(\Omega; \C^{+}_\omega)\to \R$ defined by
$$
J^{+}_\Omega(u)=\frac{1}{2}\int_{\C^{+}_\omega} |\nabla u|^2 \ dx - \int_{\C^{+}_\omega} u \ dx.
$$
A variant of Lemma \ref{energy-function} with $\C^{+}_{\Omega}$ and $J^{+}_{\Omega}$ instead of $\C_{\Omega}$ and $J_{\Omega}$ still holds true, with no change in the proof. We always denote by $u_{\Omega}$ the corresponding energy function, characterized as minimum point of $J^{+}_\Omega$ in $H^1_0(\Omega; \C^{+}_\omega)$ and we observe that such $u_{\Omega}$ always possesses the properties stated in Proposition \ref{prop1:sect2}. \vspace{1.8pt}

Next, we are going to study the problem of minimizing the functional 
$$
\E(\Omega;\C^{+}_\omega):=J^{+}_{\Omega}(u_{\Omega})
$$ 
among quasi-open sets of uniformly bounded measure. To this aim for any $c>0$ we denote by $\A^{+}_{\omega,c}$  the class of quasi-open sets in $\C_\omega^+$ of measure less or equal than $c$, namely
$$
\A^{+}_{\omega,c}:= \{\Omega\subset \C^{+}_\omega\,;~ \Omega \ \hbox{quasi-open} \ \hbox{and} \ |\Omega|\leq c\},
$$
and define
$$
\mathcal{O}_c(\C_\omega^+):=\inf\{\E(\Omega;\C_\omega^+)\,;~ \Omega\in \A^+_{\omega,c}\}.
$$
In the next theorem we give a relation between $\mathcal{O}_c(\C_\omega^+)$ and $\mathcal{O}_{2c}(\C_\omega)$ and show that $\mathcal{O}_c(\C_\omega^+)$ is achieved.
 
\begin{theorem}\label{teo1:sect7}
Let $c>0$ and $\omega\subset\R^{N-1}$ be a Lipschitz bounded domain. Then 
\beq\label{eqtesi:prop1sect6}
\mathcal{O}_{c}(\mathcal{C}_\omega^+)=\frac{1}{2}\mathcal{O}_{2c}(\mathcal{C}_\omega),
\eeq
and $\mathcal{O}_c(\C_\omega^+)$ is attained.
\end{theorem}
\begin{proof}
We start by proving that
\beq\label{eq1:sect7}
\mathcal{O}_{c}(\mathcal{C}_\omega^+) \leq \frac{1}{2}\mathcal{O}_{2c}(\mathcal{C}_\omega).
\eeq
By Theorem \ref{teo:mainteoexistmin}, Proposition \ref{prop1:sect5}, Proposition \ref{prop2:sect5} and Corollary \ref{cor1:sect5} we know that $\O_{2c}(\C_{\omega})$ is attained by a  bounded open set $\Omega^*$ of $\C_\omega$ such that $|\Omega^*|=2c$. Without loss of generality, by Lemma \ref{lem1:sect2}, we can assume that $\Omega^*$ is convex in the $x_N$-direction and symmetric with respect to the hyperplane $\{x_N=0\}$. 

Now, setting 
\beq\label{eq2:sect7}
\Omega^*_+:=\Omega^*\cap\{x_N>0\} 
\eeq
we have that $\Omega^*_+ \subset \C_\omega^+$ is an open set and $|\Omega^*_+|= c$. Let $u_{\Omega^*}\in H_0^1(\Omega^*; \C_\omega)$ be the energy function of $\Omega^*$ which is even with respect to $x_N$. Then the restriction $w:={u_{\Omega^*}}_{\big\vert_{\Omega^*_+}}$ belongs to $H_0^1(\Omega^*_+; \C_\omega^+)$ and 
$$J_{\Omega^*_+}^+(w)=\frac{1}{2}\int_{\C_\omega^+} |\nabla w|^2 \ dx - \int_{\C_\omega^+} w \ dx=\frac{1}{2}\left[\frac{1}{2}\int_{\C_\omega} |\nabla u_{\Omega^*}|^2 \ dx- \int_{\C_\omega} u_{\Omega^*} \, dx\right]=\frac{1}{2}J_{\Omega^*}(u_{\Omega^*}).$$
Hence, by definition and by construction we deduce that
$$
\E(\Omega^*_+; \mathcal{C}_\omega^+)\leq \frac{1}{2}\E(\Omega^*; \mathcal{C}_\omega)=\frac{1}{2}\mathcal{O}_{2c}(\mathcal{C}_\omega),
$$
and hence \eqref{eq1:sect7} holds. Let us prove the opposite inequality. 
%%
%% 
%% Now we prove \eqref{eqtesi:prop1sect6}. Assume by contradiction that \eqref{eqtesi:prop1sect6} is false. By \eqref{eq1:sect7} the only possibility is
%% 
%% \beq\label{eq2:teo1sect7}
%% \mathcal{O}_{c}(\mathcal{C}_\omega^+)<\frac{1}{2}\mathcal{O}_{2c}(\mathcal{C}_\omega).
%% \eeq
Arguing as in the proof of Lemma \ref{lem1:sect3}-(i), we find a sequence of open sets $(A_n^+)_n \subset\C_\omega^+$ and a sequence of positive real numbers $(\delta_n)_n$ such that:
\begin{itemize}
\item[(i)] $|A_n^+|=c+\delta_n$, for all $n\in\N$;
\item[(ii)] $\delta_n\to 0^+$, as $n\to +\infty$;
\item[(iii)] $\E(A_n^+;\C_\omega^+)\to \mathcal{O}_{c}(\mathcal{C}_\omega^+)$, as $n\to +\infty$.
\end{itemize}
Let $u_{A_n^+}\in H_0^1(A_n^+; \C_\omega^+)$ be the energy function of $A_n^+$. Then $u_{A_n^+}$ admits an extension to $\C_\omega$, denoted $v_{n}$, which is even with respect to $x_{N}$ and belongs to $H^1(\C_\omega)$. 
Setting $$A_n:=\{v_n>0\},$$ we have that $A_n$ is a quasi-open set of $\C_\omega$ (see Sect. \ref{S:preliminaries}) and $|A_n|=2c+2\delta_n$, because $u_{A_n^+}>0$ on $A_{n}^{+}$. Moreover, by definition, it holds that $v_n\in H_0^1(A_n; \C_\omega)$ and thus we get that
\beq\label{eq3:teo1sect7}
\O_{2c+2\delta_n}(\C_\omega)\leq\E(A_n; \C_\omega)\leq J_{A_n}(v_n)= 2 J_{A_n^+}^+(u_{A_n^+})=2\E(A_n^+; \C_\omega^+).
\eeq
Taking the limit as $n\to+\infty$ in \eqref{eq3:teo1sect7}, from (ii), (iii) and Proposition \ref{prop3:sect5}, we obtain
$$\O_{2c}(\C_\omega)\leq 2 \O_c(\C_\omega^+)$$
which, together with \eqref{eq1:sect7}, yields 
%% contradicting \eqref{eq2:teo1sect7}. This contradiction proves 
\eqref{eqtesi:prop1sect6}. 
The equality \eqref{eqtesi:prop1sect6} also proves that $\O_c(\C_\omega^+)$ is achieved. 
Indeed, if $\Omega^*$ is a minimizer for $\E(\cdot;\C_\omega)$ in $\A_{\omega,2c}$ and is convex and symmetric in the $x_N$-direction (which is always possible in view of Lemma \ref{lem1:sect2}), then $\Omega^*_+$ defined in \eqref{eq2:sect7} is a minimizer for $\E(\cdot;\C^{+}_\omega)$ in $\A^{+}_{\omega,2c}$. The proof is complete.
\end{proof}

\section*{Acknowledgments}
\thanks{We wish to thank Prof. Giorgio Tortone for the useful discussions about the proof of Theorem \ref{teoreg}.}

\end{document}